\documentclass[a4paper,11pt,reqno]{amsart}   
\usepackage{etex}
\usepackage[OT2, T1]{fontenc}
\usepackage{url}
\usepackage{amsmath}
\usepackage{lmodern}
\usepackage{array}
\usepackage{graphicx}
\usepackage{amsfonts}
\usepackage{amssymb}

\usepackage{soul}
\usepackage{xargs}                      
\usepackage[pdftex,dvipsnames]{xcolor}  

\usepackage{exscale, relsize}
\usepackage{stackengine}
\usepackage{scalerel}
\usepackage{setspace}

\definecolor{imperialred}{RGB}{237, 41, 57}
\definecolor{royalblue}{RGB}{64, 106, 212}
\definecolor{link}{RGB}{11,0,128}
\definecolor{olivegreen}{RGB}{128, 128, 0}
\usepackage[colorlinks=true, citecolor=royalblue,linkcolor=olivegreen,urlcolor=link]{hyperref}
\usepackage{amsthm}
\usepackage{esint}

\usepackage{paralist}
\usepackage{colonequals}		
\usepackage{mathabx}		
\usepackage[left]{lineno}
\usepackage{amscd}
\usepackage[all,cmtip]{xy}	
\usepackage{sseq}			
\usepackage{verbatim}
\usepackage{parskip}
\usepackage{microtype}
\usepackage[in]{fullpage}
\usepackage{graphicx}
\usepackage{mathrsfs} 			
\usepackage{bm} 				
\usepackage{cleveref}
\usepackage{enumitem}
\usepackage{moreenum}			
\usepackage{stmaryrd}
\usepackage{blindtext}
\usepackage{tikz-cd}
\usepackage{tikz}
\usetikzlibrary{bending,matrix,arrows,decorations.pathmorphing,backgrounds,positioning,fit,petri,arrows.meta}
\usepackage{textcomp}
\usepackage{indentfirst}
\usetikzlibrary{arrows}
\usepackage{enumitem}
\tikzset{commutative diagrams/.cd,arrow style=tikz,diagrams={>=latex'}}

\usepackage{filecontents}
\usepackage[alphabetic,lite]{amsrefs} 	
\usepackage{stmaryrd}	
\usepackage[foot]{amsaddr}
\usepackage{subfiles}
\usepackage{ragged2e}

\DeclareSymbolFont{cyrletters}{OT2}{wncyr}{m}{n}
\DeclareMathSymbol{\Sha}{\mathalpha}{cyrletters}{"58}

\newcommand{\gA}{\alpha}
\newcommand{\gB}{\beta}
\newcommand{\gG}{\gamma}

\newcommand{\gK}{\kappa}

\newcommand{\GG}{\Gamma}

\newcommand{\GGL}{\Lambda}

\newcommand{\bA}{\mathbb{A}}

\newcommand{\bG}{\mathbb{G}}

\newcommand{\cA}{\mathcal{A}}

\newcommand{\cH}{\mathcal{H}}

\newcommand{\cL}{\mathcal{L}}

\newcommand{\cT}{\mathcal{T}}

\newcommand{\cV}{\mathcal{V}}

\newcommand{\cX}{\mathcal{X}}
\newcommand{\cY}{\mathcal{Y}}

\newcommand{\fc}{\mathfrak{c}}
\newcommand{\fd}{\mathfrak{d}}

\newcommand{\fg}{\mathfrak{g}}
\newcommand{\fh}{\mathfrak{h}}
\newcommand{\fm}{\mathfrak{m}}

\newcommand{\fp}{\mathfrak{p}}
\newcommand{\fq}{\mathfrak{q}}

\newcommand{\fP}{\mathfrak{P}}

\newcommand{\sD}{\mathscr{D}}
\newcommand{\sE}{\mathscr{E}}
\newcommand{\sF}{\mathscr{F}}
\newcommand{\sG}{\mathscr{G}}
\newcommand{\sH}{\mathscr{H}}
\newcommand{\sI}{\mathscr{I}}

\newcommand{\sM}{\mathscr{M}}

\newcommand{\sO}{\mathscr{O}}

\newcommand{\sT}{\mathscr{T}}

\newcommand{\sW}{\mathscr{W}}

\newcommand{\ra}{\rightarrow}    
\newcommand{\la}{\leftarrow}
\newcommand{\lra}{\longrightarrow}                 
\newcommand{\lla}{\longleftarrow}                  
\newcommand{\hra}{\hookrightarrow}                 
\newcommand{\surjects}{\twoheadrightarrow}         
\newcommand{\injects}{\hookrightarrow}             
\newcommand{\Ra}{\Rightarrow}                      
\newcommand{\longeq}{{=\joinrel=\joinrel=}}        

\newcommand{\isom}{\simeq}                         
\newcommand{\isab}{\overset{\sim}{\lla}}           
\newcommand{\isoto}{\overset{\sim}{\lra}}          
\newcommand{\ov}{\overline}                        
\newcommand{\un}{\underline}                       
\newcommand{\wt}{\widetilde}                       
\newcommand{\wh}{\widehat}                         

\makeatletter
\DeclareRobustCommand\widecheck[1]{{\mathpalette\@widecheck{#1}}}
\def\@widecheck#1#2{%
\setbox\z@\hbox{\m@th$#1#2$}%
\setbox\tw@\hbox{\m@th$#1%
\widehat{%
\vrule\@width\z@\@height\ht\z@
\vrule\@height\z@\@width\wd\z@}$}%
\dp\tw@-\ht\z@
\@tempdima\ht\z@ \advance\@tempdima2\ht\tw@ \divide\@tempdima\thr@@
\setbox\tw@\hbox{%
\raise\@tempdima\hbox{\scalebox{1}[-1]{\lower\@tempdima\box
\tw@}}}%
{\ooalign{\box\tw@ \cr \box\z@}}}
\makeatother

\newcommand{\hva}{\!\stackon[-8pt]{\;$V^{\raisebox{1pt}{$\mathsmaller{\mathsmaller{a}}$}}$}{\vstretch{1.5}{\hstretch{1.7}{\widehat{\phantom{}}}}}}
\newcommand{\hRa}{\!\stackon[-8pt]{\;$R^{\raisebox{1pt}{$\mathsmaller{\mathsmaller{a}}$}}$}{\vstretch{1.5}{\hstretch{1.7}{\widehat{\phantom{}}}}}}
\newcommand{\hka}{\!\stackon[-8pt]{\;$K^{\raisebox{1pt}{$\mathsmaller{\mathsmaller{a}}$}}$}{\vstretch{1.5}{\hstretch{1.7}{\widehat{\phantom{}}}}}}
\newcommand{\hvai}{\!\!\stackon[-8pt]{\;$V^{\raisebox{1pt}{$\mathsmaller{\mathsmaller{a}}$}}$}{\vstretch{1.5}{\hstretch{1.7}{\widehat{\phantom{}}}}} {[\f{1}{a}]}}
\newcommand{\vaih}{\!\!\stackon[-8pt]{\;$V[\f{1}{a}]^{\raisebox{1pt}{$\mathsmaller{\mathsmaller{a}}$}}$}{\vstretch{1.5}{\hstretch{3}{\widehat{\phantom{}}}}} }
\newcommand{\via}{V{[\f{1}{a}]}}

\newcommand{\nhka}{\wh{K}^a}
\newcommand{\nhva}{\wh{V}^a}

\newcommand{\pr}{^{\prime}}                        
\newcommand{\prpr}{^{\prime\prime}}                
\newcommand{\op}{^{\mathrm{op}}}                   
\DeclareMathOperator*{\dvee}{\vee\mkern-4mu\vee}   
\newcommand{\sep}{^{\mathrm{sep}}}

\providecommand{\p}[1]{\left(#1\right)}            
\providecommand{\hkh}[1]{\left\{#1\right\}}        
\providecommand{\fps}[1]{\llbracket#1\rrbracket}   
\providecommand{\lps}[1]{(\!(#1)\!)}               
\providecommand{\abs}[1]{\left\lvert#1\right\rvert}

\providecommand{\f}[2]{\frac{#1}{#2}}              

\newcommand{\et}{\mathrm{\acute{e}t}}              
\newcommand{\Zar}{\mathrm{Zar}}                    
\newcommand{\fpqc}{\mathrm{fpqc}}                  

\newcommand{\reg}{^{\mathrm{reg}}}                 
\newcommand{\sh}{\mathrm{sh}}

\newcommand{\etp}{\pi_1^{\mathrm{\acute{e}t}}}     
\newcommand{\sHom}{\mathscr{H}\! om}               

\newcommand{\mtg}{\underline{\mathrm{Tor}}(G)}     
\newcommand{\urad}{\mathrm{rad}^u}                 
\newcommand{\Par}{\un{\mathrm{Par}}}

\providecommand{\ucolon}{{\upshape:} }

\providecommand{\SPD}[2]{\cite[\href{https://stacks.math.columbia.edu/tag/#1}{#1}, \href{https://stacks.math.columbia.edu/tag/#2}{#2}]{SP}}
\providecommand{\SP}[1]{\cite[\href{https://stacks.math.columbia.edu/tag/#1}{#1}]{SP}}

\newcommand{\q}{\quad}
\newcommand{\qq}{\quad\quad}
\newcommand{\qqq}{\quad\quad\quad}

\newcommand{\qqqqq}{\quad\quad\quad\quad\quad}

\newcommand{\qqqqqqq}{\quad\quad\quad\quad\quad\quad\quad}
\newcommand{\qqqqqqqq}{\quad\quad\quad\quad\quad\quad\quad\quad}

\newcommand{\ce}{\colonequals}         

\renewcommand{\b}{\textbf}             
\newcommand{\x}{\text}                 
\newcommand{\tst}{\textstyle}          

\DeclareMathOperator{\coker}{coker}		
\DeclareMathOperator{\im}{Im}			
\DeclareMathOperator{\Spec}{Spec}		
\DeclareMathOperator{\rad}{rad}			
\DeclareMathOperator{\Hom}{Hom}			
\DeclareMathOperator{\Frac}{Frac}		
\DeclareMathOperator{\id}{id}			
\DeclareMathOperator{\Gal}{Gal}	
\DeclareMathOperator{\Sym}{Sym}			
\DeclareMathOperator{\Res}{Res}		
\DeclareMathOperator{\GL}{GL}		
\DeclareMathOperator{\End}{End}		
\DeclareMathOperator{\Aut}{Aut}		
\DeclareMathOperator{\Lie}{Lie}		
\DeclareMathOperator{\Isom}{Isom}		
\DeclareMathOperator{\Sets}{\textbf{Set}}

\newcommand{\ba}{\begin{aligned}}
\newcommand{\ea}{\end{aligned}}
\newcommand{\be}{\begin{equation}}
\newcommand{\ee}{\end{equation}}
\newcommand{\pf}{\begin{proof}}
\newcommand{\bpf}{\begin{proof}}
\newcommand{\epf}{\end{proof}}
\newcommand{\bthm}{\begin{thm}}
\newcommand{\ethm}{\end{thm}}

\newcommand{\bprop}{\begin{prop}}
\newcommand{\eprop}{\end{prop}}
\newcommand{\bcor}{\begin{cor}}
\newcommand{\ecor}{\end{cor}}

\newcommand{\brem}{\begin{rem}}
\newcommand{\erem}{\end{rem}}

\newcommand{\brems}{\begin{rems} \hfill \begin{enumerate}[label=\b{\thenumberingbase.},ref=\thenumberingbase]}

\newcommand{\erems}{\end{enumerate} \end{rems}}
\newcommand{\begs}{\begin{egs} \hfill \begin{enumerate}[label=\b{\thenumberingbase.},ref=\thenumberingbase]}

\newcommand{\eegs}{\end{enumerate} \end{egs}}
\DeclareMathOperator{\Sch}{\textbf{Sch}}        
\newcommand{\eremstweak}{\end{enumerate} \end{rems-tweak}}
\newcommand{\eremst}{\end{enumerate} \end{rems-tweak}}
\newcommand{\blem}{\begin{lemma}}
\newcommand{\elem}{\end{lemma}}

\newcommand{\bconj}{\begin{conj}}
\newcommand{\econj}{\end{conj}}

\newcommand{\bprob}{\begin{Problem}}
\newcommand{\eprob}{\end{Problem}}

\newcommand{\bpp}{\begin{pp}}
\newcommand{\epp}{\end{pp}}
\newcommand{\bq}{\begin{Q}}
\newcommand{\eq}{\end{Q}}
\newcommand{\benum}{\begin{enumerate}[label={{\upshape(\alph*)}}]}
\newcommand{\benuma}{\begin{enumerate}[label={{\upshape(\arabic*)}}]}
\newcommand{\benumr}{\begin{enumerate}[label={{\upshape(\roman*)}}]}
\newcommand{\eenum}{\end{enumerate}}
\newcommand{\bc}{}
\newcommand{\bd}{\begin{defn}}
\newcommand{\ed}{\end{defn}}

\newcommand{\beg}{\begin{eg}}
\newcommand{\eeg}{\end{eg}}

\newcommand{\bcl}{\begin{claim}}
\newcommand{\ecl}{\end{claim}}

\makeatletter
\renewcommand*\env@matrix[1][\arraystretch]{%
  \edef\arraystretch{#1}%
  \hskip -\arraycolsep
  \let\@ifnextchar\new@ifnextchar
  \array{*\c@MaxMatrixCols c}}
\makeatother

\usepackage{aliascnt}

\newaliascnt{numberingbase}{subsection}

\theoremstyle{plain}
\newtheorem{thm}[numberingbase]{Theorem}
\Crefname{thm}{Theorem}{Theorems}

\Crefname{rethm}{Theorem}{Theorem}
\newtheorem{prop}[numberingbase]{Proposition}
\Crefname{prop}{Proposition}{Propositions}
\newtheorem{Q}[numberingbase]{Question}
\Crefname{Q}{Question}{Questions}
\newtheorem{Problem}[subsection]{Problem}
\Crefname{Problem}{Problem}{Problems}
\newtheorem{conj}[numberingbase]{Conjecture}
\Crefname{conj}{Conjecture}{Conjectures}
\newtheorem{cor}[numberingbase]{Corollary}
\Crefname{cor}{Corollary}{Corollaries}
\newtheorem{lemma}[numberingbase]{Lemma}

\Crefname{subprop}{Proposition}{Propositions}

\Crefname{subcor}{Corollary}{Corollaries}

\Crefname{sublem}{Lemma}{Lemmas}

\theoremstyle{remark}
\newtheorem{claim}[equation]{Claim}
\Crefname{claim}{Claim}{Claims}

\Crefname{subrem}{Remark}{Remarks}

\theoremstyle{definition}
\newtheorem{defn}[numberingbase]{Definition}
\Crefname{defn}{Definition}{Definitions}

\Crefname{conv}{Convention}{Conventions}
\newtheorem{eg}[numberingbase]{Example}
\Crefname{eg}{Example}{Examples}
\newtheorem{rem}[numberingbase]{Remark}
\Crefname{rem}{Remark}{Remarks}
\newtheorem*{rems}{Remarks}
\newtheorem*{egs}{Examples}

\theoremstyle{plain}
\newtheorem{thm-tweak}[subsection]{Theorem}
\Crefname{thm-tweak}{Theorem}{Theorems}
\newtheorem{lemma-tweak}[subsection]{Lemma}
\Crefname{lemma-tweak}{Lemma}{Lemmas}
\newtheorem{cor-tweak}[subsection]{Corollary}
\Crefname{cor-tweak}{Corollary}{Corollaries}
\newtheorem{prop-tweak}[subsection]{Proposition}
\Crefname{prop-tweak}{Proposition}{Propositions}
\newtheorem{conj-tweak}[subsection]{Conjecture}
\Crefname{conj-tweak}{Conjecture}{Conjectures}

\theoremstyle{definition}
\newtheorem{defn-tweak}[subsection]{Definition}
\Crefname{defn-tweak}{Definition}{Definitions}
\newtheorem{eg-tweak}[subsection]{Example}
\Crefname{eg-tweak}{Example}{Examples}
\newtheorem*{rems-tweak}{Remarks}
\newtheorem{rem-tweak}[subsection]{Remark}
\Crefname{rem-tweak}{Remark}{Remarks}

\newtheoremstyle{subsection-tweak}
   {11pt}
   {3pt}%
   {}
   {}%
   {\bfseries}
   {}%
   {.5em}
   {\thmnumber{\@{#1}{}\@{#2}.}%
    \thmnote{~{\bfseries#3.}}}

\theoremstyle{subsection-tweak}
\newtheorem{pp}[numberingbase]{}

\theoremstyle{subsection-tweak}
\newtheorem{pp-tweak}[subsection]{}

\numberwithin{equation}{numberingbase}

\makeatletter
\def\@tocline#1#2#3#4#5#6#7{
    \begingroup 
    \@ifempty{#4}{%
    }{%
    }%

    \parindent\z@ \leftskip#3\relax \advance\leftskip\@tempdima\relax
    #5\hskip-\@tempdima
      \ifcase #1
       \or\or \hskip 2em \or \hskip 1em \else \hskip 3em \fi%
      #6\nobreak\relax
    \dotfill\hbox to\@pnumwidth{\@tocpagenum{#7}}\par
    \nobreak
    \endgroup
  }
 \def\l@section{\@tocline{1}{0pt}{1pc}{}{}}

\renewcommand{\tocsection}[3]{%
  \indentlabel{\@ifnotempty{#2}{\makebox[1.3em][l]{%
    \ignorespaces#1 \bfseries{#2}.\hfill}}}\bfseries{#3}
    \vspace{-1.5pt}}

\makeatother

\makeatletter
\newcommand\appendix@section[1]{%
  \refstepcounter{section}%
  \orig@section*{Appendix \@Alph\c@section. #1}%
}
\let\orig@section\section
\g@addto@macro\appendix{\let\section\appendix@section}
\makeatother

\DeclareMathVersion{normal2}

\begin{document}

\title{THE GROTHENDIECK--SERRE CONJECTURE OVER VALUATION RINGS}


\author{NING GUO}
      \address{St. Petersburg branch of V. A. Steklov Mathematical Institute, Fontanka 27, 191023 St. Petersburg, Russia}
      \address{D\'{e}partement de Math\'{e}matiques, Universit\'{e} Paris-Saclay, Orsay Cedex, France}
      \email{guo.ning@eimi.ru}
      \date{\today}

      \begin{abstract}
      In this article, we establish the Grothendieck--Serre conjecture over valuation rings: for a reductive group scheme $G$ over a valuation ring $V$ with fraction field $K$, a $G$-torsor over $V$ is trivial if it is trivial over $K$.
      This result is predicted by the original Grothendieck--Serre conjecture and the resolution of singularities.
      The novelty of our proof lies in overcoming subtleties brought by general nondiscrete valuation rings.
      By using flasque resolutions and inducting with local cohomology, we
      prove a non-Noetherian counterpart of Colliot-Th\'el\`ene--Sansuc's case of tori.
      Then, taking advantage of techniques in algebraization, we obtain the passage to the Henselian rank one case.
      Finally, we induct on Levi subgroups and use the integrality of rational points of anisotropic groups to reduce to the semisimple anisotropic case, in which we appeal to properties of parahoric subgroups in Bruhat--Tits theory to conclude.
      In the last section, by using extension properties of reflexive sheaves on formal power series over valuation rings and patching of torsors, we prove a variant of Nisnevich's purity conjecture.
      \end{abstract}
      
      \small\maketitle
      
      \hypersetup{
          linktoc=page,     
      }
      \renewcommand*\contentsname{}
      \tableofcontents

\section{The Grothendieck--Serre conjecture and Zariski's local uniformization}
      Originally conceived by J.-P. Serre \cite{Ser58}*{p.~31, Rem.} and A. Grothendieck \cite{Gro58}*{pp.~26--27, Rem.~3} in 1958, the prototype of the Grothendieck--Serre conjecture predicted that for an algebraic group $G$ over an algebraically closed field $k$, a $G$-torsor over a nonsingular $k$-variety is Zariski-locally trivial if it is generically trivial.
      With its subsequent generalization to regular base schemes by A. Grothendieck \cite{Gro68}*{Rem.~1.11.a} and the localization by spreading out, the conjecture became the following.
      \bconj[Grothendieck--Serre]\label{GSconj}
      For a reductive group scheme $G$ over a regular local ring $R$ with fraction field $K$, the following map between nonabelian \'etale cohomology pointed sets has trivial kernel\ucolon
      \[
      H^1_{\et}(R,G)\ra H^1_{\et}(K,G);
      \]
      in other words, a $G$-torsor over $R$ is trivial if its restriction over $K$ is trivial.
      \econj
      Diverse variants and cases of \Cref{GSconj} were derived in the last decades.
      A nice survey of the topic is \cites{Ces22}.
      For state-of-the-art results, a more general variant of \Cref{GSconj} over regular semilocal rings containing fields was established by Panin and Fedorov--Panin (\cite{Pan20,FP15}); \v{C}esnavi\v{c}ius \cite{Ces21} settled the unramified quasi-split case (the prior split case is \cite{Fed22}); recently, Guo--Liu \cite{GL23} proved the conjecture for constant group schemes and the smooth projective case was proved by Guo--Panin--Stavrova \cite{GP23,PS23a, PS23b}.
      The goal of this article is to settle the analogue of \Cref{GSconj} when $R$ is instead assumed to be a valuation ring.
      This variant is expected because of the following consequence of the resolution of singularities conjecture, a weak form of Zariski's local uniformization.
      \bconj[Zariski]\label{local-unif}
      Every valuation ring is a filtered direct limit of regular local rings.
      \econj
      Even though \Cref{local-unif} is weaker than Zariski's local uniformization, all its known results come from resolutions or alternations. 
      For a variety $X$ over a field $k$, when $\mathrm{char} k=0$, the local uniformization was resolved by Zariski \cite{Zar40}; when $\mathrm{char}\,k>0$, it was proved for $3$-folds \cite{Abh66, Cut09, CP08, CP09} and surfaces \cite{Abh56}. 
      Temkin \cite{Tem13} achieved the local uniformization after taking a purely inseparable extension of function fields.
      For a valuation ring $V$ whose fraction field $K$ has no degree $p$ extensions (\emph{e.g.}, $K$ is algebraically closed) where $p$ is the residue characteristic, \Cref{local-unif} follows from $p$-primary alterations \cite{Tem17}.
      When $\dim X\geq 4$ and $\mathrm{char}\,k>0$, the local uniformization is widely open.
      
      By assuming \Cref{local-unif}, a limit argument \cite{Gir71}*{VII, 2.1.6} reduces the Grothendieck--Serre over valuation rings to  \Cref{GSconj}.
      In particular, Conjectures \ref{GSconj} and \ref{local-unif} predict the following main result.
      \bthm\label{GSVal}
      For a reductive group scheme $G$ over a valuation ring $V$ with fraction field $K$, the map
      \[
      \x{\tag{$\diamondsuit$}  $H^1_{\et}(V,G)\ra H^1_{\et}(K,G) \label{GSV}$\qq  is injective.}
      \]
      \ethm
      The special case of \Cref{GSVal} when $G$ is an orthogonal group for a nondegenerate quadratic form and $V$ is a valuation ring in which $2$ is invertible was proved in \cite{CTS87}*{6.4} and \cite{CLRR80}*{Thm.~4.5}. 
      
      Besides its connection to the resolution of singularities, the considered variant \Cref{GSVal} offers a few glimpses of the behavior of torsors in the nonarchimedean geometry (more precisely, the rigid-analytic geometry), where the building blocks are affinoids over fraction fields of certain valuation rings (indeed, nonarchimedean fields) and valuation rings usually emerge as rings of definition in Huber pairs.
      Not to mention, the simplest objects in perfectoid spaces, perfectoid fields, are required to be \emph{nondiscrete} valued fields, whose valuation rings are non-Noetherian.
      Also, the following proposition shows that the Grothendieck--Serre over valuation rings yields patching of torsors with respect to arc-covers (\emph{Cf.}~\cite{BM21}).
      \bprop[\Cref{arc}]
      For a valuation ring $V$ of rank $n>0$, the prime $\fp\subset V$ of height $n-1$, and a reductive $V$-group scheme $G$, the following map
      \[
       \x{$\mathrm{Im}(G(V_{\fp})\ra G(\kappa(\fp)))\cdot \mathrm{Im}(G(V/\fp)\ra G(\kappa(\fp)))\surjects G(\kappa(\fp))$\qq is surjective.}
      \]
      \eprop

      The non-Noetherianness of general valuation rings introduces considerable subtleties, even when $G$ is a torus.
      Namely, in this case we can no longer adopt the method of \cite{CTS87}*{4.1} and need to devise alternative arguments.
      For instance, a crucial ingredient of \emph{ibid.} is the exact sequence of \'etale sheaves
       \be\label{fes}
       0\ra \bG_{m, S}\ra i_{\ast}(\bG_{m,\xi})\ra \oplus_{x\in S^{(1)}}i_{x \ast}(\un{\b{Z}}_{x})\ra 0,
       \ee
      where $S$ is a semilocal regular scheme with the union of generic points $i: \xi\ra S$ and $x$ ranges over the points of codimension $1$.
      Being used in the proof of \emph{op.~cit.}, 2.2, however, the short exact sequence (\ref{fes}) fails for general valuation rings.
      For a valuation ring with fraction field $K$ and value group $\GG$, we have 
      $$0\ra V^{\times}\ra K^{\times}\ra \GG \ra 0,$$
      where the abelian group $\GG$ is typically infinitely generated, rendering the arguments in \cite{CTS78, CTS87} knotty to emulate.
      To circumvent this, after using a flasque resolution of tori, we apply local cohomology techniques to induct on the Krull dimension of the valuation ring.
      This reduces us to the following:
       \[
        \x{\tag{$\ast$} for a flasque torus $F$ over a valuation ring $(V,\fm_V)$ of finite rank, we have\qq  $H^2_{\fm_V}(V,F)=0$.}
       \]
      For a flasque torus with character group $\GGL$, by definition (\S\ref{fr}), the Galois action on $\GGL$ has special properties, so certain Galois cohomology of $\GGL$ vanishes, which leads to the vanishing of local cohomology ($\ast$) and therefore the case of tori:
      \bprop[\Cref{GS-tori}]
          For a torus $T$ over a valuation ring $V$ with fraction field $K$,
      \[\textstyle \x{the map \q $H^1_{\et}(V,T)\hra H^1_{\et}(K,T)$ \q is injective.}
      \]
      For a multiplicative type group $M$ of finite type over $V$, the map between pointed sets of fpqc cohomology
      \[
      \x{$H^1_{\fpqc}(V,M)\hra H^1_{\fpqc}(K,M)$\q is injective.}
      \]
      \eprop
      This case of tori, in turn, yields the simplest case of the product formula stated in (\ref{pdx}) below (or, see \Cref{pd-tori}), which is essential for further reduction of \Cref{GSVal}.

      A practical advantage of Henselian rank-one valuation rings is that several techniques of Bruhat--Tits theory, especially in \cite{BrT2}*{\S4-5}, become available.
      The goal of \S\ref{wk-ap} and \S\ref{passage-h} is to reduce \Cref{GSVal} to this case: after a limit argument that leads to the case of finite rank, we induct on the rank $n$ of a valuation ring $V$ by patching torsors.
      The induction hypothesis implies that our $G$-torsor over $V$ is a gluing of trivial torsors.
      For this gluing, we choose an $a\in V$ such that the $a$-adic completion $\hva$ is a rank-one Henselian valuation ring with $\hka\ce \Frac \hva$; so that, $\via$ is a valuation ring of rank $n-1$.
      Similar to the Beauville--Laszlo's gluing of bundles, our patching is reformulated as the product formula
      \begin{equation}\label{pdx}
          \tst G(\hka)=\mathrm{Im}\bigl(G(\via)\ra G(\hka)\bigr)\cdot G(\hva).
      \end{equation}
  The strategy for proving this formula is a ``d\'evissage'' that establishes approximation properties of certain subgroups of $G_{\nhva}$.
      In this procedure, techniques of algebraization \cite{BC22}*{\S2} play an important role, especially for a Harder-type approximation (see \S\ref{wk-ap}) and the following higher rank counterpart of \cite{Pra82}.
      \bprop[\Cref{aniso-int-rat}]\label{a-i-r}
      For a reductive anisotropic group scheme $G$ over a Henselian valuation ring $V$ with fraction field $K$, we have $ G(V)=G(K)$.
      \eprop
      Based on its special case when $K=\hka$ is complete due to Maculan \cite{Mac17}*{Thm.~1.1}, our approach to \Cref{a-i-r} is a reduction to completion that rests on techniques of algebraization to approximate schemes characterizing the anisotropicity of $G_{\nhva}$. 
      Indeed, \Cref{a-i-r} is an anisotropic version of the product formula (\ref{pdx}).
      \Cref{a-i-r} is helpful, not only for the reduction to the Henselian rank-one case, but also for the induction on Levi subgroups when reducing to the semisimple anisotropic case in \S\ref{passage-ss-ani}.
      After these reductions, we transfer \Cref{GSVal} into the injectivity of a map of Galois cohomologies.
      We conclude by taking advantage of properties of parahoric subgroups in Bruhat--Tits theory, see \Cref{final-proof}.
      
      In addition to techniques of algebraization, another crucial element of our reduction to the Henselian rank-one case is a lifting property of maximal tori of reductive group schemes.
      \blem[\Cref{extend-tor}]\label{lift-tor}
      For a reductive group scheme $G$ over a local ring $(R,\kappa)$ with a maximal $\kappa$-torus $T$, if the cardinality of $\kappa$ is at least $\dim(G^{\mathrm{ad}})$, then $G$ has a maximal $R$-torus $\sT$ such that 
      \[
      \tst \sT_{\kappa}=T.
      \]
      \elem
      This strengthens a result of Grothendieck \cite{SGA3II}*{XIV, 3.20} that a maximal torus of a reductive group scheme exists Zariski-locally on the base.
      By a correspondence of maximal tori and regular sections, the novelty is to lift regular sections instead of merely proving their existence Zariski-locally.
      Depending on inspection of the reasoning for \emph{ibid.}, the key point is \cite{Bar67}, which guarantees that Lie algebras over fields with large cardinalities contain regular sections.
      For lifting regular sections, we need the functorial property of Killing polynomials.
      Indeed, Killing polynomials over rings were defined ambiguously in
      in the original literature, see \cite{SGA3II}*{XIV, 2.2}.
      Therefore, to establish \Cref{lift-tor}, we first add the supplementary details \S\ref{C} for Killing polynomials over rings.
      Subsequently, for a Lie algebra with locally constant nilpotent rank, we use the functoriality of Killing polynomials to deduce the openness of the regular locus. This openness permits us to lift regular sections, which amounts to lifting maximal tori.   

      In \S\ref{purity}, we acquire a variant of Nisnevich's purity conjecture \cite[1.3]{Nis89}, whose statement is the following.
      \bconj[Nisnevich's purity]\label{nisconj}
       For a reductive group scheme $G$ over a regular local ring $R$ with a regular parameter $f\in \fm_R\backslash \fm_R^2$, every Zariski-locally trivial $G$-torsor over $R[\f{1}{f}]$ is trivial, that is, we have
       \[
       \tst H^1_{\Zar}(R[\f{1}{f}],G)=\{\ast\}.
       \]
      \econj 
      This conjecture generalizes Quillen's conjecture \cite[Comments]{Qui76} when $G=\GL_{n}$ and was proved by Gabber \cite{Gab81} for $G=\GL_n$ and $\mathrm{PGL}_n$ when $\dim R\leq 3$.
      In this article, we consider a variant: for a valuation ring $V$ and its ring of formal power series $V\fps{t}$, we let $R=V\fps{t}$ and $f=t$, hence $R[\f{1}{f}]=V\lps{t}$.
      \bprop[\Cref{nis}]\label{int:nis}
      For a reductive group scheme $G$ over a valuation ring $V$, every Zariski-locally trivial $G$-torsor over $V\lps{t}$ is trivial, that is, we have
\[
   H^1_{\Zar}(V\lps{t},G)=\{\ast\}.
\]

      \eprop
      This \Cref{int:nis} follows from the injectivity of the map $H^1_{\et}(V\lps{t},G)\ra H^1_{\et}(K\lps{t},G)$ proved in \Cref{lps}. 
      In fact, by cohomological properties of reflexive sheaves (see \ref{ref:pre}), every \'etale $\GL_n$-torsor over $V\lps{t}$ is trivial. 
      With an embedding $G\hra \GL_{n}$, we obtain \Cref{int:nis} by patching torsors.

      \bpp[Notation and conventions]\label{convention}
      For various notions and properties about valuation rings and valued fields, see Appendix \S\ref{apdx}.
      We adopt the notion in \cite{SGA3IIInew} for reductive group schemes: they are group schemes smooth affine over their base schemes, such that each geometric fiber is connected and contains no normal subgroup that is an iterated extension of $\bG_a$.  
      For a valuation ring $V$, we denote by $\fm_V$ the maximal ideal of $V$. 
      When $V$ has finite rank $n$, for the prime $\fp\subset V$ of height $n-1$ and $a\in \fm_V\backslash \fp$,  we denote by $\hva$ the $a$-adic completion of $V$. 
      For a module $M$ finitely generated over a topological ring $A$, we endow $M$ with the \emph{canonical topology} as the quotient of the product topology via $\pi\colon A^{\oplus n}\surjects M$. 
      By \cite{GR18}*{8.3.34}, this topology on $M$ is independent of the choice of $\pi$.
      \epp
      \subsection*{Acknowledgements} I especially thank my advisor K\k{e}stutis \v{C}esnavi\v{c}ius for his kindness, helpful advice, and extensive comments for revising. 
      I thank Jean-Louis Colliot-Thélène for his several insightful comments on generalizing the case of tori to groups of multiplicative type.
      I thank Ofer Gabber and Philippe Gille for helpful conversations about properties of anisotropic groups. 
      I thank Jean-Louis Colliot-Thélène, Uriya First, Philippe Gille, David Harari, David Harbater, Fei Liu, Matthew Morrow, Ivan Panin, and Olivier Wittenberg for their careful reading and useful comments on this article.
      I thank Fei Xu, Yang Cao, and Xiaozong Wang for useful information. 
      I also thank Kazuhiro Ito, Hiroki Kato, Arnab Kundu, Yisheng Tian, Yifei Zhao, and Jiandi Zou for helpful  discussions and remarks. 
      I thank referees for helpful comments and catching imprecisions in previous drafts for me to improve.
      This article is supported by the EDMH doctoral program.
      This work was done under support of the grant №075-15-2022-289 for creation and development of Euler International Mathematical Institute.
       \section{The case of tori}
      The goal of this section is to prove the Grothendieck--Serre conjecture over valuation rings for tori, a non-Noetherian counterpart of Colliot-Th\'el\`ene--Sansuc's result \cite{CTS87}*{4.1}, then we extend it to groups of multiplicative type (\Cref{GS-tori}~\ref{mult-brauer-ii}).
      Colliot-Th\'el\`ene and Sansuc defined flasque resolutions of tori over arbitrary base schemes, yielding several cohomological properties of tori over regular schemes.
      In particular, they proved that for a torus $T$ over a semilocal regular ring $R$ with total ring of fractions $K$,
      \begin{equation}\label{gs-old-tori}
        \x{ the map \q $H^1_{\et}(R,T)\hra H^1_{\et}(K,T)$\qq is injective,}
      \end{equation}
      which is a stronger version of the Grothendieck--Serre conjecture for tori, see \cite{CTS87}*{4.1}.
      Nevertheless, if we substitute $R$ in (\ref{gs-old-tori}) with a valuation ring $V$, then the method in \emph{ibid.} does not work any more because of the non-Noetherianness of       $V$.
      Seeking an alternative argument in this case, we induct on the rank of $V$ and use local cohomology.
      This case of tori obtained in \Cref{GS-tori} is crucial for subsequent steps of the proof of \Cref{GSVal}, such as for patching torsors (see Propositions \ref{decomp-gp} and \ref{rank-one-kernel-trivial}).
      
     \bpp[Group schemes of multiplicative type]\label{mult}
      For a scheme $S$ and an $S$-group scheme $G$, the \emph{Cartier dual} of $G$ is an fpqc sheaf $\sD_S(G)\ce \sHom_{\x{$S$-gr.}}(G,\bG_{m,S})$.
      Recall \cite{SGA3II}*{IX, 1.1} that $G$ is \emph{of multiplicative type}, if every $s\in S$ has an fpqc neighborhood $U$ such that $G_U\simeq \sD_U(M_U)=\sHom_{\x{$U$-gr.}}(M_U,\bG_{m,U})$ for a commutative group $M$. 
      An $S$-group $G$ of multiplicative type is \emph{isotrivial}, if there exists a finite \'etale surjective morphism $S\pr\ra S$ such that $\sD_{S\pr}(G_{S\pr})$ is a constant commutative group on each connected component of $S\pr$ (\cite{SGA3II}*{IX, 1.4.1}). 
      Assume that $S$ is \emph{connected}.
       One can replace $S\pr$ by one of its connected component and apply \SP{0BN2} to find an $S$-morphism $S\prpr\ra S\pr$ of schemes for a Galois cover $S\prpr$ of $S$ (by \cite{SGA1new}*{V, 5.11}, $S\prpr$ is a connected $\un{\GG}_S$-torsor for a finite group $\GG$). 
       Then, since $\GG$ has finitely many quotients, there is a \emph{minimal Galois cover} $\wt{S}/S$ such that $\sD_{\wt{S}}(G_{\wt{S}})$ is constant: the minimality of $\wt{S}/S$ means that there are no nontrivial Galois subcovers $\wt{S}\ra \wt{S\pr} \ra S$ such that $\sD_{\wt{S\pr}}(G_{\wt{S\pr}})$  is constant.
       We also say that $\wt{S}/S$ is a \emph{minimal Galois cover splitting} $G$ (or, \emph{such that} $G_{\wt{S}}$ \emph{splits}).
      Moreover, since $S$ is assumed to be connected, for every geometric point $\ov{s}\colon \Spec \Omega\ra S$ of $S$ with fundamental group $\pi\ce \pi_1^{\et}(S,\ov{s})$, where $\Omega$ is an algebraically closed field, there is an anti-equivalence \cite{SGA3II}*{X, 1.2}
      \[ \hkh{\begin{matrix}
                     \x{isotrivial multiplicative}\\
                     \x{type $S$-groups}
                   \end{matrix} } \isoto
            \hkh{\begin{matrix}
                   \x{$\pi$-modules with} \\
                   \x{continuous actions}
                 \end{matrix} } \q G\mapsto \sM(G)\ce\sD_{\ov{s}}(G_{\ov{s}})=\mathrm{Hom}_{\x{$\Omega$-gr.}}(G_{\ov{s}}, \bG_{m,\ov{s}}).
      \]
      In particular, the category of isotrivial $S$-tori is anti-equivalent to the category of finite type $\b{Z}$-lattices with continuous $\pi$-actions.
      So, every isotrivial $S$-torus $T$ of rank $n$ corresponds to an equivalence class of 
      \[
      \text{representations \q $\rho_T\colon \pi\ra \GL_n(\b{Z})$\q such that $\ker\rho_T\subset \pi$ is an open normal subgroup.}
      \]
      If $\rho_T$ and $\rho_T\pr$ are in the same equivalence class, then $\ker \rho_T=\ker \rho_T\pr$.
      The finite quotient $\GG\ce \pi/\ker \rho_T$ then yields a minimal Galois cover  $\wt{S}/S$ splitting $T$ with Galois group $\GG$ and $\pi_1^{\et}(\wt{S})\simeq \ker \rho_T$. 
      Hence, \emph{all minimal Galois covers splitting $T$ are isomorphic to each other via the Galois group ${\GG}$-action}.
      \epp
      \blem\label{lift-split-rank} 
      For an irreducible geometrically unibranch scheme $S$ of function field $K$ and an $S$-torus $T$, 
      \[
      \text{$T$ contains $\bG_{m,S}^k$\q if and only if \q $T_K$ contains $\bG_{m,K}^k$.}
      \]
      \bpf 
      It suffices to assume that $\bG_{m,K}^k\subset T_K$ and to deduce that $\bG_{m,S}^k\subset T$.
      Let $\ov{\eta}$ be a geometric point over the generic point $\Spec K\overset{\eta}{\ra} S$.
      We have $\sM(T)=\Hom_{\ov{\eta}\x{-gr.}}(T_{\ov{\eta}}, \bG_{m,\ov{\eta}})= \sM(T_K)$.
      Note that $\bG_{m,K}^k$ corresponds to a quotient lattice $\Lambda$ of $\sM(T_K)$ such that $\Lambda$ is of rank $k$ with trivial $\pi_1^{\et}(K)$-action.
      On the other hand, by \SP{0BQI}, the natural map $\pi_1^{\et}(K)\surjects \pi_1^{\et}(S)$ is surjective.  
      Therefore, $\sM(T)$ has a quotient lattice that has rank $k$ with trivial $\pi_1^{\et}(S)$-action.
      This implies that $\bG_{m,S}^k\subset T$. 
      \epf 
      \elem
      Recall \cite{EGAI}*{2.1.8} that a scheme $S$ is \emph{locally integral}, if for every $s\in S$, the local ring $\sO_{S,s}$ is integral.
      Hence, by definition, every connected component of $S$ is both an open and closed subset of $S$.
      With this notion, we generalize Grothendieck's result \cite{SGA3II}*{X, 5.16} by relaxing its Noetherian constraint.

      \blem\label{isotrivial} 
      For a locally integral, geometrically unibranch scheme $S$, every $S$-group scheme $M$ of multiplicative type and of finite type is isotrivial.
      In particular, for every torus $T$ over a normal domain $R$, there is a minimal Galois cover $\wt{R}$ of $R$ such that $T_{\wt{R}}$ splits.
      \elem 
      \bpf 
      Since every connected component of $S$ is open, we may assume that $S$ is connected. 
      Then, $M$ is fpqc locally of the form $\sD(H)$ for a finite type abelian group $H$ (determined by $M$).
      For $P\ce \un{\Isom}_{S\x{-gr.}}(M,\sD_S(H))$, our goal is to find a finite \'etale cover $S\pr\ra S$ such that $P(S\pr)\neq \emptyset$.
      By \cite{SGA3II}*{X, 5.8, 5.10~(i)}, $P$ is representable by a clopen subscheme of $\un{\Hom}_{S\x{-gr.}}(M,\sD_S(H))$ and there is an \'etale surjective morphism $\wt{S}\ra S$ such that $P_{\wt{S}}$ is a disjoint union of copies of $\wt{S}$. 
      In particular, $P$ is $S$-\'etale.
      By \cite{EGAIV4}*{18.8.15, 18.10.7}, $\wt{S}$ is locally integral and geometrically unibranch.
      We prove the following.
      \bcl 
      Every irreducible component $P_i$ of $P$ is \emph{finite} \'etale over $S$.
      \ecl
      \bpf[Proof of the claim] 
      Let $\eta\in S$ be the generic point and let $\xi_i$ be the generic point of $P_i$.
      By \cite{EGAIV2}*{2.3.4}, the $S$-flatness of $P$ implies that every $\xi_i$ lies over $\eta$. 
      Therefore, $(P_i)_{\eta}$ is the closure of $\xi_i$ in $P_{\eta}$.
      The quasi-finiteness of $P\ra S$ implies that $P_{\eta}$ is discrete, so we have $(P_i)_{\eta}=\{\xi_i\}$.
      On the other hand, since $S$ is integral and geometrically unibranch, by \cite{EGAIV4}*{18.10.7}, all $P_i$ are geometrically unibranch, and 
      \[
      \tst P=\bigsqcup_{\xi_i\in P_{\eta}}P_i.
      \]
      Therefore, every $P_i$ is clopen in $P$.
      Since it suffices to show that each $(P_i)_{\wt{S}}$ is $\wt{S}$-finite,  note that every connected component of $\wt{S}$ is open, we may assume that $\wt{S}$ is connected so that $P_{\wt{S}}\cong \bigsqcup_{\Psi}{\wt{S}}$ for a set $\Psi$.
      Each $P_i\subset P$ satisfies that $(P_i)_{\wt{S}}\cong \bigsqcup_{\Phi_i}\wt{S}$ for a subset $\Phi_i\subset \Psi$.
      As $(P_i)_{\eta}=\{\xi_i\}$ is a single point, this forces that $\Phi_i$ is finite.
      Consequently, the base change $(P_i)_{\wt{S}}$ is finite over $\wt{S}$, so $P_i$ is $S$-finite.
      \epf 
      As $S$ is connected and all $P_i\ra S$ are finite \'etale, take $S\pr\ce P_i$, whose image is $S$.
      The canonical embedding $S\pr\hra P$ then induces a section of $P_{S\pr}\ra S\pr$, so we get $M_{S\pr}\simeq \sD_{S\pr}(H)$, as desired.
      \epf 
      \bprop\label{bc-galois}
      Let $X$ be a connected scheme, let $T$ be an isotrivial $X$-torus, and let $Y\ra X$ be a minimal Galois cover splitting $T$.
      For a morphism $f\colon X\pr\ra X$ of connected schemes, every connected component of $Y\pr\ce Y\times_X X\pr$ is a minimal Galois cover splitting $T_{X\pr}$.
      \eprop
      \bpf
      Let $\GG\ce \Aut_X(Y)$ be the Galois group of $Y/X$, then $Y$ is a $\un{\GG}_X$-torsor on $X$, and $Y\pr$ is a $\un{\GG}_{X\pr}$-torsor on $X\pr$. In particular, $\GG$ acts transitively on each $X\pr$-fiber of $Y\pr$, hence induces isomorphisms among connected components of $Y\pr$. 
      We choose a geometric point $\eta\pr\ra Y\pr$, and denote its composites as $\eta\ra Y$, $\xi\pr\ra X\pr$, and $\xi\ra X$, respectively.
      Recall \SP{0BND} that the fiber functors $F_{\xi}\colon \mathrm{F\'Et}_{X}\isoto {\mathrm{Finite}\x{-}\pi_1^{\et}(X,\xi)\x{-sets}}$ and $F_{\xi\pr}\colon \mathrm{F\'Et}_{X\pr}\isoto {\mathrm{Finite}\x{-}\pi_1^{\et}(X\pr,\xi\pr)\x{-sets}}$ are equivalences of categories.
      Besides, $f$ induces a continuous homomorphism $f_{\ast}\colon \pi_1^{\et}(X\pr,\xi\pr)\ra \pi_1^{\et}(X,\xi)$ of profinite groups, fitting into the following commutative diagram
      \[
      \begin{tikzcd}
\mathrm{F\'Et}_{X} \arrow[d, "F_{\xi}"'] \arrow[rr, "\text{base change}"] &  & \mathrm{F\'Et}_{X'} \arrow[d, "F_{\xi'}"] \\
{\mathrm{Finite}\x{-}\pi_1^{\et}(X,\xi)\x{-sets}} \arrow[rr, "f^{\ast}"]                    &  & {\mathrm{Finite}\x{-}\pi_1^{\et}(X\pr,\xi\pr)\x{-sets}}.          
\end{tikzcd}
      \]
      Thus, we have $F_{\xi\pr}(Y\pr)=f_{\ast}F_{\xi}(Y)=F_{\xi}(Y)=\GG$ set-theoretically and the following short exact sequence 
      \[
      1\ra \pi_1^{\et}(Y,\eta)\ra \pi_1^{\et}(X,\xi)\ra \GG\cong\Aut_{\GG\x{-set}}(F_{\xi}(Y)) \ra 1.
      \]
      By the commutative diagram above, the $\pi_1^{\et}(X\pr,\xi\pr)$-action on $F_{\xi\pr}(Y\pr)$ is equal to the $\pi_1^{\et}(X\pr,\xi\pr)$-action on $F_{\xi}(Y)$ via the composite $\pi_1^{\et}(X\pr,\xi\pr)\overset{f_{\ast}}{\ra}\pi_1^{\et}(X,\xi)\surjects \GG$, whose image is denoted by $\GG\pr\subset \GG$.
      The surjection $\pi_1^{\et}(X\pr,\xi\pr)\surjects \GG\pr$ gives rise to the $\pi_1^{\et}(X\pr,\xi\pr)$-set structure on $F_{\xi\pr}(Y\pr)$. 
      Precisely, the $\pi_1^{\et}(X\pr,\xi\pr)$-action on $F_{\xi\pr}(Y\pr)$ is just the restriction $\GG\pr\times \GG\ra \GG$ of $\GG\times \GG\ra \GG$, leading to the coset decomposition for $\GG\pr\subset \GG$
      \[
      \tst \GG=\bigsqcup_{\gamma\in \GG\pr\backslash \GG}(\GG\pr\cdot \gamma)
      \]
      so that all left $\GG\pr$-actions on $\GG\pr\cdot \gamma$ are simply transitive and all $\GG\pr\cdot \gamma$ have the same $\GG\pr$-set structure.
      Hence, the equivalence $F_{\xi\pr}\colon \mathrm{F\'Et}_{X\pr}\isoto {\mathrm{Finite}\x{-}\pi_1^{\et}(X\pr,\xi\pr)\x{-sets}}$ (combined with \SP{03SF}) 
      implies that $(\GG\pr\cdot \gamma)_{\gamma\in \GG\pr\backslash \GG}$ correspond to Galois covers $(Y\pr_{\gamma})_{\gamma\in \GG\pr\backslash \GG}$ of $X\pr$ that are isomorphic to each other.
      Further, the finite $\pi_1^{\et}(X\pr,\xi\pr)$-set $F_{\xi\pr}(Y\pr)$ corresponds to $Y\pr$, which decomposes into connected components
      \[
     \tst  Y\pr=\bigsqcup_{\gamma\in \GG\pr\backslash \GG}Y\pr_{\gamma},
      \]
      where $Y\pr_{\gamma}$ are Galois covers of $X\pr$ with Galois group $\GG\pr$.
      If $\eta\pr\ra Y\pr$ factors through $Y\pr_{\gamma_0}$, then the following
      \[
      1\ra \pi_1^{\et}(Y_{{\gamma}_0}\pr, \eta\pr)\ra \pi_1^{\et}(X\pr,\xi\pr)\ra \GG\pr=\mathrm{Gal}(Y_{\gamma_0}\pr/X\pr)\ra 1
      \]
      is a short exact sequence.
      Since the torus $T$ induces a representation $\rho_T\colon \pi_1^{\et}(X,\xi)\ra \GL(\b{Z}^n)$ with the image $\GG$, where $\b{Z}^n\simeq \Hom_{\x{$\xi$-gr.}}(T_{\xi},\bG_m)$, its base change $T_{X\pr}$ induces a representation $f_{\ast}\circ \rho_T\colon \pi_1^{\et}(X\pr,\xi\pr)\ra \GL(\b{Z}^n)$. By construction of $\GG\pr$, we have $\GG\pr=\mathrm{Im}(f_{\ast}\circ \rho_T)$. 
      So the desired minimality of $Y_{\gamma_0}\pr$ amounts to the equality $\GG\pr=\pi_1^{\et}(X\pr, \xi\pr)/\pi_1^{\et}(Y_{\gamma_0}\pr,\eta\pr)$, which follows from the last displayed short exact sequence. 
      \epf

      \bpp[Flasque resolution of tori]\label{fr}
      The concepts of quasitrivial and flasque tori are rooted in two special Galois modules that serve as character groups: permutation and flasque modules.
      For a finite group $G$, let $\cL_G$ be the category of $G$-modules that are finite type $\b{Z}$-lattices.
      If a module $M\in \cL_G$ has a $\b{Z}$-basis on which $G$ acts via permutations, then $M$ is a \emph{permutation module}; in this case, $M\simeq \oplus_{i}\b{Z}[G/H_i]$ for certain subgroups $H_{i}\subset G$.
      If a module $N\in \cL_G$ satisfies $H^1(G, \mathrm{Hom}_{\b Z}(N,Q))=0$ for any permutation module $Q$, then $N$ is a \emph{flasque module}.
      For example, a trivial $G$-module $Q_0\in \cL_G$ is a permutation module and $H^{1}(G,\Hom_{\b{Z}}(N,Q_0))=0$ for any flasque $G$-module $N$.
      For a scheme $S$ and an $S$-torus $T$, if every connected component $Z$ of $S$ has a Galois cover $Z\pr\ra Z$ with Galois group $G$ splitting $T$ such that the $G$-module $\sD_S(T)(Z\pr)$ is flasque (resp., permutation), then $T$ is \emph{flasque} (resp., \emph{quasitrivial}).
      When $S$ is connected, every quasitrivial torus is a finite product of Weil restrictions $\mathrm{Res}_{S\pr_i/S}(\bG_{m})$ for finite \'etale  connected covers $S\pr_i\ra S$.
      As proved in \cite{CTS87}*{Thm.~1.3}, for a torus $T$ over a scheme $S$ whose every connected component is open, there is a short exact sequence of $S$-tori, that is, a \emph{flasque resolution} of $T$:
      \begin{equation}\label{flasque-resolution}
        \x{$1\ra F\ra P\ra T\ra 1$,\qqq where $F$ is flasque and $P$ is quasitrivial.}
      \end{equation}
      \epp
       \blem\label{2-vanish}
      For a flasque torus $F$ over a valuation ring $V$ of finite rank, the local cohomology vanishes:
      \[
      H^2_{\fm_V}(V,F)=0.
      \]
      \elem
      \bpf
      We denote $X=\Spec V$ and $Z=\Spec(V/\fm_V)$.
      Let $n\geq 1$ be the rank of $V$, then $X\backslash Z$ is the spectrum of a valuation ring of rank $n-1$.
      By excision \cite{Mil80}*{III, 1.28}, we may replace $X$ by its Henselization $X^{\mathrm{h}}$.
      For a variable $X$-\'etale scheme $X\pr$ with preimage $Z\pr\ce X\pr\times_{X} Z$, let $\cH^q_{Z}(-,F)$ be the \'etale sheafification of the presheaf  $X\pr\mapsto H^q_{Z\pr}(X\pr, F)$.
      By the local-to-global $E_2$ spectral sequence 
      \[
      \x{\qqqqqqqq\qqqqqqq$H^p_{\et}(X,\cH^q_Z(X,F)) \Ra H^{p+q}_Z(X,F)$,\qqqqq\qq (\cite{SGA4II}*{V, 6.4})}
      \]
      to show that $H^2_{Z}(X,F)=0$, it suffices to obtain the following vanishings
      \[
      H^0_{\et}(X,\cH^2_{Z}(X, F))=H^1_{\et}(X,\cH^1_Z(X, F))=H^2_{\et}(X, \cH^0_Z(X, F))=0.
      \]
      Subsequently, in the following two paragraphs, we calculate $\cH^q_Z(X,F)$ for $0\leq q\leq 2$.
      
      Let $\ov{x}\ra X$ be a geometric point. If $\ov{x}$ factors through $X\backslash Z$, then $\cH^q_Z(X,F)_{\ov{x}}=0$. 
      Now, we take $\ov{x}$ as a fixed geometric point over $\fm_V$, so $\cH^q_Z(X,F)_{\ov{x}}=H^q_{\ov{\fm}_V}(V^{\sh},F)$, where $V^{\sh}$ is the strict Henselization of $V$ with the maximal ideal $\ov{\fm}_V$.
      The local map $V\ra V^{\sh}$ of local rings is faithfully flat (\SP{07QM}) and preserves value groups (\SP{0ASK}).
      Therefore, for the prime $\fp\subset V$ of height $n-1$, there is a unique prime ideal $\fP\subset V^{\sh}$ lying over $\fp$ (that is, $\fp V^{\sh}=\fP$). 
      By \cite{SGA4II}*{V, 6.5}, we have the exact sequence
      \begin{equation}\label{str-loc-ext}
       \cdots\ra H^i_{\et}(V^{\sh}, F) \ra H^i_{\et}((V^{\sh})_{\fP},F)\ra H^{i+1}_{\ov{\fm}_V}(V^{\sh},F)\ra H^{i+1}_{\et}(V^{\sh},F)\ra \cdots.
      \end{equation}
      First, we compute $H^{q}_{\ov{\fm}_V}(V^{\sh},F)$ when $q=0$ and $2$. 
      The injectivity of $H^0_{\et}(V^{\sh},F)\hra H^0_{\et}((V^{\sh})_{\fP},F)$ and  the vanishings of $H^1_{\et}((V^{\sh})_{\fP},F)$ and $H^i_{\et}(V^{\sh},F)$ for $i=1,2$ (see \SP{03QO}) imply the following 
      \begin{equation}\label{stalk-vanish}
      H^0_{\ov{\fm}_V}(V^{\sh},F)=H^2_{\ov{\fm}_V}(V^{\sh},F)=0.
      \end{equation}
      This (\ref{stalk-vanish}) leads to $\cH^0_Z(X,F)=\cH^2_Z(X,F)=0$, so we get $H^0_{\et}(X,\cH^2_Z(X,F))=H^2_{\et}(X,\cH^0_Z(X,F))=0$.

      Next, we calculate $H^{1}_{\ov{\fm}_V}(V^{\sh},F)$. 
      From (\ref{str-loc-ext}) we obtain the following short exact sequence:  
      \[
      0\ra H^0_{\et}(V^{\sh},F)\ra H^0_{\et}((V^{\sh})_{\fP}, F)\ra H^1_{\ov{\fm}_V}(V^{\sh}, F)\ra H^1_{\et}(V^{\sh}, F)=0.
      \] 
      For the Cartier dual $\sD_X(F)$ of $F$, let $\Lambda\ce \sD_X(F)(V^{\sh})$ and $\Lambda^{\vee}\ce \Hom_{\b{Z}}(\Lambda,\b{Z}) $.
      By Cartier duality,
      \[
      \begin{matrix}[1.5]
        H^0_{\et}(V^{\sh},F)\cong F\p{V^{\sh}}\cong \sHom_{\x{$V$-gr.}}(\sD_X(F),\bG_m)(V^{\sh})=\Hom_{\b{Z}}(\Lambda, ({V^{\sh}})^{\times})\cong\Lambda^{\vee}\otimes_{\b{Z}}({V^{\sh}})^{\times}, \\
        \x{\!\!\!\!\!\!\!\!\!\!and similarly,\qqqqq $H^0_{\et}((V^{\sh})_{\fP},F)\cong \Lambda^{\vee}\otimes_{\b{Z}}(V^{\sh})_{\fP}^{\times}$.}
      \end{matrix}
      \]
     The value group $\GG_{V^{\sh}/\fP}$ of $V^{\sh}/\fP$, by \Cref{basic-one}~\ref{exact-seq-val-grps}, is isomorphic to $(V^{\sh})_{\fP}^{\times}/({V^{\sh}})^{\times}$.
     Therefore,
      \[
      H^1_{\ov{\fm}_V}(V^{\sh},F)=(\GGL^{\vee}\otimes_{\b{Z}}(V^{\sh})_{\fP}^{\times})/(\GGL^{\vee}\otimes_{\b{Z}}(V^{\sh})^{\times})\cong \GGL^{\vee}\otimes_{\b{Z}}\GG_{V^{\sh}/\fP}.
      \]
      Since $X$ is Henselian local and $\cH^1_Z(X,F)$ is an abelian sheaf on $X$, by \cite{SGA4II}*{VIII, 8.6}, we have  
      \begin{equation}\label{h1}
        H^1_{\et}(X,\cH^1_Z(X,F)) \cong H^1(\etp(V),H^1_{\ov{\fm}_V}(V^{\sh},F)) \cong  H^1(\etp(V), \mathrm{Hom}_{\b{Z}}(\GGL,\GG_{V^{\sh}/\fP})).
      \end{equation}
      To see the action of $\etp(V)$ on $\mathrm{Hom}_{\b{Z}}(\GGL, \GG_{V^{\sh}/\fP})$, by \Cref{isotrivial}, we first note that the $\etp(V)$-action on $\GGL$ factors through its quotient $\Gal(Y/X)$, where $Y$ is the minimal Galois cover of $X$ splitting $F$. 
      Besides, 
      \[
         \GG_{V^{\sh}/\fP}\stackrel{\x{\SP{05WS}}}{=\joinrel=\joinrel=\joinrel=} \GG_{\p{V/\fp}^{\sh}} \overset{\x{\SP{0ASK}}}{=\joinrel=\joinrel=\joinrel=} \GG_{V/\fp},
      \]
      so $\etp(V)$ acts trivially on $\GG_{V^{\sh}/\fP}\cong\Frac(V/\fp)^{\times}/(V/\fp)^{\times}$.
      Thus, the $\etp(V)$-action on $\mathrm{Hom}_{\b{Z}}(\GGL, \GG_{V/\fp})$ factors through $\Gal(Y/X)$.
      Since $\etp(V)$ is a projective limit of finite groups $\Gal(X_{\alpha}/X)$, where $X_{\alpha}$ ranges over Galois covers of $X$, a limit argument \cite{Ser02}*{I, \S2.2, Cor.~1} reduces (\ref{h1}) to
      \begin{equation}\label{h2}
       \tst  H^1_{\et}(X,\cH^1_Z(X,F))\isom \varinjlim_{\alpha}  H^1\bigl(\Gal(X_{\alpha}/X), \mathrm{Hom}_{\b{Z}}(\GGL,\GG_{V/\fp})^{\etp(X_{\alpha})}\bigr).
      \end{equation}
      We express $\Gamma_{V/\fp}$ as a direct limit of finite type $\b{Z}$-submodules $(\Gamma_i)_{i\in I}$. 
      Since $\GGL$ is $\b{Z}$-finitely presented,
      \begin{equation}\label{colimit}
       \tst \varinjlim_{i\in I}\Hom_{\b{Z}}(\Lambda, \Gamma_{i})\isoto \Hom_{\b{Z}}(\Lambda,\Gamma_{V/\fp}).
      \end{equation}
      Combining the isomorphism (\ref{colimit}) with a limit argument \cite{Ser02}*{I, \S2.2, Prop.~8}, we reduce (\ref{h2}) to
      \begin{equation*}
           \tst \varinjlim_{\alpha}H^1\bigl(\Gal(X_{\alpha}/X), \varinjlim_{i\in I}\Hom_{\b{Z}}(\GGL,\GG_i)^{\etp(X_{\alpha})}\bigr)  =  \varinjlim_{\alpha}\varinjlim_{i\in I}H^1\bigl(\Gal(X_{\alpha}/X), \Hom_{\b{Z}}(\GGL,\GG_i)^{\etp(X_{\alpha})}\bigr).
      \end{equation*}
      It suffices to calculate for a large $\alpha_{0}$ such that $X_{\alpha_0}$ splits $F$.
      In this situation, $\etp(X_{\alpha_0})$ acts trivially on $\mathrm{Hom}_{\b{Z}}(\GGL,\GG_{i})$. 
      Since $F$ is a flasque torus, its character group $\GGL$ is a flasque $\Gal(X_{\alpha_0}/X)$-module.
      As aforementioned, $\Gal(X_{\alpha_0}/X)$ acts trivially on $\GG_{V/\fp}$, so the $\GG_{i}$ are finite type $\b{Z}$-lattices with trivial $\Gal(X_{\alpha_0}/X)$-action.
      The example in \S\ref{fr} implies $H^1(\Gal(X_{\alpha_0}/X),\Hom_{\b{Z}}(\Lambda,\Gamma_i))=0$, which verifies that
      \[
       H^1_{\et}(X, \cH^1_Z(X,F))=0. \qedhere
      \]
      \epf
\bprop\label{GS-tori}\label{mult-brauer}
      For a valuation ring $V$ and a finite type $V$-group scheme $M$ of multiplicative type, 
      \benumr
      \item\label{mult-brauer-i} $H^2_{\fpqc}(V,M)\hra H^2_{\fpqc}(\Frac V,M)$ is injective; in particular, the restriction of Brauer group
      \[
      \text{$\mathrm{Br}(V)\hra \mathrm{Br}(\Frac V)$\q  is injective;}
      \]
      \item\label{mult-brauer-ii} $H^1_{\fpqc}(V,M)\hra H^1_{\fpqc}(\Frac V,M)$ is injective.
      \eenum
      \eprop
      \bpf
      As $V$ is a filtered direct union of valuation subrings of finite rank (\cite{BM21}*{2.22}), a limit argument \cite{SGA4II}*{VII, 5.7} reduces us to the case when $V$ has finite rank $n$.
      Note that for a quasitrivial $V$-torus $P$,  we have $P\isom \prod_{S\pr_i}\mathrm{Res}_{S\pr_i/\Spec V}\bG_{m}$ for  finite \'etale connected $V$-schemes $S_{i}\pr$, so \cite{SGA3IIInew}*{XIX, 8.4} gives an isomorphism $H^{1}_{\et}(V,P)\cong \prod_{S_{i}\pr}H^{1}_{\et}(S_{i}\pr, \bG_m)$. 
      The Grothendieck--Hilbert's 90 \cite{SGA3II}*{VIII, 4.5} identifies $H^1_{\et}(S_i\pr,\bG_{m})\simeq H^1_{\Zar}(S_i\pr, \bG_m)$, which are trivial by \cite{Bou98}*{II, \S5, no.~3, Prop.~5}. 
      So we have 
      \[
      \x{$H^1_{\et}(V,P)=\{\ast\}$\q for every quasitrivial $V$-torus $P$.}
      \]
      \benumr
      \item First, we reduce to the case for flasque tori. By the short exact sequence \cite{CTS87}*{1.3.2}
      \[
      \x{$1\ra M\ra F\ra P\ra 1$,}
      \]
      where $F$ is flasque and $P$ is quasitrivial, we obtain the commutative diagram with exact rows
      \[
      \begin{tikzcd}
        H^1_{\fpqc}(V,P)\ar{r}  & H^2_{\fpqc}(V,M) \ar{r} \ar{d} & H^2_{\fpqc}(V,F)  \ar{d} \\
          & H^2_{\fpqc}(\Frac V,M) \ar{r} & H^2_{\fpqc}(\Frac V,F),
      \end{tikzcd}
      \]
      where $H^{1}_{\fpqc}(V,P)=H^{1}_{\et}(V,P)=\{\ast\}$.
      Hence, it suffices to prove the assertion for the flasque $F$.

      Next, we induct on the rank $n$ of $V$.
      The case of $V=\Frac V$ is trivial, so when $n\geq 1$, for the prime $\mathfrak{p}$ of $V$ of height $n-1$, we assume that the assertion holds for $V_{\mathfrak{p}}$ (which has rank $n-1$).
      Denote $X=\Spec V$ and $Z=\Spec(V/\fm_V)$. 
      By \cite{SGA4II}*{V, 6.5}, we have the long exact sequence:
      \begin{equation}\label{loc-seq}
        \x{$\cdots \ra H^2_{Z}(X, F)\ra H^2_{\fpqc}(X,F)\ra H^2_{\fpqc}(X-Z,F)\ra H^3_{Z}(X,F)\ra \cdots$.}
      \end{equation}
      We conclude by the induction hypothesis and $H^2_{Z}(X,F)=0$ proved in \Cref{2-vanish}.      
      \item We first reduce to the case when $M$ is a torus. The isotriviality of $M$ yields a short exact sequence 
      \[
      1\ra T\ra M\ra \mu\ra 1, 
      \]
      where $T$ is a $V$-torus and $\mu$ is a finite multiplicative type $V$-group. 
      For the commutative diagram
\[
      \begin{tikzcd}
      \mu(V) \ar{r} \ar{d}  & H^1_{\fpqc}(V,T)\ar{r} \ar{d} & H^1_{\fpqc}(V,M) \ar{r} \ar{d} & H^1_{\fpqc}(V,\mu)  \ar{d} \\
      \mu(\Frac V) \ar{r} & H^1_{\fpqc}(\Frac V,T)\ar{r}  & H^1_{\fpqc}(\Frac V,M) \ar{r} & H^1_{\fpqc}(\Frac V,\mu)
      \end{tikzcd}
      \]
      with exact rows, the valuative criterion for properness of $\mu$ leads to $\mu(V)=\mu(\Frac V)$ and the injectivity of $H^1_{\fpqc}(V,\mu)\hra H^1_{\fpqc}(\Frac V,\mu)$. 
      Thus, a diagram chase reduces us to showing that 
      \[
      \x{$H^1_{\et}(V,T)\ra H^1_{\et}(\Frac V,T)$\qq is injective.}
      \]
     
      A flasque resolution of $T$ as (\ref{flasque-resolution}) leads to the following commutative diagram with exact rows
      \[
      \begin{tikzcd}
        H^1_{\et}(V,P)\ar{r}  & H^1_{\et}(V,T) \ar{r} \ar{d} & H^2_{\et}(V,F)  \ar{d} \\
          & H^1_{\et}(\Frac V,T) \ar{r} & H^2_{\et}(\Frac V,F),
      \end{tikzcd}
      \]
      where $H^{1}_{\et}(V,P)=\{\ast\}$.
      Since the map $H^2_{\et}(V,F)\hra H^2_{\et}(\Frac V,F)$ is injective by \ref{mult-brauer-i},  the map
      \[
      \x{\qq $H^1_{\et}(V,T)\hra H^1_{\et }(\Frac V,T)$\qq is injective.}\qedhere
      \]
      \eenum
      \epf
      \bcor
      For a flasque torus $F$ over a valuation ring $V$ with fraction field $K$, the map
      \[
      \x{\qqq $H^1_{\et}(V,F)\isoto H^1_{\et}(K,F)$ \qq is an isomorphism.}
      \]
      \ecor
      \bpf
      The injectivity follows from \Cref{mult-brauer}~\ref{mult-brauer-ii}. A limit argument reduces us to the case when $V$ has finite rank, then we iteratively use \Cref{2-vanish} with the following exact sequence (\emph{cf.}~\ref{loc-seq}) 
      \[
         H^1_{\et}(V,F)\ra H^1_{\et}(\Spec V\backslash \{\fm_V\},F)\ra H^2_{\fm_V}(V,F)=0,
      \]
      to reduce the rank of valuation rings by removing closed points, so the surjectivity follows. 
      \epf
      
      \section{Algebraizations and a Harder-type approximation}\label{wk-ap}
      The upshot of this section is \Cref{open-normal}, a higher-height analogue of Harder's weak approximation \cite{Har68}*{Satz.~2.1} to reduce \Cref{GSVal} to the case of Henselian rank-one valuation rings. 
      To prove this, we take advantage of techniques of algebraization from \cite{BC22}*{\S2} and Conrad's topologization of points.
      \bpp[Topologizing $R$-points of schemes]\label{topologizing}
      For a topological ring $R$ and an $R$-scheme (or  $R$-algebraic stack) $X$, 
      the problem of topologizing $X(R)$ functorially in $X$ compatible with the topology of $R$ has been  studied in recent years. 
      Precisely, we expect a topology on $X(R)$ satisfying some of the following
      \benumr
      \item\label{top-i} each $R$-morphism $X\ra X\pr$ induces a continuous map $X(R)\ra X\pr(R)$;
      \item\label{top-ii} for every integer $n\geq 0$, we have a canonical homeomorphism $\bA^n(R)\simeq R^n$;
      \item\label{top-iii} each closed immersion $X\hra X\pr$ induces an embedding $X(R)\hra X\pr(R)$; 
      \item\label{top-iv} each open immersion $X\hra X\pr$ induces an open embedding $X(R)\hra X\pr(R)$; and
      \item\label{top-v} for all $R$-morphisms $X\pr\ra X\la X\prpr$ of $R$-schemes, the identifications \[
      \text{\qq$(X\pr\times_X X\prpr)(R)=X\pr(R)\times_{X(R)}X\prpr(R)$\qq are homeomorphisms.}
      \]
      \eenum
      For all affine schemes $X$ of finite type over $R$, Conrad proved \cite{Con12}*{Prop.~2.1} that there is a \textbf{unique} way to topologize $X(R)$ such that \ref{top-i}--\ref{top-iii} and \ref{top-v} are satisfied. 
      Such topologization is realized by taking a closed immersion $X\hra \bA^n_R$ and endowing $X(R)$ with the subspace topology from $R^n$. 
      The resulting topology is not dependent on the choice of embeddings.
      For schemes $X$ locally of finite type over $R$, topologizing $X(R)$ is reduced to the affine case by patching open affine subschemes of $X$, which calls for several extra constraints on $R$.
      Namely, under the assumption that $R$ is local and $R^{\times}\subset R$ is open with continuous inversion (\emph{e.g.}, Hausdorff topological fields and arbitrary valuation rings with valuation topology), Conrad showed \cite{Con12}*{Prop.~3.1} that there is a \textbf{unique} way to topologize $X(R)$ satisfying \ref{top-i}--\ref{top-v} for all schemes $X$ locally of finite type over $R$.
      Subsequently, {\v{C}}esnavi{\v{c}}ius generalized Conrad's result to algebraic stacks (\emph{cf.}~\cite{MB01}*{Section.~2} for the case of Hausdorff topological fields).
      Without the local assumption, if $R^{\times}\subset R$ is open with continuous inversion, then $X(R)$ can be topologized for (ind-)quasi-affine or (sub)projective $R$-schemes $X$, see \cite{BC22}*{\S2.2.7}. 
      Note that all aforementioned results are generalizations of Conrad's version,  
      hence they are compatible when restricting the families of $X$ or of $R$.
      Since we only consider schemes, our topologization \emph{only} involves the following formation of Conrad.
      \epp
      \blem[\cite{Con12}*{Prop.~3.1}]\label{con12} 
      Let $R$ be a local topological ring such that $R^{\times}\subset R$ is open with  continuous inversion.
      There is a unique way to topologize $X(R)$ satisfying \ref{top-i}--\ref{top-v} for all schemes $X$ locally of finite type over $R$.
      Moreover, if $R$ is Hausdorff and $X$ is $R$-separated, then $X(R)$ is Hausdorff.
      \elem
      \blem[\cite{Con12}*{Ex.~2.2}]\label{con12-2.2} 
      For any continuous map $R\pr\ra R$ of topological rings and any affine scheme $X$ of finite type over $R$, the natural homomorphism $X(R)\ra X(R\pr)$ is continuous.
      Moreover, if $R\pr\subset R$ is closed (resp., open) subring, then $X(R)\hra X(R\pr)$ is a closed (resp., open) embedding. 
      \elem
      \bd 
      For a topological ring $R$ and a scheme $X$ locally of finite type over $R$, if $X(R)$ can be topologized as in \S\ref{topologizing}, then we say that $X(R)$ has a topology \emph{induced from} $R$.
      In particular, if there is an ideal $I\subset R$ such that the topology on $R$ is $I$-adic, then the induced topology on $X(R)$ is called $I$\emph{-adic}. 
      \ed
      Now, we apply Conrad's formation to our case when $R$ is a valued field. 
      Recall \S\ref{valuation topologies} and \Cref{top-prop} that for every valued field $(K,\nu)$, there is a valuation topology determined by $\nu$ and it is Hausdorff.
      By \Cref{nonarch}, a valued field $(K,\nu)$ is nonarchimedean, if the valuation topology on $K$ is induced by a nontrivial rank-one valuation, or equivalently, the valuation ring $V(\nu)$ of $K$ has a prime of height one.
      \blem\label{ces-par}
      Let $(K,\nu)$ be a valued field and let $X$ be a scheme locally of finite type over $K$.
      \benumr 
      \item\label{top-field} The set $X(K)$ has a topology induced from the valuation topology on $K$.
      \item\label{sep-Haus} If $X$ is separated over $K$, then $X(K)$ is Hausdorff for the valuation topology.
      \item\label{cl-open} For the valuation ring $V\subset K$ and an affine finite type $V$-scheme $Y$, the natural map $Y(V)\hra Y(K)$ is a closed and open embedding for the valuation topology.
      \item\label{dense} If $K$ is Henselian nonarchimedean and $X$ is $K$-smooth, then for the completion $\wh{K}$ of $K$ and the topologies on $X(K)$ and on $X(\wh{K})$ induced from $K$ and $\wh{K}$ respectively, the following map
      \[
\text{$X(K)\ra X(\wh{K})$ \q has dense image.}
      \]
      \eenum
      \elem
      \bpf
      For \ref{top-field} and \ref{sep-Haus}, note that by \Cref{top-prop}, $K$ is Hausdorff so $K^{\times}\subset K$ is open.
      It is clear that the inversion on $K^{\times}$ is continuous for the subspace topology. 
      It suffices to use \Cref{con12} to topologize $X(K)$; moreover, if $X$ is separated over $K$, then $X(K)$ is Hausdorff for the valuation topology.
      The assertion \ref{cl-open} follows from \Cref{con12-2.2} and \Cref{top-prop}  that the ball $V\subset K$ is closed and open.

      For \ref{dense}, we recall \S\ref{comparison} that the topology on $K$ is indeed $a$-adic for an $a\in V$ such that $\sqrt{(a)}$ is of height one.
      Thus $\wh{K}$ is the $a$-adic completion $\hka$.
      We then apply \cite{BC22}*{2.2.10~(iii)} and check the conditions: 
      \begin{itemize}
      \item[-] Let the topological ring $B$ be $K$ with $a$-adic topology. Then $\wh{B}=\hka$ and $(\hka)^{\times}\subset \hka$ is an open subring with continuous inversion.
      \item[-] Let the nonunital open subring $B\pr$ be the ideal $(a)$ of the valuation ring $V$. The induced topology on $(a)$ has an open neighborhood base of zero consisting of ideals $(a^n)_{n\geq 1}\subset (a)$ (\Cref{a-completion}~\ref{a-prime}).
      \item[-] The nonunital ring $(a)$ is Henselian in the sense of Gabber (\cite{BC22}*{2.2.1}), that is, every polynomial $f(T)=T^N(T-1)+a_NT^N+\cdots+a_1T+a_0$ where $a_i\in (a)$ and $N\geq 1$ has a (unique) root in $1+(a)$.
      Because $V$ is Henselian, by \SP{0DYD}, the pair $(V,(a))$ is also Henselian.
      Hence, Gabber's criterion shows that $(a)$ is Henselian, so the conditions in \cite{BC22}*{2.2.10~(iii)} are satisfied. \qedhere
      \end{itemize}
      \epf
      \blem\label{etale-open} 
      For a Henselian valued field $F$, 
      \benumr
      \item\label{et-op} every smooth morphism $f\colon X\ra Y$ between $F$-schemes locally of finite type induces an open map of topological spaces $f_{\mathrm{top}}\colon X(F)\ra Y(F)$;
      \item\label{smooth-kernel} for a monomorphism of $F$-flat locally finitely presented group schemes $G\pr\hra G$ where $G\pr$ is $F$-smooth, and the $F$-algebraic space $G\prpr\ce G/G\pr$, the map $G(F)\ra G\prpr(F)$ is open.
      \eenum
      \elem
      \bpf
      For \ref{et-op}, see \cite{GGMB14}*{3.1.4} and note that the `topological Henselianity' there yields the desired openness by \emph{loc.~cit.}, 3.1.2.
      For \ref{smooth-kernel}, see \cite{Ces15d}*{4.3~(a) and 2.8~(2)}, where $R$ is our $F$.
      \epf
      In addition to the topological properties above, the following lemma will be used repeatedly in the sequel.
      \blem\label{open-product-closure} 
      For a topological group $G$, an open subgroup $H\subset G$, and a subset $S\subset G$, we have $$S\cdot H=\ov{S}\cdot H.$$
      \elem 
      \bpf 
      Since $\ov{S}\cdot H\subset \ov{S\cdot H}$, it suffices to see that $S\cdot H=\ov{S\cdot H}$.
      The subset $G\backslash (S\cdot H)$ is a union of $g_iH$ for some $g_i\in G$, hence is open. In particular, $S\cdot H$ is closed, so the assertion follows.
      \epf

      \bpp[Regular sections, Cartan subalgebras and subgroups of type (C)]\label{C}
      Let $R$ be a ring and let $\fh$ be a Lie algebra over $R$ as a locally free module of rank $n$. 
      The Lie algebra structure (Lie bracket) is a morphism $A\colon \fh\ra \End_R(\fh)$.
      For any $R$-algebra $R\pr$, the $i$-th coefficient of the characteristic polynomial of degree $n$ for $B\in \End_{R\pr}(\fh_{R\pr})$ is of the form $(-1)^{n-i}\mathrm{Tr}(\wedge^{n-i}B)$, so the $i$-th coefficient of the characteristic polynomial is a morphism $\End_R(\fh)^{\otimes i}\ra R$. 
      Composing $A^{\otimes i}$ with the last morphism, we get 
      \[c_i\colon \fh^{\otimes i}\ra R,\]
      hence $c_i\in (\fh^{\vee})^{\otimes i}\subset \GG(\un{\Sym}_R(\fh^{\vee}))$.
      We define the \emph{Killing polynomial} of $\fh$ as $P_{\fh}(t)\ce t^{n}+c_{1}t^{n-1}+\cdots+c_{n}\in \GG(\un{\Sym}_R(\fh^{\vee}))[t]$.
      By construction, the formation of Killing polynomials commutes with base change.
      When $R$ is a field $k$, the largest $r$ such that $P_{\fh}(t)$ is divisible by $t^r$ is the \emph{nilpotent rank} of $\fh$.
      The nilpotent rank of the Lie algebra of a reductive group scheme is \'etale-locally constant (see \cite{SGA3II}*{XV, 7.3} and \cite{SGA3IIInew}*{XXII, 5.1.2, 5.1.3}).
      Every $a\in \fh$ satisfying $c_{n-r}(a)\neq 0$ is called a \emph{regular element}. 
      Let $G$ be a reductive group scheme over a scheme $S$. 
      For the Lie algebra $\fg$ of $G$, if a subalgebra $\fd\subset \fg$ is Zariski-locally a direct summand such that its geometric fiber $\fd_{\ov{s}}$ at each $s\in S$ is nilpotent and equals to its own normalizer, then $\sigma$ is a \emph{Cartan subalgebra} of $\fg$ (\cite{SGA3II}*{XIV, 2.4}). 
      We say an $S$-subgroup $D\subset G$ is \emph{of type (C)}, if $D$ is $S$-smooth with connected fibers, and $\mathrm{Lie}(D)\subset \fg$ is a Cartan subalgebra.
      A section $\sigma$ of $\fg$ is a \emph{regular section}, if $\sigma$ is in a Cartan subalgebra such that $\sigma(s)\in \fg_s$ is a regular element for all $s\in S$.
      A section of $\fg$ with regular fibers is \emph{quasi-regular}, hence regular sections are quasi-regular.
      \bpp[Schemes of maximal tori]\label{max-tor}
      For a reductive group scheme $G$ defined over a scheme $S$, the functor
      \[
      \x{$\mtg\colon \Sch_{/S}\op\ra \Sets$,\qq $S\pr\mapsto \{\x{maximal tori of $G_{S\pr}$}\}$.}
      \]
      is representable by an $S$-affine smooth scheme (\cite{SGA3II}*{XIV, 6.1}).
      For an $S$-scheme $S\pr$ and a maximal torus $T\in \mtg(S\pr)$ of $G_{S\pr}$, by \cite{SGA3IIInew}*{XXII, 5.8.3}, the morphism defined by conjugating $T$
         \begin{equation}\label{conj-tor}
            \tst \x{$G_{S\pr}\ra \un{\mathrm{Tor}}(G_{S\pr}),\qq g\mapsto gTg^{-1}$}
         \end{equation}
         induces an isomorphism $G_{S\pr}/\un{\mathrm{Norm}}_{G_{S\pr}}(T)\cong \un{\mathrm{Tor}}(G_{S\pr})$.
         Here, $\un{\mathrm{Norm}}_{G_{S\pr}}(T)$ is an $S\pr$-smooth scheme (see \cite{SGA3II}*{XI, 2.4bis}).
         Now, we establish the following lifting property of $\mtg$.
            \epp      
      \epp
    \blem\label{extend-tor}
    Let $G$ be a reductive group scheme over a local ring $R$ with residue field $\kappa$ and $Z$ the center of $G$. If the cardinality of $\kappa$ is at least $\dim (G/Z)$, then the following map is surjective:
      \[
        \tst  \mtg(R)\surjects \mtg(\kappa).
      \]
      \elem
      \bpf
      An isomorphism \cite{SGA3II}*{XII, 4.7~c)} of schemes $\mtg\isom \un{\mathrm{Tor}}(G/Z)$ reduces us to the semisimple adjoint case, where the maximal tori of $G$ are exactly the subgroups of type (C) (\cite{SGA3II}*{XIV, 3.18}).
      These subgroups are bijectively assigned by $D\mapsto \mathrm{Lie}(D)$ to the Cartan subalgebras of $\fg\ce \mathrm{Lie}(G)$, see \emph{ibid.}, 3.9.
      It suffices to lift a Cartan subalgebra $\fc_{\kappa}\subset \fg_{\kappa}$ to that of $\fg$.
      Since $\sharp \kappa\geq \dim(G/Z)=\dim(G)$, by \cite{Bar67}*{Thm.~1}, $\fc_{\kappa}$ is of the form $\mathrm{Nil}(a_{\kappa})\ce \bigcup_{n}\ker(\mathrm{ad}(a_{\kappa}^n))$ for some $a_{\kappa}\in \fc_{\kappa}$. 
      Hence \cite{SGA3II}*{XIII, 5.7} implies that each $a_{\kappa}\in \fc_{\gK}$ is a regular element of $\fg_{\gK}$.
      We take a section $a$ of $\fg$ passing through $a_{\gK}$ and claim that $\cV\ce \{\x{$s\in \Spec R$ \,|\, $a_{s}\in \fg_{s}$ is regular}\}$ is an open subset of $\Spec R$.
      We may assume that $R$ is reduced.
      Since the nilpotent rank of $\fg$ is locally constant, the Killing polynomial of $\fg$ at every $s\in \Spec R$ is uniformly of the form $P_{\fg_s}(t)=t^r(t^{n-r}+(c_1)_st^{n-r-1}+\cdots+(c_{n-r})_{s})$ such that $(c_{n-r})_s$ is nonzero.
      Thus, the regular locus in $\fg$ is the principle open subset $\{c_{n-r}\neq 0\}\subset \b{W}(\fg)$.
      The morphism $\b{W}(\fg)\ra \Spec R$ is flat, so $\cV\neq \emptyset$ is open, forcing that $\cV=\Spec R$.
      In particular, the regular elements $a_{\gK}\in \fc_{\gK}$ lifts to a quasi-regular section $a\in \fg$, which by \cite{SGA3IIInew}*{XIV, 3.7}, is regular.
      By definition of regular sections, there is a Cartan subalgebra of $\fg$ containing $a$ and is the desired lifting of $\fc_{\gK}$.
      \epf
      Next, we combine this lifting property with techniques of algebraization to deduce the density \Cref{alg}. 
      The next pages will deal with localizations, $a$-adic topology and completions of valuation rings.
      It is therefore recommended that readers refer to the Appendix \ref{apdx}, especially \S\ref{a-adic-top} and \Cref{a-completion}.
      \bpp[Rings of Cauchy sequences]\label{cs} 
      To the best of our knowledge, it is Gabber who first considered rings of Cauchy sequences (see also its generalization to Cauchy nets \cite{BC22}*{2.1.12}). 
      In this article, we take only one particular form to suit our need.
      Concretely, for a ring $A$ and a $t\in A$ such that $1+t\subset A^{\times}$, consider the truncated Cauchy sequences $(a_N)_{N\geq n}$ in $A[\f{1}{t}]$ for an $n\geq 0$. 
      With termwise addition and multiplication, all truncated Cauchy sequences form a ring $\mathrm{Cauchy}^{\geq n}(A[\f{1}{t}])$.  
      With this concept, one can translate the approximation process into certain operations on rings of Cauchy sequences and thus grasp the approximation properties through the algebrogeometric properties of the ring $\mathrm{Cauchy}^{\geq n}(A[\f{1}{t}])$.
      \epp
      \bpp[Setup]\label{setup}
      In the sequel, consider the subcase of \S\ref{cs}: let $A=V$ be a valuation ring of rank $n$ and let $t=a$ lie in $\fm_V\backslash \fp$ for the prime $\fp$ of height $n-1$.
      By \Cref{a-completion}, $V[\f{1}{a}]$ and the $a$-adic completion $\hva$ are valuation rings of ranks $n-1$ and $1$ respectively, and the $a$-adic completion $\vaih$ of $\via$ is $\hka\ce \Frac \hva$.
      By \Cref{completion-NA} and \Cref{hens-val-field}, $\hka$ is  nonarchimedean and $\hva$ is a \emph{Henselian local ring}.
      For every $\hka$-scheme $X$ locally of finite type, we will endow $X(\hka)$ with the $a$-adic topology.
      \epp
      \blem\label{local-ring} 
      For the setup \S\ref{setup}, the $\varinjlim_{m\geq 0}\mathrm{Cauchy}^{\geq m}(\via)$ is a local ring with residue field $\hka$.
      \elem
      \bpf
      Taking $a$-adic completion of $\via$ yields the following surjection map
      \[
      \textstyle \cA\ce \varinjlim_{m\geq 0}\mathrm{Cauchy}^{\geq m}(\via)\surjects \hka,
      \]
      whose kernel is denoted by $I$.
      For any sequence $(b_N)_N\in I$, its tail lies in $\mathrm{Im}(a^mV\ra \via)$ for all $m>0$, so the tail of $(1+b_N)_N$ consists of units in $V$ that lie in $\mathrm{Im}((1+a^mV)\ra \via)$.
      Since $\via$ is local, the tail of $(1+b_N)_N$ is termwise invertible in $\via$ and the inverses form a Cauchy sequence. 
      Since $I\subset \cA$ is an ideal such that  $\cA/I$ is a field and $1+I$ is invertible,  $\cA$ is a local ring with residue field $\hka$.
      \epf

      \beg\label{a-open} 
      Consider the setup \S\ref{setup}.
      Then  \Cref{top-prop} implies that $\hva\subset \hka$ is open and closed.
      Let $G$ be a reductive $V$-group scheme and recall $\mtg$ (\S\ref{max-tor}).
      By \Cref{ces-par}~\ref{cl-open}, the subsets $G(\hva)\subset G(\hka)$ and $\mtg(\hva)\subset \mtg(\hka)$ are $a$-adically open and closed.
      \eeg

      \blem\label{alg}
      Consider the setup \S\ref{setup}.
      For a reductive $V$-group scheme $G$,
      \[
      \tst \x{the image of\q $\mtg(\via)\ra \mtg(\hka)$\q is $a$-adically dense.}
      \]
      \bpf
      As shown in \Cref{local-ring}, the ring $\varinjlim_{m\geq 0}\mathrm{Cauchy}^{\geq m}(\via)$ is local with residue field $\hka$.
      Since $\mtg$ is finitely presented and affine over $\via$, the lifting \Cref{extend-tor} leads to a surjection below
      \[
      \tst \qq\varinjlim_{m\geq 0}\bigl(\mtg\bigl(\mathrm{Cauchy}^{\geq m}(\via)\bigr)\bigr)\isom \mtg\bigl(\varinjlim_{m\geq 0}\bigl(\mathrm{Cauchy}^{\geq m}(\via)\bigr)\bigr)\surjects \mtg(\hka).
      \]
      Due to this surjection, all elements in $\mtg(\hka)$ are limits of Cauchy sequences in $\mtg(V[\f{1}{a}])$, hence the image of the map $\mtg(\via)\ra \mtg(\hka)$ is $a$-adically dense in $\mtg({\hka})$.
      \epf
      \elem
      Roughly speaking, this density permits us to ``replace'' maximal tori of $G_{\nhka}$ by those of $G_{\via}$.
      Next, we obtain openness of certain maps, then take images to construct an open normal subgroup of $G(\hka)$ contained in the closure of the image of $G(\via)\ra G(\hka)$.
      First, recall some criteria for openness.
      
      \blem\label{via-torus}
      Consider the setup \S\ref{setup}. 
      Let $T$ be a torus over $\via$.
      \benumr
      \item\label{comp-prod} There is a minimal Galois cover $R$ of $\via$ splitting $T$ (see \S\ref{mult}), and we have isomorphisms
      \[
      \tst R\otimes_{\via}\hka\simeq \hRa\simeq \prod_{i=1}^rL_i,
      \]
      where $\hRa$ is the $a$-adic completion of $R$ for the topology induced from $\via$. 
      Each $L_i/\hka$ is a minimal Galois extension splitting $T_{\nhka}$ and is $a$-adically complete;
      in particular, any minimal Galois extension $L_0/K$ splitting $T_{\nhka}$ is isomorphic to $L_i$ for all $i$, that is, $L_0\simeq L_i\simeq L_j$ for $i\neq j$.
      \item\label{torus-open-first} For a minimal Galois field extension $L_0/\hka$ splitting $T_{\nhka}$, the image $U$ of the  norm map 
      \[
      N_{L_0/\nhka}\colon T(L_0)\ra T(\hka)
      \]
       is $a$-adically open in $T(\hka)$ and contained in the closure $\ov{T(\via)}$ of $\mathrm{Im}(T(\via)\ra T(\hka))$.
      \eenum
      \elem
      \bpf\hfill
      \benumr
      \item 
      The existence of a minimal Galois cover $R/V[\f{1}{a}]$ splitting $T$ follows from \Cref{isotrivial}.
      Since $R$ is a finite flat $\via$-module, it is free and we have $\hRa\simeq R\otimes_{\via}\hka\simeq \prod_{i=1}^rL_i$, where $L_i$ are $a$-adically complete fields.
      By \Cref{bc-galois} and \S\ref{mult} we conclude.
      \item First, we prove that $U$ is $a$-adically open.
      For the norm map $\mathrm{Res}_{L_0/\nhka}(T_{L_0})\ra T_{\nhka}$, its kernel $\cT$ is a torus: after some base change, $T_{\nhka}$ splits as $\bG_m^k$, so the associated $\b{Z}$-module of the corresponding base change of $\cT$ is the following $\b{Z}$-lattice with a trivial Galois action
          \[
          \tst \x{$\mathrm{Coker}\bigl(\b{Z}^k\ra \b{Z}[\Gal(L_0/\hka)]^k, (n_i)\mapsto (n_i\cdot \id)\bigr)\isom \b{Z}[\Gal(L_0/\hka)-\{\id\}]^k$}.
          \]
      So, by \cite{SGA3II}*{IX, 2.1~e)}, as a torus, the kernel $\cT$ is $\hka$-smooth.
      By \Cref{etale-open}~\ref{smooth-kernel}, the map 
      \[
    \tst \qqq \text{$N_{L_0/\nhka}\colon T(L_0)\ra T(\hka)$,\q i.e. \q $(\mathrm{Res}_{L_0/\nhka}T_{L_0})(\hka)\ra \bigl((\mathrm{Res}_{L_0/\nhka}T_{L_0})/\cT\bigr)(\hka)$}
      \]
      is $a$-adically open so the image $U= N_{L_0/\nhka}(T(L_0))\subset T(\hka)$ is $a$-adically open.

      Next, we prove that $U\subset \ov{T(\via)}$.
      The isomorphism $\hRa\cong \prod_{i=1}^{r}L_i$ obtained in \ref{comp-prod} implies that the image of $R^{\times}\ra \prod_{i=1}^{r}L_i^{\times}$ is $a$-adically dense.
      As $T_{R}$ is split, the image of the composite
      \[
     \tst T(R)\ra \prod_{j=1}^{r}T(L_j)\overset{\mathrm{pr}_1}{\ra} T(L_1)\cong T(L_0)
      \]
      is $a$-adically dense. 
      Composing this with $N_{L_0/\nhka}$, we see that $T(R)$ has dense image in $U=N_{L_0/\nhka}(T(L_0))$.
      The composite $T(R)\ra T(L_0) \ra T(\hka)$ factors through the norm map $N_{R/\via}\colon T(R)\ra T(\via)$, so the image of $T(\via)$ is dense in $U$, that is $\tst U\subset \ov{T(\via)}$.\qedhere
      \eenum
      \epf
      Subsequently, we approximate the $\hka$-points of a maximal torus of $G_{\nhka}$ by using $\via$-points.
      \blem\label{open-in-torus}
      Consider the setup \S\ref{setup}.
      For a reductive $V$-group scheme $G$, the closure $\ov{G(\via)}$ of the image of $G(\via)\ra G(\hka)$, a maximal torus $T$ of $G_{\wh{K}^a}$ with minimal splitting field $L_0$, and
      \[
      \x{the norm map \q $N_{L_0/\wh{K}^a}\colon T(L_0)\ra T(\hka)$,}
      \]
      the image $U=N_{L_0/\nhka}(T(L_0))$ is an $a$-adically open subgroup of $T(\hka)$ and is contained in $\ov{G(\via)}$.
      \elem
      \bpf 
      The $a$-adically open aspect of the assertion follows from \Cref{via-torus}~\ref{torus-open-first} because the arguments there, by base change, apply to all $\hka$-tori as well.
      The proof for $U\subset \ov{G(\via)}$ proceeds as follows.
      \benumr
      \item 
         Since $\hka$ is Henselian, by a criterion for openness \Cref{etale-open}~\ref{smooth-kernel}, the following map from (\ref{conj-tor})
         \[
         \tst \x{$\phi\colon G(\hka)\ra \mtg(\hka),\qq g\mapsto gTg^{-1}$\q is $a$-adically open.}
         \]
         Consequently, $\phi$ sends every $a$-adically open neighborhood $W$ of $\id\in G(\hka)$ to an $a$-adically open neighborhood of $T$.
         The density \Cref{alg} of $\mtg(\via)$ in $\mtg(\hka)$  implies that
         \[
         \tst \phi(W)\cap \mathrm{Im}(\mtg(\via)\ra \mtg(\hka))\neq \emptyset.
         \]
         Hence, there are a torus $T\pr\in \mtg(\via)$ and a $g\in W$ such that $gTg^{-1}=T\pr_{\nhka}\in \phi(W)$.
      \item 
      For any $u\in U$, the map $\sigma_u\colon G(\hka)\ra G(\hka)$ defined by $g\mapsto g^{-1}ug$ is continuous. 
      Let $W\ce \sigma_u^{-1}(U)$.
      By the construction in (i), there are a $w\in W$ and a torus $T\pr\in \mtg(\via)$ such that $w T w^{-1}=T\pr_{\nhka}$.
      Note that $u\in w U w^{-1}=\gG N_{L_0/\nhka}(T(L_0))\gG^{-1}$, which by transport of structure, is equal to $N_{L_0/\nhka}(T\pr_{\nhka}(L_0))$.
      By \Cref{via-torus}, the last term is contained in $\ov{\im(T\pr(\via)\ra T\pr(\hka))}$, so is contained in $\ov{G(\via)}$. \qedhere
          \eenum
      \epf
      \bcor\label{lift-tor-int} 
     Consider the setup \S\ref{setup} and a reductive $V$-group scheme $G$, we have 
     \[
     \mathrm{Im}\bigl(\mtg(\hva)\ra \mtg(\hka)\bigr)\subset \ov{\mathrm{Im}\bigl(\mtg(V)\ra \mtg(\hka)\bigr)}.
     \]
     More precisely, for every maximal torus $T$ of $G_{\nhva}$ and every $a$-adically open neighborhood $W$ of $\mathrm{id}\in G(\hka)$, there exist a maximal torus $T_0$ of $G$ and a $g\in W$ such that $(T_0)_{\nhka}=gT_{\nhka}g^{-1}$.
     \ecor
     \bpf
     By the argument (i) for \Cref{open-in-torus},  $\phi(W)\cap \mtg(\hva)$ is an $a$-adically open neighborhood of $T_{\nhka}\in \mtg(\hka)$. 
     Since $V\simeq \via\times_{\nhka}\hva$ (\Cref{a-completion}~\ref{fp-a-comp}) and $\mtg$ is affine, we get
     \[
     \tst \mtg(V)\isoto \mtg(\via)\times_{\mtg(\nhka)}\mtg(\hva).
     \]
     By \Cref{alg}, the image of $\mtg(\via)\ra \mtg(\hka)$ is $a$-adically dense, so we have
     \[
     \tst \phi(W)\cap \mtg(\hva)\cap \mathrm{Im}(\mtg(\via))\neq \emptyset,
     \]
     giving a maximal torus $T_0\in \mtg(V)$ and $g\in W$ such that $(T_0)_{\nhka}=gT_{\nhka}g^{-1}\in \phi(W)$.
     \epf
      Next, we prove \Cref{open-normal} by constructing an open subgroup in the closure of $G(\via)$.
      By lumping together the approximations in toral cases (\Cref{open-in-torus}),  the resulting open subgroup is normal.
      This normality is crucial for the dynamic argument for root groups for the product formula \Cref{decomp-gp}.
      \bprop\label{open-normal}
      Consider the setup \S\ref{setup}.
      For a reductive $V$-group scheme $G$, the closure $\ov{G(\via)}$ of the image of $G(\via)\ra G(\hka)$ contains an $a$-adically open normal subgroup $N$ of $G(\hka)$.
      \eprop
      \bpf In the proof, all open subsets without the word `Zariski' refer to  $a$-adically open subsets.
      \benumr
      \item Fix a maximal torus $T\subset G_{\nhka}$. We denote by $\fg$ the Lie algebra of $G_{\nhka}$ and by $\fh$ the Lie algebra of $T$. 
      For each $g\in G_{\nhka}$ and the subspace $\fg^{\mathrm{ad}(g)}\subset \fg$ fixed by $\mathrm{ad}(g)$, by \cite{SGA3II}*{XIII, 2.6~b)}, $\dim \fg^{\mathrm{ad}(g)}\geq \dim T$. 
      Let \emph{regular locus} $G\reg\subset G_{\nhka}$ be the subscheme of all $g\in G_{\nhka}$ that satisfy $\dim(\fg^{\mathrm{ad}(g)})=\dim T$. 
      By \cite{SGA3II}*{XIII, 2.7}, $G\reg$ is Zariski open.
      By the following equation 
      \[
      \dim(\fg^{\mathrm{ad}(g)})=\dim(\fh^{\mathrm{ad}(g)})+\dim((\fg/\fh)^{\mathrm{ad}(g)}),
      \]
      an element $t\in T$ is regular in $G_{\nhka}$ (namely, $t\in T^{\mathrm{reg}}\ce G^{\mathrm{reg}}\cap T$) if and only if $(\fg/\fh)^{\mathrm{ad}(t)}=0$. 
      \item Recall $L_0$ and the open subgroup $U\subset T(\hka)$ in \Cref{open-in-torus}, we claim that $U\cap T^{\mathrm{reg}}(\hka)\neq \emptyset$. 
      Consider the norm map $\mathrm{Nm}\colon \Res_{L_0/\nhka}(T_{L_0})\ra T$.
      Note that $T_{L_0}\isom \bG_{m,L_0}^k$ is isomorphic to a Zariski dense open subset of $\bA^k_{L_0}$, so $\Res_{L_0/\nhka}(T_{L_0})$  is also a Zariski dense open subset of $\bA^{mk}_{\nhka}$ for $m\ce[L_0:\hka]$.
      The field $\hka$ is infinite, so we have $\bigl({\Res_{L_0/\nhka}(T_{L_0})}\bigr)(\hka)\cap \mathrm{Nm}^{-1}(T^{\mathrm{reg}})(\hka)\neq \emptyset$.
      Applying $\mathrm{Nm}$ to this nonempty intersection, we proved our claim that $U\cap T^{\mathrm{reg}}(\hka)\neq \emptyset$.
      \item For a fixed $t_{0}\in U\cap T\reg(\hka)$, by (i), we have $(\fg/\fh)^{\mathrm{ad}(t_0)}=0$. 
      So \cite{SGA3II}*{XIII, 2.2} implies that
      \[
      f\colon G_{\nhka}\times T\ra G_{\nhka}, \qq (g,t)\mapsto gtg^{-1}
      \]
      is smooth at $(\id, t_0)$.
      Thus, there is a Zariski open neighborhood $W$ of $(\id,t_0)$ such that $f|_W\colon W\ra G_{\nhka}$ is smooth. 
      By \Cref{etale-open}~\ref{et-op}, $W(\hka)\ra G(\hka)$ is open.
      Thus the open neighborhood $W\pr\ce W(\hka)\cap (G(\hka)\times U)$ of $(\id, t_0)$ has open image under $f_{\mathrm{top}}$.
      The $G_{\nhka}$-translations $\tau_{h}\colon (g,t)\mapsto (hg,t)$ for $h\in G_{\nhka}$ induce automorphisms of $G_{\nhka}\times T$, so $f$ is also smooth at $(h,t_0)$. 
      Similar to above, all $G(\hka)$-translations of $W\pr$ have open images under $f_{\mathrm{top}}$.
      So, there is an open subset $U_0\subset U$ such that $E\ce f(G(\hka)\times U_0)$ is open.
      Let $N$ be the subgroup of $G(\hka)$ generated by $E$.
      The openness of $E$ implies that $N$ is an open subgroup of $G(\hka)$. 
      \item As $E$ is stable under $G(\hka)$-conjugation, $N$ is normal in $G(\hka)$.
      For each $g\in G(\hka)$, we denote $T^{g}\ce gTg^{-1}$.
      Then $U^{g}\ce N_{L_0/\nhka}(T^g(L_0))$ satisfies $U^{g}=gUg^{-1}$. 
      \Cref{open-in-torus} applies to $T^{g}$ and gives $U^{g}\subset \ov{G(\via)}$.
      Thus $E\subset \bigcup_{g\in G(\nhka)} U^g \subset \ov{G(\via)}$. Since $E$ generates $N$, we obtain 
\[ \tst N\subset \ov{G(\via)}. \qedhere \]
      \eenum
      \epf
    \bcor\label{clopen}
     With the notations in \Cref{open-normal},  $\ov{G(\via)}$ is an open subgroup of $G(\hka)$ and
     \[
     \tst \ov{G(\via)}\cdot G(\hva)=\mathrm{Im}\bigl(G(\via)\ra G(\hka)\bigr)\cdot G(\hva).
     \]
     \ecor
     \bpf
       The image of $G(\via)\ra G(\hka)$ is a subgroup of $G(\hka)$, hence so is its closure $\ov{G(\via)}$.
       Since $\ov{G(\via)}$ contains the open subset $N$, it is an open subgroup of $G(\hka)$.
       Recall \Cref{a-open} that the subgroup $G(\hva)\subset G(\hka)$ is open and closed. 
       By \Cref{open-product-closure}, the desired equation follows.
     \epf
      \section{Passage to the Henselian rank one case: patching by a product formula}\label{passage-h}
      The aim of this section is to reduce \Cref{GSVal} to the case when $V$ is a Henselian valuation ring of rank one.
      The key of our reduction \Cref{rank-one-kernel-trivial} is the product formula \Cref{decomp-gp} for patching torsors:
      \[
          \tst  \x{$G(\hka)=\mathrm{Im}\bigl(G(\via)\ra G(\hka)\bigr)\cdot G(\hva).$ }
      \]
      To show this product formula, we use the Harder-type weak approximation \Cref{open-normal}.

      First, we recall a criterion for anisotropicity \cite{SGA3IIInew}*{XXVI, 6.14}, which is practically useful.
      \blem\label{ani-cri}
        A reductive group scheme $G$ over a semilocal \emph{connected} scheme $S$ is anisotropic if and only if $G$ has no proper parabolic subgroup and $\rad(G)$ contains no copy of $\bG_{m,S}$.
      \elem
      Precisely, to determine whether $G$ is anisotropic, we consider the functor parametrizing parabolic subgroups
      \[
      \Par(G)\colon \Sch_{/S}\op\ra \Sets,\qq S\pr\mapsto \{\x{parabolic subgroups of $G_{S\pr}$}\},
      \]
      which is representable by a smooth projective $S$-scheme (see \cite{SGA3IIInew}*{XXVI, 3.5})\footnote{For the formation of $\Par(G)$, the base scheme $S$ does not have to be connected.}.
      Note that $G$ is also an element in $\Par(G)(S)$; we denote this non-proper parabolic subgroup by $\ast\in \Par(G)(S)$.
      
      Recall \Cref{nonarch} and \S\ref{comparison} that a valued field $K$ is \emph{nonarchimedean} if its valuation ring $V$ has a height-one prime ideal $\fp_1$.
      The completion $\wh{K}$ equals the $a$-adic completion $\hka$ of $K$ for an $a\in \fp_1\backslash\{0\}$.

      \blem\label{parden}
      For a Henselian nonarchimedean valued field $K$ with its completions $\wh{K}$, a reductive $V$-group scheme $G$, and the valuation topology on $\Par(G)(\wh{K})$ induced from $\wh{K}$,
      \benumr
      \item\label{par-dense} the image of $\Par(G)(K)\ra \Par(G)(\wh{K})$ is dense;
      \item\label{par-nontri} let $V\subset K$ and $\wh{V}\subset \wh{K}$ be the valuation rings, if $\Par(G)(\wh{V})\neq \{\ast\}$, then $\Par(G)(V)\neq \{\ast\}$.  
      \eenum 
      \elem 
      \bpf
      The assertion \ref{par-dense} follows from \Cref{ces-par}~\ref{dense}.  
      If $\Par(G)(\wh{V})\neq \{\ast\}$, then the valuative criterion for the separatedness of $\Par(G)$ implies that $\Par(G)(\wh{K})$ contains an $x\neq \ast$. 
      By \Cref{ces-par}~\ref{sep-Haus}, $\Par(G)(\wh{K})$ is Hausdorff so $x$ has an open neighborhood $U_x$ that excludes $\ast$. 
      The density of the image of $\Par(G)(K)\ra \Par(G)(\wh{K})$ shown in \ref{par-dense} yields an $y\in \Par(G)(K)$ whose image is contained in $U_x$. 
      Therefore, $y\neq \ast$ and $\Par(G)(K)\neq \{\ast\}$.
      By the valuative criterion for the properness of $\Par(G)$ over $V$, we conclude.
      \epf 

      The following \Cref{aniso-int-rat} generalizes \cite{Pra82}*{Theorem~(BTR)} to valuation rings of higher rank. 
    For a reductive group scheme $H$ over a scheme $S$, the $S$-\emph{split rank} of $G$ is the largest $k$ such that $\bG_{m,S}^k\subset G$.
     In particular, for any $S$-scheme $S\pr$, the $H_{S\pr}$ is anisotropic if and only if it has zero $S\pr$-split rank.

      \bprop\label{aniso-int-rat}
      Let $G$ be a reductive group scheme over a valuation ring $V$ with fraction field $K$.
      \benum
      \item\label{parabolic} A parabolic subgroup $P\subset G$ is minimal if and only if the parabolic subgroup $P_K\subset G_K$ is minimal.
      \item\label{ani-frac} The $V$-split rank of $G$ equals the $K$-split rank of $G_{K}$.
      \item \label{par-mini} If $K$ is Henselian nonarchimedean, then for the completion $\wh{V}$ of $V$ and a minimal parabolic subgroup $P\subset G$, the base change $P_{\wh{V}}$ is a minimal parabolic subgroup of $G_{\wh{V}}$.
       \item\label{hen-ani}
       If $K$ is  Henselian nonarchimedean, then for the completion $\wh{V}$ of $V$, 
\[
   \x{the $V$-split rank of $G$ equals the $\wh{V}$-split rank of $G_{\wh{V}}$.}
\]
  \item\label{ani-eq} If $K$ is Henselian and $V\neq K$, then $G$ is anisotropic if and only if $G(V)=G(K)$.
      \eenum
      \eprop
      \bpf \hfill
      \benum
      \item 
      If $P_K$ is minimal, then any minimal parabolic subgroup $Q$ of $G$ contained in $P$ satisfies $Q_K= P_K$.
      The valuative criterion for the separatedness of $\Par(G)$ over $V$ implies that $Q=P$, so $P$ is minimal.
      Now, we assume that $P\subset G$ is minimal.
      If there is a minimal parabolic subgroup $Q$ of $G_K$ contained in $P_K$, then the valuative criterion for the properness of $\Par(G)$ lifts $Q$ to a parabolic $\wt{Q}\subset G$, which must be minimal.
      Then, by \cite{SGA3IIInew}*{XXVI, 5.7~(ii)}, two minimal parabolics $\wt{Q}$ and $P$ are conjugated by an element of $G(V)$, which forces that $P_K=Q$ is minimal.  
      \item 
      When $G$ is a $V$-torus, we note that \Cref{lift-split-rank} suffices.
      In the general case, we reduce to this case of tori.
      Let $L$ be a Levi subgroup of a minimal parabolic $P\subset G$ and denote by $\rad(L)_{\mathrm{split}}$ the maximal $V$-split subtorus  of  $\rad(L)$.
      By \cite{SGA3IIInew}*{XXVI, 6.16}, the $V$-split rank of $G$ is equal to $\dim (\rad(L)_{\mathrm{split}})$.
      By \ref{parabolic}, $P_K$ is still a minimal parabolic subgroup of $G_K$ thereby \emph{op.~cit.} applies: the $K$-split rank of $G$ is equal to $\dim (\rad(L_K)_{\mathrm{split}})$.
      So we are reduced to the known toral case (\cite{SGA3IIInew}*{XXII, 4.3.6}) $\dim(\rad(L)_{\mathrm{split}})=\dim((\rad(L)_K)_{\mathrm{split}})$ for the $V$-torus $\rad(L)$.
      \item Let $L$ be a Levi subgroup of $P$, then $L_{\wh{V}}$ is a Levi subgroup of $P_{\wh{V}}$.
      By \cite{SGA3IIInew}*{XXVI, 1.20}, the set $\Par(L)(\wh{V})$ is the set of parabolics of $G_{\wh{V}}$ that are contained in $P_{\wh{V}}$ and $\Par(L)(V)$ is the set of parabolics of $G$ that are contained in $P$. 
      Hence, we conclude by \Cref{parden}~\ref{par-nontri}.
      \item
      For a Levi subgroup $L$ of a minimal parabolic subgroup $P$ of $G$, by \ref{par-mini},  $L_{\wh{V}}$ is a Levi subgroup of the minimal parabolic subgroup $P_{\wh{V}}$ of $G_{\wh{V}}$.
      Therefore, a similar argument in \ref{ani-frac} reduces us to the case when $G$ is a $V$-torus $T$.
      Taking the quotient of $T$ by its maximal split subtorus $T_{\mathrm{split}}$, we may assume that $T$ is anisotropic.
      Consider the following functor (\cite{SGA3II}*{X, 5.6})
      \[
      \un{X}^{\ast}(T)\colon \Sch\op_{/V}\ra \Sets, \qq R\mapsto \Hom_{\x{$R$-gr.}}(T_R,\mathbb{G}_{m,R}),
      \]
      which is representable by an \'etale locally constant group scheme.
      Since $T$ is isotrivial (\Cref{isotrivial}), by \cite{SGA3IIInew}*{XXVI, 6.6}, the property $\un{X}^{\ast}(T)(R)\neq 0$ is equivalent to that $T_R$ contains a copy of $\bG_{m,R}$.
      If $\un{X}^{\ast}(T)(\wh{V})\neq 0$, then by \Cref{a-completion}~\ref{comp-res},  the sets $\un{X}^{\ast}(T)(V/\fm_V)=\un{X}^{\ast}(T)(\wh{V}/\fm_{\wh{V}})$ contain nonzero elements.
      Since $V$ is Henselian and $\un{X}^{\ast}(T)$ is $V$-smooth, we have the surjection
      \[ 
        \x{$\un{X}^{\ast}(T)(V)\surjects \un{X}^{\ast}(T)(V/\fm_V)\neq 0$.}
       \]
      Thus $T$ contains a copy of $\bG_{m,V}$, which is in contradiction to the anisotropic assumption on $T$.  
      This contradiction shows that $\un{X}^{\ast}(T)(\wh{V})=0$, namely, $T_{\wh{V}}$ is also anisotropic, hence we conclude.
      \item
      If we have $G(K)=G(V)$, then it is impossible for $G$ to contain a $\bG_{m,V}$ because $K^{\times}=\bG_{m}(K)\subset G(K)$ strictly contains $V^{\times}=\bG_m(V)\subset G(V)$. 
      Therefore, $G$ is anisotropic.
      Now assume that $G$ is anisotropic and we show that $G(K)=G(V)$.
      By \cite{BM21}*{2.22}, $V$ is a filtered direct union of valuation subrings $V_{i}$ of finite rank, such that each $V_{i}\ra V$ is a local ring map.
      By \cite{EGAIV4}*{18.6.14~(ii)}, $V$ is a filtered direct union of Henselian valuation subrings $V_{i}^{\mathrm{h}}$ of finite rank.
      Similarly, $K$ is a filtered direct union of $K_i^{\mathrm{h}}\ce\Frac(V_i^{\mathrm{h}})$.
      Since $G$ is finitely presented over $V$, there is an index $i_{0}$ and an affine group scheme $G_{i_0}$ smooth and finitely presented (\cite{Nag66}*{Thm.~3'}) over $V_{i_{0}}^{\mathrm{h}}$ such that $G_{i_0}\times_{V_{i_0}}V\simeq G$.
      Further, by \cite{Con14}*{3.1.11}, $G_{i_0}$ and hence $(G_{i})_{i\geq i_0}$ are reductive group schemes.
      It is clear that all $(G_i)_{i\geq i_0}$ are anisotropic.
      By a limit argument \SP{01ZC}, we have $G(V)=\varinjlim_{i\geq i_0}G(V_i^{\mathrm{h}})$ and $G(K)=\varinjlim_{i\geq i_0}G(K_i^{\mathrm{h}})$.
      Subsequently, it remains to prove the case when $V$ is Henselian of finite rank $n$.

      First, we prove the case when $V$ is of rank one.
      For $a\in \fm_{V}\backslash\{0\}$, we form the $a$-adic completion $\hva$ of $V$ with $\hka\ce \Frac \hva$.
      By \ref{hen-ani}, $G_{\nhva}$ is anisotropic.
      For the nonarchimedean complete valued field $\hka$, by \cite{Mac17}*{Thm.~1.1}, $G(\hva)$ is a maximal  bounded\footnote{Recall from \cite{BrT2}*{1.7.3~f) or 4.2.19} (\emph{cf.} \cite{BLR90}*{Ch.~1, Def.~2}) that for a valued field $(K,\nu)$ and a $K$-scheme $X$, a subset $P\subset X(K)$ is \emph{bounded}, if for all $f\in K[X]$, we have $\mathrm{inf}_{x\in P}\nu(f(x))>-\infty$.
      For instance, the subset $\b{Z}_p\subset \b{Q}_p$ is bounded because $\nu(\b{Z}_p)\geq 0$; the subset $\{p^{-n}\}_{n\geq 1}$ is not bounded because $\nu(p^{-n})=-n$ tends to $-\infty$.} subgroup of $G(\hka)$.
      On the other hand, a result of Bruhat--Tits--Rousseau \cite{Rou77}*{Thm.~5.2.3} (or \cite{BrT2}*{p.~156, Rem.}) shows that $G(\hka)$ is bounded.
      Consequently, we have $G(\hva)=G(\hka)$.
      The rank-one assumption ensures that $V\hookrightarrow \hva$ is injective (\cite{FK18}*{Ch.~0, Thm.~9.1.1~(2)}), so the map $G(V)\hra G(\hva)$ is injective.
      The equality $V=K\times_{\nhka}\hva$ (\Cref{a-completion}~\ref{fp-a-comp}) and the affineness of $G$ yield a bijection
      \[
      \tst G(V)\isoto G(K)\times_{G(\nhka)}G(\hva)\cong G(K). 
      \]
      When $V$ is of rank $n>1$, we assume the assertion holds for the case of rank $\leq n-1$ and prove by induction.
      For the prime $\fp\subset V$ of height $n-1$, by \Cref{basic-one}~\ref{hens-preserve}, the localization $V_{\fp}$ and the quotient $V/\fp$ are Henselian valuation rings.  
      Due to \Cref{a-completion}~\ref{rank-val}, the rank of $V/\fp$ is one and the rank of $V_{\fp}$ is $n-1$.
      Since $V$ is Henselian, sections of $\Par(G)$ and $\un{X}^{\ast}(\rad(G))$ over $V/\fm_V$ lift over $V$.
      Hence, $G_{V/\fm_V}$ is anisotropic and so is $G_{V/\fp}$.
      As $G$ is anisotropic, by \ref{ani-frac}, so are $G_{K}$ and $G_{V_{\fp}}$. 
      By the settled rank-one case and the induction hypothesis, we have
      \begin{equation}\label{ani-one}
       \x{$G(V/\fp)=G(V_{\fp}/\fp V_{\fp})$\qq and \qq $G(V_{\fp})=G(K)$.}
      \end{equation}
      The affineness of $G$ and the isomorphism $V\isoto V_{\fp}\times_{V_{\fp}/\fp V_{\fp}}V/\fp$ lead to the isomorphism
      \begin{equation}\label{ani-fib-pro}
            G(V)\isoto G(V_{\fp})\times_{G(V_{\fp}/\fp V_{\fp})} G(V/\fp).
      \end{equation}
      Therefore,  the combination of (\ref{ani-fib-pro}) and  (\ref{ani-one}) gives us the desired equation $G(V)=G(K)$.    \qedhere
      \eenum
      \epf
      The following lemma provides the version for tori of the product formula.
      \blem\label{pd-tori}      
        For a valuation ring $V$ of rank $n>0$, the prime $\fp\subset V$ of height $n-1$, an element $a\in \fm_{V}\backslash \fp$, the $a$-adic completion $\hva$ with $\hka\ce \Frac \hva$, and a $V$-torus $T$, we have the product formula
        \[
         \tst  T(\hka)=\mathrm{Im}\bigl(T(\via)\ra T(\hka)\bigr)\cdot T(\hva).
        \]
      \elem
      \bpf
          The left-hand side contains the right-hand side, so it remains to show that every element of $T(\hka)$ is a product of elements of $\mathrm{Im}\bigl(T(\via)\ra T(\hka)\bigr)$ and $T(\hva)$. 
          Consider the commutative diagram
     \[
      \tst \tikz {
      \node  (A) at (0.5,0.6) {$0$};
      \node  (B) at (2,0.6) {$T(V)$};
      \node  (C) at (4,0.6) {$T(\via)$};
      \node  (D) at (6.5,0.6) {$H^1_{V/(a)}(V,T)$};
      \node  (E) at (9.3,0.6) {$H^1(V,T)$};
      \node  (F) at (12,0.6) {$H^1(\via,T)$};
      \node  (a) at (0.5,-0.5)   {$0$};
      \node  (b) at (2,-0.5) {$T(V^{\mathrm{h}})$};
      \node  (c) at (4,-0.5) {$T(V^{\mathrm{h}}[\f{1}{a}])$};
      \node  (d) at (6.5,-0.5) {$H^1_{\small{V^{\mathrm{h}}/(a)}}(V^{\mathrm{h}},T)$};
      \node  (e) at (9.3,-0.5) {$H^1(V^{\mathrm{h}}, T)$};
      \node  (f) at (12,-0.5) {$H^1(V^{\mathrm{h}}[\f{1}{a}],T)$};
      \node  (aa) at (0.5,-1.7) {$0$};
      \node  (bb) at (2,-1.7) {$T(\hva)$};
      \node  (cc) at (4,-1.7) {$T(\hka)$};
      \node  (dd) at (6.5,-1.7) {$H^1_{\small{\nhva/(a)}}(\hva, T)$};
      \node  (ee) at (9.3,-1.7) {$H^1(\hva,T)$};
      \node  (ff) at (12,-1.7) {$H^1(\hka, T)$,};
      \draw [draw = black, thin,
      arrows={
      - Stealth }]
      (A) edge  (B);
      \draw [draw = black, thin,
      arrows={
      - Stealth }]
      (B) edge  (C);
      \draw [draw = black, thin,
      arrows={
      - Stealth }]
      (C) edge (D);
      \draw [draw = black, thin,
      arrows={
      - Stealth }]
      (D) edge (E);
      \draw [draw = black, thin,
      arrows={
      - Stealth }]
      (E) edge (F);
      \draw [draw = black, thin,
      arrows={
      - Stealth }]
      (B) edge (b);
      \draw [draw = black, thin,
      arrows={
      - Stealth }]
      (C) edge (c);
      \draw [draw = black, thin,
      arrows={
      - Stealth }]
      (D) edge (d);
      \draw [draw = black, thin,
      arrows={
      - Stealth }]
      (E) edge (e);
      \draw [draw = black, thin,
      arrows={
      - Stealth }]
      (F) edge (f); 
      \draw [draw = black, thin,
      arrows={
      - Stealth }]
      (a) edge (b);
      \draw [draw = black, thin,
      arrows={
      - Stealth }]
      (b) edge (c);
      \draw [draw = black, thin,
      arrows={
      - Stealth }]
      (c) edge (d);
      \draw [draw = black, thin,
      arrows={
      - Stealth }]
      (d) edge (e);
      \draw [draw = black, thin,
      arrows={
      - Stealth }]
      (e) edge (f);
      \draw [draw = black, thin,
      arrows={
      - Stealth }]
      (aa) edge (bb);
      \draw [draw = black, thin,
      arrows={
      - Stealth }]
      (bb) edge (cc);
      \draw [draw = black, thin,
      arrows={
      - Stealth }]
      (cc) edge (dd);
      \draw [draw = black, thin,
      arrows={
      - Stealth }]
      (dd) edge (ee);
      \draw [draw = black, thin,
      arrows={
      - Stealth }]
      (ee) edge (ff);
      \draw [draw = black, thin,
      arrows={
      - Stealth }]
      (b) edge (bb);
      \draw [draw = black, thin,
      arrows={
      - Stealth }]
      (c) edge (cc);
      \draw [draw = black, thin,
      arrows={
      - Stealth }]
      (d) edge (dd);
      \draw [draw = black, thin,
      arrows={
      - Stealth }]
      (e) edge (ee);
      \draw [draw = black, thin,
      arrows={
      - Stealth }]
      (f) edge (ff);
      }
      \]
       where $V^{\mathrm{h}}$ is the Henselization of $V$ and the rows are exact sequences of local cohomology \cite{SGA4II}*{V, 6.5.3}. 
       By \SP{0F0L}, $V^{\mathrm{h}}$ is also the $a$-Henselization of $V$, hence the $a$-adic completion of $V^{\mathrm{h}}$ is $\hva$ (see \cite{FK18}*{0, 7.3.5}).  
       By the tori case \Cref{GS-tori}, the three horizontal morphisms in the two rightmost squares are injective. 
       The excision \cite{Mil80}*{III, 1.28} combined with a limit argument yield an isomorphism $H^{1}_{V/(a)}(V,T)\cong H^{1}_{V^{\mathrm{h}}/(a)}(V^{\mathrm{h}},T)$. 
       Therefore, a diagram chase gives the following decomposition
       \begin{equation}\label{pd1}
         \tst  T(V^{\mathrm{h}}[\f{1}{a}])=\mathrm{Im}\p{T(\via)\ra T(V^{\mathrm{h}}[\f{1}{a}])}\cdot T(V^{\mathrm{h}}).
       \end{equation}
       By \cite{BC22}*{2.2.17}, the image of $T(V^{\mathrm{h}}[\f{1}{a}])\ra T(\hka)$ is dense. 
       The openness of $T(\hva)\subset T(\hka)$ provided by \Cref{ces-par}~\ref{cl-open}, and \Cref{open-product-closure} imply that 
       \begin{equation}\label{pd2}
         \tst \mathrm{Im}\bigl(T(V^{\mathrm{h}}[\f{1}{a}])\ra T(\hka)\bigr)\cdot T(\hva)=\ov{\mathrm{Im}\bigl(T(V^{\mathrm{h}}[\f{1}{a}])\ra T(\hka)}\bigr)\cdot T(\hva)=T(\hka).
       \end{equation}
       Combining (\ref{pd1}) and (\ref{pd2}), we obtain the product formula for the case of tori.
      \epf
      \bprop\label{decomp-gp}

      For a valuation ring $V$ of rank $n>0$, the prime $\fp\subset V$ of height $n-1$, an element $a\in \fm_V\backslash\fp$, the $a$-adic completion $\hva$ of $V$ with $\hka\ce \Frac \hva$, a reductive $V$-group scheme $G$, the subgroup $G(\hva)\subset G(\hka)$ and the image $\mathrm{Im}(G(\via))$ of the map $G(\via)\ra G(\hka)$, we have
      \[
      G(\hka)=\mathrm{Im}\p{G(V\mathsmaller{[\f{1}{a}]})}\cdot G(\hva).
      \]
      \eprop
      \bpf
      The right-hand side is contained in the left-hand side, so it remains to show that every element of $G(\hka)$ is a product of elements of $\mathrm{Im}\p{G(\via)}$ and $G(\hva)$. 
      The proof is divided into two cases.

      \textbf{Case 1: without proper parabolic subgroups}

      The case when $G_{\nhva}$ is anisotropic follows from \Cref{aniso-int-rat}~\ref{ani-eq}. 
      If $G_{\nhva}$ contains no proper parabolic subgroup and $\rad(G_{\nhva})$ contains a nontrivial split torus of $G_{\nhva}$, we consider the commutative diagram
      \begin{equation}
      \begin{tikzcd}
      0  \arrow{r}
      & \rad(G)(\hva) \arrow{d} \arrow{r} & G(\hva) \arrow{d} \arrow{r} & (G/\rad(G))(\hva) \ar[equal]{d}        \arrow{r} & H^1(\hva, \rad(G)) \arrow{d} \\
      0 \arrow{r} & \rad(G)(\hka) \arrow{r} & G(\hka) \arrow{r} & (G/\rad(G))(\hka) \arrow{r} & H^1(\hka,\rad(G))
      \end{tikzcd}
      \end{equation}
      with exact rows, where the equality follows from \Cref{ani-cri} and \Cref{aniso-int-rat}~\ref{ani-eq}. 
      Since $\rad(G_{\nhva})$ is a torus, by \Cref{GS-tori}, the last vertical arrow is injective. 
      Thus, a diagram chase gives $G(\hka)=\rad(G)(\hka)\cdot G(\hva)$ so the product formula for $\rad(G)$ (\Cref{pd-tori}) leads to the assertion.

      \textbf{Case 2: with a proper parabolic subgroup}

      By \Cref{ani-cri}, the remaining case is when $G_{\nhva}$ contains a proper parabolic subgroup. 
       For a minimal parabolic subgroup $P$ of $G_{\nhva}$, denote its unipotent radical by $U\ce \urad(P)$. 
       As exhibited in \cite{SGA3IIInew}*{XXVI, 6.11},  the centralizer of a maximal split torus $T\subset P$ in $G_{\nhva}$ is a Levi subgroup $L$ of $P$.
       By \emph{ibid.}, 2.4~\emph{ff.}, there is a maximal torus $\wt{T}\subset G_{\nhva}$ containing $T$.
       The proof proceeds as the following steps.

      \emph{\textbf{Step 1}: for the maximal split subtorus $T$ of $P$, we have $T(\hka)\subset \mathrm{Im}(G(\via))\cdot G(\hva)$.} 
       The base change $\wh{T}\ce\wt{T}_{\nhka}$ is a maximal torus of $G_{\nhka}$. 
       For $\wt{T}$ we apply \Cref{lift-tor-int} to $W\ce \ov{\mathrm{Im}(G(\via))}\cap G(\hva)$, so there are a $g\in W$ and a maximal torus $T_0\subset G$ such that $(T_0)_{\nhka}=g\wh{T}g^{-1}$.
          The product formula for tori (\Cref{pd-tori}) shows that $T_0(\hka)=\mathrm{Im}\p{T_{0}(\via)}\cdot T_0(\hva)$.
          Hence we get
          \begin{equation}\label{z-inside}
          \tst \wh{T}(\hka)=g^{-1}T_0(\hka)g=g^{-1}\mathrm{Im}\p{T_0(\via)}\cdot T_0(\hva)g\subset g^{-1}\mathrm{Im}\p{G(\via)}\cdot G(\hva)g.
          \end{equation}
        Since $g\in \ov{\mathrm{Im}\p{G(\via)}}\cap G(\hva)$, (\ref{z-inside}) implies that $\wh{T}(\hka)\subset \ov{\mathrm{Im}\p{G(\via)}}\cdot G(\hva)$.
        Note that \Cref{clopen} gives us $\ov{\mathrm{Im}(G(\via))}\cdot G(\hva)=\mathrm{Im}(G(\via))\cdot G(\hva)$. 
        Consequently, we get
        \begin{equation}\label{tor-decomp}
         \tst T(\hka)\subset \wt{T}(\hka)=\wh{T}(\hka)\subset \mathrm{Im}\p{G(\via)}\cdot G(\hva).
        \end{equation}
      \emph{\textbf{Step 2}: we prove that} $U(\hka)\subset \ov{\mathrm{Im}\p{G(\via)}}$.
      The maximal split torus $T$ acts on $G_{\nhva}$ via the map
\[
   T\times G_{\nhva}\ra G_{\nhva},\qq (t, g)\mapsto tgt^{-1},
\]
inducing a weight decomposition $\mathrm{Lie}(G_{\nhva})=\bigoplus_{\alpha\in X^{\ast}(T)}\Lie(G_{\nhva})^{\alpha}$, where $X^{\ast}(T)$ is the character lattice of $T$. 
The subset $\Phi \subset X^{\ast}(T)-\{0\}$ such that $\mathrm{Lie}(G_{\nhva})^{\alpha}\neq 0$ is the relative root system of $(G_{\nhva}, T)$. 
By \cite{SGA3IIInew}*{XXVI, 6.1;  7.4}, $\mathrm{Lie}(L)$ is the zero-weight space of $\Lie(G_{\nhva})$ and the set $\Phi_{+}$ of positive roots fits into the decomposition  
\[
   \tst \x{$\mathrm{Lie}(P)=\mathrm{Lie}(L)\oplus \bigl(\bigoplus_{\alpha\in \Phi_+}\mathrm{Lie}(G_{\nhva})^{\alpha}\bigr)$\qq  with\qq  $\Lie(U)=\bigoplus_{\alpha\in \Phi_+}\Lie(G_{\nhva})^{\alpha}$.}      
\]
Let $\wt{K}/\hka$ be a Galois field extension splitting $G_{\nhva}$.
By \emph{ibid.}, 2.4~\emph{ff.}, there is a split maximal torus $T\pr\subset L_{\wt{K}}\subset P_{\wt{K}}$ of $G_{\wt{K}}$ containing $T_{\wt{K}}$. 
The centralizer of $T\pr$ in $G_{\wt{K}}$ is itself, which is also a Levi subgroup of a Borel $\wt{K}$-subgroup $B\subset P_{\wt{K}}$.
The adjoint action of $T\pr$ on $G_{\wt{K}}$ induces a decomposition $\Lie(G_{\wt{K}})=\bigoplus_{\alpha\in X^{\ast}(T\pr)}\Lie(G_{\wt{K}})^{\alpha}$, whose coarsening is the base change of $\Lie(G_{\nhva})=\bigoplus_{\alpha\in X^{\ast}(T)}\Lie(G_{\nhva})^{\alpha}$ over $\wt{K}$. 
For the root system $\Phi\pr$ with the positive set $\Phi\pr_{+}$ for the Borel $B$, \emph{ibid.}, 7.12 gives us a surjective map $\eta \colon X^{\ast}(T\pr) \surjects X^{\ast}(T)$ such that $\Phi_{+}\subset \eta(\Phi\pr_+)\subset \Phi_{+}\cup \{0\}$.
By \emph{ibid.}, 1.12, we have a decomposition
\[
\tst  U_{\wt{K}}=\prod_{\alpha\in \Phi\prpr}U_{\wt{K},\alpha},\qq \Lie(U_{\wt{K}})=\bigoplus_{\alpha\in \Phi\prpr}\Lie(G_{\wt{K}})^{\alpha},
   \] 
where $\Phi\prpr\subset \Phi\pr_{+}$ and we have isomorphisms $f_{\alpha}\colon U_{\wt{K},\alpha}\isab \bG_{a,\wt{K}}$.
Since $\Lie(L)\subset \Lie(G_{\nhva})$ is the zero-weight space for the $T$-action, the restriction to $T$ of weights in $\Lie(U_{\wt{K}})$ must be nonzero, that is $\eta(\Phi\prpr)\subset \Phi_{+}$. 
For a cocharacter $\xi\colon \bG_m\ra T$, the dual map $\eta^{\ast}\colon X_{\ast}(T)\injects X_{\ast}(T\pr)$ of $\eta$ sends $\xi$ to a cocharacter $\eta^{\ast}(\xi)\in X_{\ast}(T\pr)$ of $T_{\wt{K}}$.
The adjoint action of $\bG_{m}$ on $U$ induced by $\xi$ is denoted by
\[
  \mathrm{ad}\colon  \bG_m(\hka)\times U(\hka)\ra U(\hka),\qq (t,u)\mapsto \xi(t)u\xi(t)^{-1}.
\]
For the open normal subgroup $N\subset G(\hka)$ constructed in \Cref{open-normal}, the intersection $N\cap U(\hka)$ is open in $U(\hka)$, nonempty and stable under $T(\hka)$-action.
We consider the commutative diagram
\[
\tikz {
\node  (C) at (-5.2,0) {$\bG_m(\hka)\times (N\cap U(\hka))$};
\node  (A) at (0,0) {$T(\hka)\times (N\cap U(\hka))$};
\node  (B) at (4.5,0) {$N\cap U(\hka)$};
\node  (c) at (-5.2,-1.5) {$\bG_m(\hka)\times U(\hka)$};
\node  (a) at (0,-1.5) {$T(\hka)\times U(\hka)$};
\node  (b) at (4.5,-1.5)   {$U(\hka)$};
\node  (d) at (-5.2,-3) {$\bG_m(\wt{K})\times U(\wt{K})$};
\node  (e) at (0,-3) {$T(\wt{K})\times U(\wt{K})$};
\node  (f) at (4.5,-3)   {$U(\wt{K})$.};
\node  (l) at (-2.6, 0.2) {$\xi\times \id$};
\node  (ll) at (2.55,0.2) {$\mathrm{ad}$};
\node  (m) at (-2.6, -1.3) {$\xi\times \id$};
\node  (mm) at (2.55, -1.3) {$\mathrm{ad}$};
\node  (n) at (-2.6, -2.8) {$\xi\times \id$};
\node  (nn) at (2.55, -2.8) {$\mathrm{ad}$};
\draw [draw = black, thin,
arrows={
- Stealth }]
(e) edge  (f);
\draw [draw = black, thin,
arrows={
- Stealth }]
(d) edge  (e);
\draw [draw = black, thin,
arrows={
- Stealth }]
(b) edge  (f);
\draw [draw = black, thin,
arrows={
- Stealth }]
(a) edge  (e);
\draw [draw = black, thin,
arrows={
- Stealth }]
(c) edge  (d);
\draw [draw = black, thin,
arrows={
- Stealth }]
(A) edge  (B);
\draw [draw = black, thin,
arrows={
- Stealth }]
(a) edge (b);
\draw [draw = black, thin,
arrows={
- Stealth }]
(B) edge (b);
\draw [draw = black, thin,
arrows={
- Stealth }]
(A) edge (a);
\draw [draw = black, thin,
arrows={
- Stealth }]
(C) edge  (c);
\draw [draw = black, thin,
arrows={
- Stealth }]
(C) edge  (A);
\draw [draw = black, thin,
arrows={
- Stealth }]
(c) edge  (a);
}
\]
Let $\varpi$ be a topologically nilpotent unit (\ref{pseudo-uniformizer}) of $\hka$.
For an integer $m$, the action of $\varpi^{m}$ on $u\in U(\hka)$ is denoted by $(\varpi^m)\cdot u$.
Let $\wt{u}$ be the image of $u$ in $U(\wt{K})$.
Since $\wt{u}=\prod_{\alpha\in \Phi\prpr}f_{\alpha}({g}_{\alpha})$ with ${g}_{\alpha}\in \wt{K}$,  the image  of $(\varpi^m)\cdot u$ in $U(\wt{K})$ is $\p{\eta^{\ast}(\xi)(\varpi^m)} \wt{u} \p{\eta^{\ast}(\xi)(\varpi^m)}^{-1}$, expressed as the following
\[
\tst  \prod_{\alpha\in \Phi\prpr}\p{\eta^{\ast}(\xi)(\varpi^m)}f_{\alpha}(g_{\alpha})\p{\eta^{\ast}(\xi)(\varpi^m)}^{-1}=\prod_{\alpha\in \Phi\prpr}f_{\alpha}\bigl((\varpi^m)^{\langle \eta^{\ast}(\xi),\alpha\rangle}g_{\alpha}\bigr)=\prod_{\alpha\in \Phi\prpr}f_{\alpha}\bigl((\varpi^m)^{\langle \xi,\eta(\alpha)\rangle}g_{\alpha}\bigr). 
\]
Because $\eta(\Phi\prpr)\subset \Phi_{+}$, we can choose a cocharacter $\xi$ such that $\langle \xi, \eta(\alpha)\rangle$ are strictly positive for all $\alpha\in \Phi\prpr$. 
Then, when $m$ increases, the element $(\varpi^m)\cdot u\in U(\wt{K})$ $a$-adically converges to the identity, and so the same holds in $U(\hka)$.
Thus, since $N\cap U(\hka)$ is an open neighborhood of identity, every orbit of the $T(\hka)$-action on $U(\hka)$ intersects with $N\cap U(\hka)$ nontrivially.
So, we have $U(\hka)=\bigcup_{t\in T(\nhka)}t(N\cap U(\hka))t^{-1}=N\cap U(\hka)$,
which implies that $U(\hka)\subset N$.
By combining with \Cref{open-normal}, we get 
\begin{equation}\label{uni-in}
\tst U(\hka)\subset \ov{\mathrm{Im}\p{G(\via)}}.
\end{equation}

\emph{\textbf{Step~3}: we have $P(\hka)\subset \ov{\mathrm{Im}\p{G(\via)}}\cdot G(\hva)$.}  By \Cref{aniso-int-rat}~\ref{ani-eq}, the quotient $H\ce L/T$ satisfies $H(\hka)=H(\hva)$.
Since $T$ is split, Hilbert's theorem 90 gives the vanishing in the commutative diagram
      \begin{equation}\label{aniso-quot}
      \begin{tikzcd}
      0  \arrow{r}
      & T(\hva) \arrow{d} \arrow{r} & L(\hva) \arrow{d} \arrow{r} & H(\hva) \ar[equal]{d}        \arrow{r} & H^1(\hva, T)=0 \arrow{d} \\
      0 \arrow{r} & T(\hka) \arrow{r} & L(\hka) \arrow{r} & H(\hka) \arrow{r} & H^1(\hka,T)=0
      \end{tikzcd}
      \end{equation}
      with exact rows.
      A diagram chase yields $L(\hka)=T(\hka)\cdot L(\hva)$.
      Combining this with (\ref{tor-decomp}) and (\ref{uni-in}), by \Cref{clopen}, we conclude that
      \begin{equation}\label{par-decomp}
       \tst P(\hka)\subset \ov{\mathrm{Im}\p{G(\via)}}\cdot G(\hva).
      \end{equation}

      \emph{\textbf{Step 4}: the end of the proof.} By \cite{SGA3IIInew}*{XXVI, 4.3.2, 5.2}, there is a parabolic subgroup $Q$ of $G$ such that $P\cap Q=L$ fitting into the surjection
      \begin{equation}\label{rad-surj}
        \urad(P)(\hka)\cdot \urad(Q)(\hka)\surjects G(\hka)/P(\hka).
      \end{equation}
      Applying (\ref{uni-in}) to (\ref{rad-surj}) for $U$ and $\rad^u(Q)$ gives $G(\hka)\subset \ov{\mathrm{Im}(G(\via))}\cdot P(\hka)$, which combined with (\ref{par-decomp}) yields $G(\hka)\subset \ov{\mathrm{Im}(G(\via))}\cdot G(\hva)$.
      With the equality $\ov{\mathrm{Im}(G(\via))}\cdot G(\hva)=\mathrm{Im}(G(\via))\cdot G(\hva)$ verified in \Cref{clopen}, the desired product formula $G(\hka)=\mathrm{Im}(G(\via))\cdot G(\hva)$ follows.\qedhere
      \epf
The following corollary of independent interest shows that torsors under reductive group schemes satisfy arc-patching (see \cite{BM21}), where the arc-cover of $\Spec V$ is of the form $\Spec V/\fp\sqcup \Spec V_{\fp}$.
\bcor\label{arc}
For a valuation ring $V$ of rank $n\geq 1$, the prime $\fp\subset V$ of height $n-1$, and a reductive $V$-group scheme $G$, the following map
\[
   \x{$\mathrm{Im}(G(V_{\fp})\ra G(\kappa(\fp)))\cdot \mathrm{Im}(G(V/\fp)\ra G(\kappa(\fp)))\surjects G(V_{\fp}/\fp)$\qq is surjective.}
\]
\ecor
\bpf
By a limit argument (\SP{01ZC}, \cite{BM21}*{2.22}), we may assume that $V$ contains an element $a$ cutting out the height-one prime ideal of $V/\fp$.
Note that $\via=V_{\fp}$ and the $a$-adic completion of $V/\fp$ is $\hva$.
The affineness of $G$ and \Cref{a-completion}~\ref{fp-a-comp} $V/\fp\isoto V_{\fp}/\fp\times_{\Frac \nhva} \hva$ give us the isomorphism
\[
   G(V/\fp)\isoto G(V_{\fp}/\fp)\times_{G(\Frac \nhva)} G(\hva).
\]
By \Cref{decomp-gp}, the image of $G(V_{\fp})\times G(V/\fp)$ in $G(\Frac \hva)$ generates $G(V_{\fp}/\fp)$. 
\epf
      \bprop\label{rank-one-kernel-trivial}
      For \Cref{GSVal}, proving that \ref{GSV} has trivial kernel for rank one Henselian $V$ suffices.
      \eprop
      \bpf
      A twisting technique \cite{Gir71}*{III, 2.6.1(1)} reduces us to showing that the map \ref{GSV} has trivial kernel.
      The valuation ring $V$ is a filtered direct union of valuation subrings $V_i$ of finite rank (see, for instance, \cite{BM21}*{2.22}).
      Since direct limits commute with localizations, the fraction field $K=\Frac(V)$ is also a filtered direct union of $K_i\ce \Frac(V_i)$.
      A limit argument \cite{Gir71}*{VII, 2.1.6} gives compatible isomorphisms $\tst H^1_{\et}(V,G)\cong \varinjlim_{i\in I}H^1_{\et}(V_i,G)$ and $\tst H^1_{\et}(K,G)\cong \varinjlim_{i\in I}H^1_{\et}(K_i,G)$.
      Thus, it suffices to prove that \ref{GSV} has trivial kernel for $V$ of finite rank, say $n\geq 0$.
      When $n=0$, the valuation ring $V=K$ is a field, so this case is trivial.
      Our induction hypothesis is to assume that \Cref{GSVal} holds for two kinds of valuation rings $V\pr$: (1) for $V\pr$ Henselian of rank one; (2) for $V\pr$ of rank $n-1$.
      Indeed, (1) is only used for the case $n=1$.
      
      Let $\cX$ be a $G$-torsor lying in the kernel of $H^1_{\et}(V,G)\ra H^1_{\et}(K,G)$. For the prime $\mathfrak{p}\subset V$ of height $n-1$, we choose an element $a\in \fm_V\backslash\mathfrak{p}$ and consider the $a$-adic completion $\hva$ of $V$ with fraction field $\hka$.
      The induction hypothesis gives  the triviality of $\cX|_{\via}$ hence a section $s_1\in \cX(\via)$.
      Consequently, $\cX$ is trivial over $\hka$ and by induction hypothesis again, trivial over $\hva$ with $s_2\in \cX(\hva)$.
      By the product formula $G(\hka)=\mathrm{Im}(G(\via))\cdot G(\hva)$ in \Cref{decomp-gp}, there are $g_1\in G(\via)$ and $g_2\in G(\hva)$ such that $g_1 s_1$ and $g_2s_2$ have the same image in $\cX(\hka)$.
      Since $\cX$ is affine over $V$, by \Cref{a-completion}~\ref{fp-a-comp}, we have $\cX(V)\simeq \cX(\via)\times_{\cX(\nhka)}\cX(\hva)$, which is nonempty, so the triviality of $\cX$ follows.
      \epf
      
      \section{Passage to the semisimple anisotropic case}\label{passage-ss-ani}
      After the passage to the Henselian rank-one case \Cref{rank-one-kernel-trivial}, in this section, we further reduce \Cref{GSVal} to the case when $G$ is semisimple anisotropic, see \Cref{no-para}.
      For this, by induction on Levi subgroups, we reduce to the case when $G$ contains no proper parabolic subgroups.
      Subsequently, we consider the semisimple quotient of $G$, which is semisimple anisotropic.
      By using the integrality of rational points of anisotropic groups and a diagram chase, we obtain the desired reduction.
      \bprop\label{no-para}
      To prove \Cref{GSVal}, it suffices to show that \ref{GSV} has trivial kernel in the case when $V$ is a Henselian valuation ring of rank one and $G$ is semisimple anisotropic.
      \eprop
      \bpf
      First, we reduce to the case when $G$ contains no proper parabolics. 
      If $G$ contains a proper minimal parabolic $P$ with Levi $L$ and unipotent radical $\rad^u(P)$, then we consider the commutative diagram
      \[
      \tikz {
      \node  (C) at (-3.2,0) {$H_{\et}^1(V,L)$};
      \node  (A) at (-0.6,0) {$H_{\et}^1(V,P)$};
      \node  (B) at (2,0) {$H_{\et}^1(V,G)$};
      \node  (c) at (-3.2,-1.5) {$H_{\et}^1(K,L)$};
      \node  (a) at (-0.6,-1.5) {$H_{\et}^1(K,P)$};
      \node  (b) at (2.01,-1.5)   {$H_{\et}^1(K,G).$};
      \node  (L) at (-2.9,-0.73) {$l_L$};
      \node  (P) at (-0.3,-0.73) {$l_P$};
      \node  (G) at (2.3,-0.73) {$l_G$};
      \draw [draw = black, thin,
      arrows={
      - Stealth }]
      (A) edge  (B);
      \draw [draw = black, thin,
      arrows={
      - Stealth }]
      (a) edge (b);
      \draw [draw = black, thin,
      arrows={
      - Stealth }]
      (B) edge (b);
      \draw [draw = black, thin,
      arrows={
      - Stealth }]
      (A) edge (a);
      \draw [draw = black, thin,
      arrows={
      - Stealth }]
      (C) edge  (c);
      \draw [draw = black, thin,
      arrows={
      - Stealth }]
      (C) edge  (A);
      \draw [draw = black, thin,
      arrows={
      - Stealth }]
      (c) edge  (a);
      }
      \]
      By \cite{SGA3IIInew}*{XXVI, 2.3}, the left horizontal arrows are bijective.
      If a $G$-torsor $\cX$ lies in $\ker(l_G)$, then it satisfies $\cX(K)\neq \emptyset$.
      By \emph{ibid.}, 3.3; 3.20, the fpqc quotient $\cX/P$ is representable by a scheme which is projective over $V$.
      The valuative criterion of properness gives $(\cX/P)(V)=(\cX/P)(K)\neq \emptyset$, so we can form a fiber product $\cY\ce \cX\times_{\cX/P}\Spec V$ from a $V$-point of $\cX/P      $.
      Since $\cY(K)\neq \emptyset$, the class $[\cY]\in \ker(l_P)$.
      On the other hand, the image of $[\cY]$ in $H_{\et}^1(V, G)$ coincides with $[\cX]$.
      Consequently, the triviality of $\ker(l_L)$ amounts to the triviality of $\ker(l_G)$.
      By \emph{ibid.}, 1.20 and \Cref{rank-one-kernel-trivial}, we are reduced to proving \Cref{GSVal} where $V$ is Henselian of rank one and $G$ has no proper parabolic subgroup, more precisely, to showing that $\ker(H^1(V,G)\ra H^1(K,G))=\{\ast\}$ for such $V$ and $G$.
      
      For the radical $\rad(G)$ of $G$, the quotient $G/\rad(G)$ is $V$-anisotropic, and by \Cref{aniso-int-rat}, satisfies $(G/\rad(G))(V)=(G/\rad(G))(K)$, fitting into the following commutative diagram with exact rows
      \[
      \tikz {
      \node  (A) at (-0.6,0) {$(G/\rad(G))(V)$};
      \node  (B) at (3,0) {$H_{\et}^1(V,\rad(G))$};
      \node  (C) at (6,0) {$H_{\et}^1(V,G)$};
      \node  (D) at (9,0) {$H_{\et}^1(V,G/\rad(G))$};
      \node  (a) at (-0.6,-1.5) {$(G/\rad(G))(K)$};
      \node  (b) at (3,-1.5)   {$H_{\et}^1(K,\rad(G))$};
      \node  (c) at (6,-1.5) {$H_{\et}^1(K,G)$};
      \node  (d) at (9,-1.5) {$H_{\et}^1(K,G/\rad(G)).$};
      \node  (l) at (3.7,-0.75)  {$\scriptstyle{l(\rad(G))}$};
      \node  (ll) at (6.4,-0.75) {$\scriptstyle{l(G)}$};
      \node  (lll) at (9.9,-0.75) {$\scriptstyle{l(G/\rad(G))}$};
      \draw [draw = black, thin,
      arrows={
      - Stealth }]
      (A) edge  (B);
      \draw [draw = black, thin,
      arrows={
      - Stealth }]
      (B) edge  (C);
      \draw [draw = black, thin,
      arrows={
      - Stealth }]
      (C) edge (D);
      \draw [draw = black, thin,
      arrows={
      - Stealth }]
      (a) edge (b);
      \draw [draw = black, thin,
      arrows={
      - Stealth }]
      (b) edge (c);
      \draw [draw = black, thin,
      arrows={
      - Stealth }]
      (c) edge (d);
      \draw [draw = black, thin,
      arrows={
      - Stealth }]
      (B) edge (b);
      \draw [draw = black, thin,
      arrows={
      - Stealth }]
      (C) edge (c);
      \draw (-0.64cm,-3.2em) -- (-0.64cm,-0.8em);
      \draw (-0.56cm,-3.2em) -- (-0.56cm,-0.8em);
      \draw [draw = black, thin,
      arrows={
      - Stealth }]
      (D) edge (d);
      }
      \]
      If $\ker(l(G/\rad(G)))$ is trivial, then by the case of tori \Cref{GS-tori} and Four Lemma, we conclude.
      \epf
      
      \section{Proof of the main theorem}\label{pf}
      In this section, we finish the proof of our main result \Cref{GSVal}.
      By the reduction of \Cref{no-para}, it suffices to deal with semisimple anisotropic group schemes over Henselian valuation rings of rank one.
      In this situation, we argue by using techniques in Bruhat--Tits theory and Galois cohomology to conclude.
      \bthm\label{final-proof}
      For a Henselian rank-one valuation ring $V$ and a semisimple anisotropic $V$-group $G$,
      \[
      \ker(H_{\et}^1(V,G)\ra H_{\et}^1(\Frac V,G))=\{\ast\}.
      \]
      \ethm
      \bpf
      Denote $K\ce \Frac V$ and let $\wt{V}$ be a strict Henselization of $V$ at $\fm_V$ with fraction field $\wt{K}$ as a subfield of a separable closure $K\sep$.
      For the three Galois groups $\GG\ce \Gal(\wt{V}/V)$, $\GG_{\wt{K}}\ce \Gal(K\sep/\wt{K})$ and $\GG_K\ce \Gal(K\sep/K)$, since $\GG\cong \Gal(\wt{K}/K)$, we have $\GG_K/\GG_{\wt{K}}\simeq \GG$.
      An application of the Cartan--Leray spectral sequence yields an isomorphism $H^1_{\et}(V,G)\isom H^1(\GG,G(\wt{V}))$.
      By \cite{SGA4II}*{VIII, 2.1}, we have $H^1_{\et}(K,G)\isom H^1(\GG_K, G(K\sep))$.
      With these bijections, the composite of the following maps $\alpha$ and $\beta$
      \[
      H^1(\GG,G(\wt{V}))\stackrel{\alpha}{\ra} H^1(\GG, G(\wt{K}))\stackrel{\beta}{\ra} H^1(\GG_K,G(K\sep))
      \]
      corresponds to the map $H^1_{\et}(V,G)\ra H^1_{\et}(K,G)$. Hence it suffices to show that $\alpha$ and $\beta$ have trivial kernels.
      For $\beta\colon H^1(\GG,G(\wt{K}))\ra H^1(\GG_K,G(K\sep))$, invoke the inflation-restriction exact sequence \cite{Ser02}*{5.8~a)}
      \[
      0\ra H^1(G_1/G_2, A^{G_2})\ra H^1(G_1,A)\ra H^1(G_2,A)^{G_1/G_2},
      \]
      for which $G_2$ is a closed normal subgroup of a group $G_1$ and $A$ is a $G_1$-group. It suffices to take
      \[
      \x{$G_1\ce \GG_K$, \qq $G_2\ce \GG_{\wt{K}}$,\qq and $A\ce G(K\sep)$.}
      \]
      For $\alpha\colon H^1(\GG, G(\wt{V}))\ra H^1(\GG, G(\wt{K}))$, let $z\in H^1(\GG,G(\wt{V}))$ be a cocycle in $\ker \alpha$, which signifies that
      \begin{equation}\label{cocycle}
        \x{there is an $h\in G(\wt{K})$ such that for every $s\in \GG$,\qq $z(s)=h^{-1} s(h)\in G(\wt{V})$.}
      \end{equation}
      Now we come to Bruhat--Tits theory and consider $G(\wt{V})$ and $hG(\wt{V})h^{-1}$ as two subgroups of $G(\wt{K})$.
      Let $\wt{\sI}(G)$ denote the building of $G_{\wt{K}}$.
      Since $G_{\wt{K}}$ is semisimple, the extended building $\wt{\sI}(G)^{\mathrm{ext}}\ce \wt{\sI}(G)\times (\Hom_{\x{$\wt{K}$-gr.}}(G,\bG_{m,\wt{K}})^{\vee}\otimes_{\b{Z}}\b{R})$ has       trivial vectorial part and equals to $\wt{\sI}(G)$.
      The elements of $G(\wt{K})$ act on the building $\wt{\sI}(G)$. 
      For each facet $F\subset \wt{\sI}(G)$, we consider its stabilizer $P^{\dagger}_F$ and its connected pointwise stabilizer $P^0_F$.
      In fact, there are group schemes $\mathfrak{G}^{\dagger}_{F}$ and $\mathfrak{G}^0_F$ over $\wt{V}$ such that $\mathfrak{G}^{\dagger}_F(\wt{V})=P^{\dagger}_F$ and $\mathfrak{G}^0_F(\wt{V})=P^0_F$, see \cite{BrT2}*{4.6.28}.
      Note that the residue field of $\wt{V}$ is separably closed and the closed fiber of $G_{\wt{V}}$ is reductive, so, by \cite{BrT2}*{4.6.22, 4.6.31}, there is a special point $x$ in the building $\wt{\sI}(G)$ such that the Chevalley group $G_{\wt{V}}$ is the stabilizer $\mathfrak{G}^{\dagger}_x=\mathfrak{G}_x^0$ of $x$ with connected fibers.
      By definition \cite{BrT2}*{5.2.6}, $G(\wt{V})$ is a parahoric subgroup of $G(\wt{K})$.
      Therefore, its conjugate $hG(\wt{V})h^{-1}$ is also a parahoric subgroup $P^0_{h^{-1}\cdot x}$.
      Since $G(\wt{V})$ is $\GG$-invariant, every $s\in \GG$ acts on $hG(\wt{V})h^{-1}$ as follows
      \[
      \tst \x{$s(hG(\wt{V})h^{-1})=s(h)G(\wt{V})s(h^{-1})\stackrel{(\ref{cocycle})}{\longeq}hG(\wt{V})h^{-1}$.}
      \]
      The $\GG$-invariance of $G(\wt{V})$ and $hG(\wt{V})h^{-1}$ amounts to that $x$ and $h\cdot x$ are two fixed points of $\GG$ in $\wt{\sI}(G)$.
      But by \cite{BrT2}*{5.2.7}, the anisotropicity of $G_{K}$ gives the uniqueness of fixed points in $\wt{\sI}(G)$.
      Thus, we have $G(\wt{V})=hG(\wt{V})h^{-1}$, which means that for every $g\in G(\wt{V})$ its conjugate $hgh^{-1}$ fixes $x$.
      This is equivalent to that $g$ fixes $h^{-1}\cdot x$ and to the inclusion of stabilizers $P_x^{\dagger}\subset P^{\dagger}_{h^{-1}\cdot x}$.
      On the other hand, every $\tau\in P^{\dagger}_{h^{-1}\cdot x}$ satisfies $h\tau h^{-1}\cdot x=x$, so $h\tau h^{-1}\in P^{\dagger}_x= G(\wt{V})$.
      Since $h$ normalizes $G(\wt{V})$, this inclusion implies that $\tau\in G(\wt{V})$ and $P^{\dagger}_{h^{-1}\cdot x}\subset G(\wt{V})$.
      Combined with $P^{\dagger}_x\subset P^{\dagger}_{h^{-1}\cdot x}$, this gives $P^{\dagger}_x=P^{\dagger}_{h^{-1}\cdot x}=G(\wt{V})$.
      Therefore, the stabilizer $P^{\dagger}_{h^{-1}\cdot x}$ is also a parahoric subgroup and equals to $P^0_{h^{-1}\cdot x}$.
      By \cite{BrT2}*{4.6.29}, the equality $P^0_x=P^0_{h^{-1}\cdot x}$ implies that $h^{-1}\cdot x=x$, so $h\in P^0_x=G(\wt{V})$, which gives the triviality of $z$.\qedhere
      \epf

\section{Torsors over $V\lps{t}$ and Nisnevich's purity conjecture}\label{purity}
In \cite[1.3]{Nis89}, Nisnevich proposed a conjecture that for a reductive group scheme $G$ over a regular local ring $R$ with a regular parameter $f\in \fm_{R}\backslash \fm_R^2$, every Zariski-locally trivial $G$-torsor over $R[\f{1}{f}]$ is trivial:
\[
  \tst H^1_{\Zar}(R[\f{1}{f}],G)=\{\ast\}.
\]
Recently, Fedorov proved this conjecture when $R$ is semilocal regular defined over an infinite field and $G$ is strongly locally isotropic (that is, each factor in the decomposition of $G^{\mathrm{ad}}$ into Weil restrictions of simple groups is Zariski-locally isotropic); he also showed that the isotropicity is necessary, see \cite{Fed21b}.

In this section, we prove a variant of Nisnevich's purity conjecture when $R$ is a formal power series $V\fps{t}$ over a valuation ring $V$, see \Cref{nis}.
For this, we devise a cohomological property \Cref{lps} of $V\lps{t}$ by taking advantage of techniques of reflexive sheaves. 
\bpp[Coherentness and reflexive sheaves]\label{ref:pre}
A scheme with coherent structure sheaf is \emph{locally coherent}; a quasi-compact quasi-separated locally coherent scheme is \emph{coherent}.
For a valuation ring $V$ with spectrum $S$, by \cite{GR18}*{9.1.27}, every essentially finitely presented affine $S$-scheme is coherent. 
For a locally coherent scheme $X$ and an $\sO_{X}$-module $\sF$, we define the \emph{dual} $\sO_X$-module of $\sF$ $$\sF^{\vee}\ce \sHom_{\sO_X}(\sF,\sO_X).$$
We say that $\sF$ is \emph{reflexive}, if it is coherent and the map $\sF\ra \sF^{\dvee}$ is an isomorphism.
A coherent sheaf $\sG$ has a presentation Zariski-locally $\sO_{X}^{\oplus m}\ra \sO^{\oplus n}_X\ra\sG\ra 0$, whose dual is the exact sequence $0\ra \sG^{\vee}\ra \sO_X^{\oplus n}\ra \sO_X^{\oplus m}$ exhibiting $\sG^{\vee}$ as the kernel of maps between coherent sheaves, hence by \SP{01BY} $\sG^{\vee}$ is coherent, a priori finitely presented.
If $\sF$ is reflexive at a point $x\in X$, then the dual of a presentation $\sO_{X,x}^{\oplus m\pr}\ra \sO_{X,x}^{\oplus n\pr}\ra \sF^{\vee}_x\ra 0$ is a left exact sequence $0\ra \sF_{x}\ra \sO_{X,x}^{\oplus n\pr}\ra \sO_{X,x}^{\oplus m\pr}$.
\epp

\blem[reflexive hull]\label{ref-hull}
Let $X$ be an integral locally coherent scheme and let $\sF$ be a coherent $\sO_X$-module, then $\sF^{\vee}$ and $\sF^{\vee\!\vee}$ are reflexive $\sO_X$-modules.
\elem
\bpf 
It suffices to show that $\sF^{\vee}$ is reflexive.
As $\sF$ is coherent, choose a finite presentation $\sO_X^{\oplus m}\ra \sO_X^{\oplus n}\ra \sF\ra 0$, take its dual and its triple dual, we have the commutative diagram with exact rows
\[
      \tikz {
      \node  (D) at (-4,0) {$0$};
      \node  (C) at (-2.4,0) {$\sF^{\vee}$};
      \node  (A) at (-0.3,0) {$\sO_X^{\oplus n}$};
      \node  (B) at (1.7,0) {$\sO_X^{\oplus m}$};
      \node  (c) at (-2.4,-1.2) {$\sF^{\vee\!\vee\!\vee}$};
      \node  (a) at (-0.3,-1.2) {$\sO_X^{\oplus n}$};
      \node  (b) at (1.7,-1.2)   {$\sO_X^{\oplus m}$.};
      \node  (L) at (-1.4,0.2) {$u$};
      \node  (P) at (-1.3,-1) {$u^{\vee\!\vee}$};
      \draw [draw = black, thin,
      arrows={
      - Stealth }]
      (A) edge  (B);
      \draw [draw = black, thin,
      arrows={
      - Stealth }]
      (D) edge  (C);
      \draw [draw = black, thin,
      arrows={
      - Stealth }]
      (a) edge (b);
      \draw [draw = black, thin,
      arrows={
      - Stealth }]
      (B) edge (b);
      \draw [draw = black, thin,
      arrows={
      - Stealth }]
      (A) edge (a);
      \draw [draw = black, thin,
      arrows={
      - Stealth }]
      (C) edge  (c);
      \draw [draw = black, thin,
      arrows={
      - Stealth }]
      (C) edge  (A);
      \draw [draw = black, thin,
      arrows={
      - Stealth }]
      (c) edge  (a);
      }
      \]
      Our goal is to show that the left most vertical arrow is an isomorphism. 
      Since the other vertical arrows are isomorphisms, a diagram chase reduces us to showing that $u^{\vee\!\vee}$ is injective. 
      Consider the dual of $u$
      \[
      u^{\vee}\colon \sO_X^{\oplus n}\ra \sF^{\vee\!\vee}
      \]
      and its tensor product with the function field $K$ of $X$, we get the following exact sequence 
      \[
      K^{\oplus n}\ra \sF^{\vee\!\vee}\otimes_{\sO_X}K\ra \coker(u^\vee)_K\ra 0.
      \]
      As $\sF$ is finitely presented, by \SP{0583}, we have $\sF^{\vee\!\vee}\otimes_{\sO_X}K\simeq \Hom_K(\sF^{\vee}\otimes_{\sO_X}K,K)$ and we view $K^{\oplus n}$ as $\Hom_K(K^{\oplus n},K)$.
      Note that $u\otimes_{\sO_X}K\colon \sF^{\vee}\otimes_{\sO_X}K\hra K^{\oplus n}$ is injective (since $u$ is injective), we find that $\coker(u^{\vee})_K=0$, that is, $\coker(u^{\vee})$ is a torsion $\sO_X$-module.
      This implies that $\sHom_{\sO_X}(\coker(u^{\vee}),\sO_X)=0$, so we take dual of the exact sequence $\sO_X^{\oplus n}\overset{u^{\vee}}{\ra} \sF^{\vee\!\vee}\ra \coker(u^{\vee})\ra 0$ to get the injectivity of $u^{\vee\!\vee}$.
\epf

\blem[\cite{GR18}*{11.4.1}]\label{ref:gabber}
For a valuation ring $V$ with spectrum $S$, a flat finitely presented morphism of schemes $f\colon X
\ra S$, a coherent $\sO_X$-sheaf $\sF$, 
a point $x\in X$ such that the fiber of $f$ containing $x$ is regular, and the integer $n\ce \dim \sO_{f^{-1}(f(x)),x}$, 
\benumr
\item if $\sF$ is $f$-flat at $x$, then $\mathrm{proj.dim}_{\sO_{X,x}}\sF_x\leq n$;
\item we have $\mathrm{proj.dim}_{\sO_{X,x}}\sF_x\leq n+1$; and
\item if $\sF$ is reflexive at $x$, then $\mathrm{proj.dim}_{\sO_{X,x}}\sF_x\leq \max(0,n-1)$.
\eenum
\elem
\bpf\hfill
\benumr
\item Since $\sO_{X}$ is coherent and $\sF_{x}$ is finitely presented, there is free resolution of $\sF_{x}$ by finite modules
\[
   \cdots \ra P_2\ra P_1\ra P_0\ra \sF_x\ra 0.
\]
It suffices to show that $L\ce \mathrm{Im}(P_n\ra P_{n-1})$ is free.
Now we have the following exact sequence
\[
   0\ra L\ra P_{n-1}\ra \cdots \ra P_1\ra P_0\ra \sF_x\ra 0.
\]
Denote $y= f(x)$. Since $\sF_{x}$ and $\ker(P_{i}\ra P_{i-1})$ are $f$-flat for $1\leq i\leq n-1$, the following sequence
\[
   0\ra L\otimes_{\sO_{X,x}}\sO_{f^{-1}(y),x}\ra  \cdots \ra  P_0\otimes_{\sO_{X,x}}\sO_{f^{-1}(y),x}\ra \sF_x\otimes_{\sO_{X,x}}\sO_{f^{-1}(y),x}\ra 0
\]
is exact. 
Denote $y\ce f(x)$.
For the maximal ideal $\fm_x$ of $\sO_{f^{-1}(y),x}$ at $x$ and the residue field $k(x)$ of $x$ in $\sO_{X,x}$, we note that $L\otimes_{\sO_{X,x}}\p{\sO_{f^{-1}(y),x}/\fm_x\sO_{f^{-1}(y),x}}=L\otimes_{\sO_{X,x}}k(x)$. 
For a free basis $(e_l)_{l\in I}$ generating $L\otimes_{\sO_{X,x}}k(x)$, by Nakayama's lemma, there is a surjective map $u\colon \bigoplus_{l\in I}\sO_{X,x}e_{l}\surjects L$.  
Since $f^{-1}(y)$ is regular, by \SP{00O9}, the module $L\otimes_{\sO_{X,x}}\sO_{f^{-1}(y),x}$ is free. 
Therefore, the map $(u\otimes 1)_{x}: ((\bigoplus_{l\in I}\sO_{X,x}e_l) \otimes_{\sO_{S}}k(y))_x\ra  (L\otimes_{\sO_{S}}k(y))_x$ is an isomorphism. 
By \cite[11.3.7]{EGAIV4}, $u$ is injective. 
Consequently, the $\sO_{X,x}$-module $L$ is free and $\mathrm{proj.dim}_{\sO_{X,x}}\sF_{x}\leq n$.
\item We prove the assertion Zariski-locally. 
There is a surjective map $\sO_{X}^{\oplus m}\ra \sF$, whose kernel $\sG$ is a torsion-free coherent $\sO_X$-module.
Since $V$ is a valuation ring, $\sG$ is $f$-flat, so by (i) we have $\mathrm{proj.dim}_{\sO_X}\sG\leq n$.
Therefore, \SP{00O5} implies that $\mathrm{proj.dim}_{\sO_X}\sF=\mathrm{proj.dim}_{\sO_X}\sG+1\leq n+1$.
\item By the analysis in \S\Cref{ref:pre}, there is an exact sequence $0\ra \sF_{x}\ra \sO_{X,x}^{\oplus k} \stackrel{\phi}{\ra} \sO_{X,x}^{\oplus l}$.
By (ii), we have
\[
   \mathrm{proj.dim}_{\sO_{X,x}}\sF_x\stackrel{\x{\SP{00O5}}}{\longeq}\max (0,\mathrm{proj.dim}_{\sO_{X,x}}(\coker \phi)-2)\leq \max (0,n-1).  \qedhere
\]
\eenum
\epf

Since $(V\fps{t},t)$ is a Henselian pair, by \cite{Ces21}*{3.1.3(b)}, reductive group schemes over $V$ and $V\fps{t}$ are in a one-to-one correspondence under extension-restriction operations. 
Hence, in the sequel, it suffices to assume that reductive group schemes are defined over $V$.
We bootstrap from the case when $G=\GL_n$.

\blem\label{extend}
For a valuation ring $V$, every vector bundle over $V\lps{t}$ extends to a vector bundle over $V\fps{t}$. 
In particular, all $\GL_{n}$-torsors (or equivalently, all vector bundles) over $V\lps{t}$ are trivial:
\[
   H^1_{\et}(V\lps{t},\GL_n)=\{\ast\}.
\]
\elem
\bpf
The Henselization $V\{t\}$ of $V[t]$ along $tV[t]$ is a filtered direct limit of \'etale ring extensions $R_{i}$ over $V[t]$ with isomorphisms $V[t]/tV[t]\isoto R_{i}/tR_i$.
By \cite[2.1.22]{BC22}, a vector bundle $\sE$ over $V\lps{t}$ descends to a vector bundle $\sE\pr$ over $V\{t\}[\f{1}{t}]$. 
By a limit argument \cite[VII, 2.1.6]{Gir71}, we have $H^1_{\et}(V\{t\}[\f{1}{t}],\GL_n)=\varinjlim_iH^1_{\et}(R_i[\f{1}{t}],\GL_n)$ so $\sE\pr$ descends to a vector bundle $\sE_{i_0}$ over $R_{i_0}[\f{1}{t}]$ for an $i_{0}$.
Due to \cite{GR18}*{10.3.24~(ii)}, $\sE_{i_0}$ extends to a finitely presented quasi-coherent sheaf ${\sW_{i_0}}$ on $R_{i_0}$.
Note that $R_{i_0}$ is coherent (\ref{ref:pre}), by \SP{01BZ}, $\sW_{i_0}$ is coherent.
By \Cref{ref-hull}, $\sH_{i_0}\ce \sW_{i_0}^{\dvee}$ is reflexive. 
For the morphism $f\colon \Spec R_{i_0}\ra \Spec V$, we exploit \Cref{ref:gabber}~(iii) to conclude that ${\sH}_{i_0}$ is free. 
Consequently, $\sE_{i_0}$ extends to the vector bundle $(\sH_{i_0})_{V\fps{t}}$ over $V\fps{t}$. Since $\sE_{i_0}=\sH_{i_0}|_{V\lps{t}}$ is trivial, $\sE$ is trivial.
\epf

 The anisotropic (indeed, the `wound') case of the following \Cref{lps}~\ref{lps-iii} was established in \cite{FG21}*{Cor.~4.2}, where the authors considered formal power series over general rings.

\bprop\label{lps}
For a valuation ring $V$ with fraction field $K$ and a $V$-reductive group scheme $G$,
\benum
\item\label{bij} the following natural map of pointed sets induced by base change is bijective:
\[
   H^1_{\et}(V\fps{t},G)\isoto H^1_{\et}(V\lps{t},G)\times_{H^1_{\et}(K\lps{t},G)}H^1_{\et}(K\fps{t},G);
\]
\item 
the map $H^1_{\et}(V\lps{t},G)\ra H^1_{\et}(K\lps{t},G)$ has trivial kernel; and 
\item\label{lps-iii} the map $H^1_{\et}(V\fps{t},G)\ra H^1_{\et}(V\lps{t},G)$ has trivial kernel.
\eenum
\eprop
\bpf\hfill
\benum
\item 
First, we show the surjectivity. 
If there are torsor classes $\alpha\in H^1_{\et}(K\fps{t},G)$ and $\beta\in H^1_{\et}(V\lps{t},G)$ whose images in $H^1_{\et}(K\lps{t},G)$ coincide, then we find a torsor class $\gamma\in H^1_{\et}(V\fps{t},G)$ whose restrictions are $\alpha$ and $\beta$.
Recall the nonabelian cohomology exact sequence \cite{Gir71}*{III, 3.2.2}
\[
    (\GL_{n,V\fps{t}}/G)(R)\ra H^1_{\et}(R,G)\ra H^1_{\et}(R,\GL_{n})
\]
such that the set of $\GL_{n}(R)$-orbits $\GL_{n}(R)\backslash (\GL_{n,V\fps{t}}/G)(R)$ embeds into $H^1_{\et}(R,G)$, where $R$ can be $V\lps{t}$, $K\lps{t}$, or $K\fps{t}$.
Recall that by \Cref{extend}, we have $H^1_{\et}(V\lps{t},\GL_n)=\{\ast\}$ and note that $H^1_{\et}(K\fps{t},\GL_n)=\{\ast\}$, so there are $\wt{\gA}\in (\GL_{n,V\fps{t}}/G)(K\fps{t})$ and $\wt{\gB}\in (\GL_{n,V\fps{t}}/G)(V\fps{t})$ whose images are $\gA$ and $\gB$ respectively and such that the images of $\wt{\gA}$ and $\wt{\gB}$ in $(\GL_{n,V\fps{t}}/G)(K\lps{t})$ are in the same $\GL_n(K\lps{t})$-orbit.
By the valuative criterion for properness of the affine Gra{\ss}mannian,
\[
    \GL_n(K\lps{t})=\GL_n(K\fps{t})\cdot \GL_n(V\lps{t})
\]
holds, so up to group translations, we may assume that the images of $\wt{\gA}$ and $\wt{\gB}$ in $(\GL_{n,V\fps{t}}/G)(K\lps{t})$ are identical. 
Because $G$ is reductive, by \cite{Alp14}*{9.7.7}, the quotient $\GL_{n,V\fps{t}}/G$ is affine over $V\fps{t}$.
Thus, the fiber product $V\fps{t}\isoto V\lps{t}\times_{K\lps{t}}K\fps{t}$ induces the following bijection of sets
\[
    \qq (\GL_{n,V\fps{t}}/G)(V\fps{t})\isoto (\GL_{n,V\fps{t}}/G)(K\fps{t})\times_{(\GL_{n,V\fps{t}}/G)(K\lps{t})}(\GL_{n,V\fps{t}}/G)(V\lps{t}).
\]
Consequently, there is a $\wt{\gG}\in (\GL_{n,V\fps{t}}/G)(V\fps{t})$ corresponding to $(\wt{\gA}, \wt{\gB})$.
In particular, the image $\gG\in H^1_{\et}(V\fps{t},G)$ of $\wt{\gG}$ is a desired torsor class that induces $\gA$ and $\gB$, hence the surjectivity of \ref{bij}.

It remains to show the injectivity. 
By \cite{GR18}*{5.8.14}, we have bijections $H^1_{\et}(V\fps{t},G)\simeq H^1_{\et}(V,G)$ and $H^1_{\et}(K\fps{t},G)\simeq H^1_{\et}(K,G)$.
Then the Grothendieck--Serre for valuation rings \Cref{GSVal} implies that $H^1_{\et}(V\fps{t},G)\ra H^1_{\et}(K\fps{t},G)$ has trivial kernel.
Therefore, the map of \ref{bij} is indeed injective hence bijective.

\item For a $G_{V\lps{t}}$-torsor $X$ trivializes over $K\lps{t}$, we take a trivial $G_{K\fps{t}}$-torsor $X\pr$ over $K\fps{t}$ with an isomorphism $\iota\colon X|_{K\lps{t}}\isoto X\pr|_{K\lps{t}}$. 
By (a), there is a $G_{V\fps{t}}$-torsor $\cX$ restricts to $X$ such that $\cX_{K\fps{t}}$ is trivial.
By the main result \Cref{GSVal} and \cite[5.8.14]{GR18}, the map $H^1_{\et}(V\fps{t},G)\hra H^1_{\et}(K\fps{t},G)$ is injective. 
Hence, the torsor $\cX$ that restricts to $X$ is trivial. 
\item By the Grothendieck--Serre over valuation rings (\Cref{GSVal}) and \cite[5.8.14]{GR18}, the map
\[
   H^1_{\et}(V\fps{t},G)\ra H^1_{\et}(K\fps{t},G)
\]
is injective.
Since $K\fps{t}$ is a discrete valuation ring, the map $H^1_{\et}(K\fps{t},G)\ra H^1_{\et}(K\lps{t},G)$ is injective.
The injective map $H^1_{\et}(V\fps{t},G)\ra H^1_{\et}(K\lps{t},G)$ factors through $H^1_{\et}(V\fps{t},G)\ra H^1_{\et}(V\lps{t},G)$, hence the later is injective. \qedhere
\eenum
\epf    
      
Now we prove a variant of the Nisnevich's purity conjecture for formal power series over valuation rings.   
\bcor\label{nis}
For a reductive group scheme $G$ over a valuation ring $V$, every Zariski-locally trivial $G$-torsor over $V\lps{t}$ is trivial, that is, we have
\[
   H^1_{\Zar}(V\lps{t},G)=\{\ast\}.
\]
\ecor
\bpf
A Zariski $G$-torsor over $V\lps{t}$ is an \'etale $G$-torsor over $V\lps{t}$ trivializing over $K\lps{t}$. 
Hence the assertion follows from \Cref{lps}~(b).
\epf
\newpage
\begin{appendix}

\section{Valuation rings and valued fields}\label{apdx}
The purpose of this appendix is to list the common properties of valuation rings and valued fields, especially those used in this article, and to be as concise and brief as possible.
We therefore try to cite the literature just for endorsement, even though some of the arguments can be carried out directly.
\bpp[Valuation rings]\label{valuation rings}
For a field $K$, a subring $V\subset K$ such that every $x\in K$ satisfies that $x\in V$ or $x^{-1}\in V$ or both is a \emph{valuation ring} of $K$ (\SPD{052K}{00IB}).
For the groups of units $K^{\times}$ and $V^{\times}$, the quotient $\GG\ce K^{\times}/V^{\times}$ is an abelian group with respect to the multiplications in $K^{\times}$.
The quotient map $\nu\colon K^{\times}\ra \GG$ induces a map $V\backslash\{0\}\subset K^{\times}\ra \GG$, also denoted by $\nu$.
This map $\nu$ is the \emph{valuation} associated to $V$.
Even though the \emph{rank} of $\GG$ (and, of $V$) is the ``order type'' of the collection of convex subgroups (\cite{EP05}*{p.~26, 29}), in practice, one may identify the rank of $V$ as its Krull dimension when it is finite (\cite{EP05}*{Lem.~2.3.1}).
The abelian group $\GG$ has an order $\geq$: for $\gamma,\gamma\pr\in \GG$, we declare that $\gamma\geq \gamma\pr$ if and only if  $\gamma-\gamma\pr$ is in the image of $\nu \colon V\backslash\{0\}\ra \GG$. 
By \SP{00ID}, $(\GG,\geq)$ is a totally ordered abelian group, called the \emph{value group} of $V$.
If $\GG\simeq \b{Z}$, then $\nu$ is a \emph{discrete valuation}. 
Conversely, given a totally ordered abelian group $(\GG,\geq,+)$, if there is a surjection $\nu\colon K^{\times}\surjects \GG$ such that for all $x,y\in K$, we have $\nu(xy)=\nu(x)+\nu(y)$ and $\nu(x+y)\geq \min\{\nu(x),\nu(y)\}$, then $\nu$ extends to a map $K\surjects \GG\cup \{\infty\}$ by declaring that $\nu(x)=\infty$ if and only if $x=0$, where $\infty$ is a symbol whose sum with any element is still $\infty$; such $\nu$ is also a \emph{valuation} on $K$ (\cite{EP05}*{p.~28}). 
If a field $K$ is equipped with a valuation $\nu$, then the pair $(K,\nu)$ is called a \emph{valued field}. 
Every valuation $\nu$ on $K$ gives rise to a valuation ring $V(\nu)\subset K$ as the following
\[
   V(\nu)\ce \{x\in K\,|\,\nu(x)\geq 0\},
\]
and every valuation ring of $K$ comes from a valuation (\cite{EP05}*{Prop.~2.1.2}).
There may exist different valuations $\nu$ and $\nu\pr$ on a field $K$, yielding different valuation rings of $K$.
Two valuations $\nu$ and $\nu\pr$ on $K$ are \emph{equivalent}, if they define the same valuation rings $V(\nu)=V(\nu\pr)$.
By \cite{EP05}*{Prop.~2.1.3}, $\nu$ and $\nu\pr$ are equivalent if and only if there is an isomorphism of ordered groups $\iota\colon \GG_{\nu}\isoto \GG_{\nu\pr}$ such that $\iota \circ\nu = \nu\pr$.
\epp


\bprop\label{basic-one} 
Let $V$ be a valuation ring of a field $K$ with value group $\GG$ and $\fp\subset V$ a prime ideal.
\benumr
\item\label{normal-domain} $V$ is a normal local domain and every finitely generated ideal of $V$ is principal;
\item\label{prime-ext} for the localization $V_{\fp}$ of $V$ at $\fp$, we have $\fp=\fp V_{\fp}$;
\item\label{loc-quot} $V_{\fp}$ is a valuation ring for $K$ and $V/\fp$ is a valuation ring for the residue field $\kappa(\fp)=V_{\fp}/\fp$;
\item\label{fp-val} we have an isomorphism $V\isoto V/{\fp}\times_{V_{\fp}/\fp }V_{\fp}$ and thus $\Spec V=\Spec V/\fp \sqcup_{\Spec (V_{\fp}/\fp)} \Spec V_{\fp}$;
\item\label{exact-seq-val-grps} for the value groups $\GG_{V_{\fp}}$ and $\GG_{V/\fp}$ of $V_{\fp}$ and of $V/\fp$ respectively, we have isomorphisms
\[
    \text{$\GG_{V/\fp}\simeq (V_{\fp})^{\times}/V^{\times}$ \q and\qq   $\GG_{V}/\GG_{V/\fp}\simeq \GG_{V_{\fp}}$,}
\]
corresponding to the short exact sequence $1\ra (V_{\fp})^{\times}/V^{\times}\ra K^{\times}/V^{\times}\ra K^{\times}/(V_{\fp})^{\times}\ra 1$;

\item\label{hens-val} the Henselization and the strict Henselization of $V$ are valuation rings with value groups $\GG$;
\item\label{hens-preserve} if $V$ is Henselian, then $V_{\fp}$ and $V/\fp$ are Henselian valuation rings.
\eenum
\eprop
\bpf
For \ref{normal-domain}, see \cite{FK18}*{Ch.~0, 6.2.2}. 
To show \ref{prime-ext}, we write every element in $\fp V_{\fp}$ as $a/b$, where $a\in \fp V$ and $b\in V\backslash \fp$. 
If $a/b\in V$ then $a/b\in \fp$. 
Since $V$ is a valuation ring, it remains the case when $b/a\in V$.
       Then $b\in \fp V$, which leads to a contradiction.
       For \ref{loc-quot}, see \cite{FK18}*{Ch.~0, Prop.~6.4.1}.
       For \ref{fp-val}, we note that $V=\{x\in V_{\fp}| (x\mod \fp V_{\fp}) \in V/\fp\}$ (\cite{FK18}*{Ch.~0, Prop.~6.4.1}). 
       The spectral aspect follows from \SP{0B7J}.
       For \ref{exact-seq-val-grps}, we first deduce from the fiber product $V\simeq V/\fp\times_{V_{\fp}/\fp} V_{\fp}$ that $\GG_{V/\fp}=\kappa(\fp)^{\times}/(V/\fp)^{\times} \simeq (V_{\fp})^{\times}/V^{\times}$ then substitute this into the short exact sequence $1\ra \Frac(V/\fp)^{\times}/(V/\fp)^{\times}\ra K^{\times}/V^{\times}\ra K^{\times}/(V_{\fp})^{\times}\ra 1$.
       For \ref{hens-val}, see \SP{0ASK}.
        For \ref{hens-preserve}, note that $V_{\fp}$ and $V/\fp$ are valuation rings due to \ref{loc-quot}.
        By \SP{05WQ}, $V/\fp$ is Henselian.
         For $V_{\fp}$, we use Gabber's criterion \SP{09XI} to check that every monic polynomial 
      \[
      \x{\qq $f(T)=T^N(T-1)+a_NT^N+\cdots+a_1T+a_0$,\qq where $a_{i}\in \fp V_{\fp}$ for $i=0,\cdots, N$ and $N\geq 1$}
      \]       
      has a root in $1+\fp V_{\fp}$. 
      Note that this criterion only involves $\fp V_{\fp}$.
      Here, by \ref{prime-ext}, $\fp V_{\fp}$ is equal to $\fp$.
      By \SP{0DYD}, the Henselianity of $V$ implies that $(V,\fp)$ is a Henselian pair, thereby we conclude.
\epf
\bpp[Valuation topologies]\label{valuation topologies}
Given a field $K$ with a valuation $\nu\colon K\surjects \GG\cup\{\infty\}$, for each $\gamma\in \GG$ and each $x\in K$, we define the \emph{open ball} $U_{\gamma}(x)\subset K$ with center $x$ and \emph{radius} $\gamma$, as the following subset 
\[
   U_{\gamma}(x)\ce \{y\in K\,|\,\nu(y-x)>\gamma\}.
\]
All open balls $(U_{\gamma}(x))_{\gamma\in \GG}$ form an open neighborhood base of $x$ and generates a topology on $K$, the \emph{valuation topology} determined by $\nu$. 
Under this topology, the valued field $(K,\nu)$ has a unique (up to isomorphisms) field  extension $(\wh{K},\wh{\nu})$ that is complete in which $K$ is dense (\cite{EP05}*{Thm.~2.4.3}), that is, the \emph{completion} of $(K,\nu)$ with respect to the valuation topology. 
Similarly, the valuation ring $\wh{V}$ of $(\wh{K},\wh{\nu})$ is the valuative completion of $V$.
The inequality $\nu(x+y)\geq \min\{\nu(x),\nu(y)\}$ leads to various topological properties, some of which are counter-intuitive.
In the sequel, we let $B_{\gamma}(x)\ce \{z\in K\,|\, \nu(z-x)\geq \gamma\}$ and $S_{\gamma}(x)\ce \{z\in K\,|\, \nu(z-x)=\gamma\}$ be the \emph{closed ball} and the \emph{sphere} with center $x$ and \emph{radius} $\gamma$ respectively.
\epp
\bprop\label{top-prop}
For a valued field $(K,\nu)$ with the valuation topology and elements $x\in K$ and $\gamma\in \GG$,
\benumr
\item\label{ineq} for $y,z\in K$, the smallest and second smallest among $\nu(x-y), \nu(y-z)$, and $\nu(z-x)$ are equal;
\item every point of the closed ball $B_{\gamma}(x)$ is a center: for all $y\in B_{\gamma}(x)$, we have $B_{\gamma}(y)=B_{\gamma}(x)$;
\item\label{clopenn} every closed  ball is open and every open ball is closed;
\item any pair of balls in $K$ are either disjoint or nested;
\item the sphere $S_{\gamma}(x)$ is both closed and open, hence it is not the boundary $\partial B_{\gamma}(x)$ of $B_{\gamma}(x)$.
\eenum
In particular, the valuation topology on $(K,\nu)$ is Hausdorff and the valuation ring $V(\nu)\subset K$ is clopen.
\eprop
\bpf 
If \ref{ineq} holds, then for any $a\neq b$ in $K$ and $\delta\ce \nu(a-b)$, we have $U_{2\delta}(a)\cap U_{2\delta}(b)\neq \emptyset$, hence $K$ is Hausdorff.
The assertion \ref{ineq} follows from the inequality $\nu(c+d)\geq \min\{\nu(c),\nu(d)\}$ for all $c,d\in K$, and the other assertions follow from \ref{ineq}, see the arguments in \cite{EP05}*{p.~45 and Rem.~2.3.3} and \cite{PGS10}*{p.~3}. 
\epf 

\bpp[Absolute values]\label{abv}
Let $K$ be a field. 
An \emph{absolute value} on $K$ is a function $|\cdot|\colon K\ra \b{R}_{\geq 0}$ such that (1) $\abs{x}=0$ if and only if $x=0$; (2) $\abs{xy}=\abs{x}\cdot\abs{y}$; and (3) $\abs{x+y} \leq \abs{x}+ \abs{y}$ (triangle inequality).
We say that $|\cdot|$ is \emph{archimedean}, if $\abs{\b{N}}\subset \b{R}_{\geq 0}$ is unbounded; $|\cdot|$ is \emph{nonarchimedean}, if $\abs{\b{N}}\subset \b{R}_{\geq 0}$ is bounded.
These notions originate from the `Archimedean property': for arbitrary positive real numbers $x$ and $y$, there is $n\in \b{N}$ such that $xn>y$. 
In fact, an absolute value $|\cdot|$ is nonarchimedean if and only if it satisfies the \emph{strong triangle inequality} $\abs{x+y}\leq \max\{\abs{x},\abs{y}\}$: one takes $M$ such that $\abs{\b{N}}<M$ and notes that
\[
   \tst \abs{x+y}^n\leq \sum_{k=0}^{n}\abs{n\choose{k}}\cdot \abs{x}^k\abs{y}^{n-k}\leq (n+1)M\cdot \max\{\abs{x},\abs{y}\}^n,
\]
whose $n$-th root when $n\ra +\infty$ yields $\abs{x+y}\leq \max\{\abs{x},\abs{y}\}$.
In particular, by checking the axioms of valuations (\ref{valuation rings}), an absolute value $|\cdot|\colon K\ra \b{R}_{\geq 0}$ is nonarchimedean if and only if there is a valuation $\nu\colon K\ra \Gamma\cup\{\infty\}$ of rank one (a value group is of rank one if and only if it is embeddable into $\mathbf{R}$ as a totally ordered abelian subgroup, so $\Gamma\subset \mathbf{R}$) such that $e^{-\nu(\cdot)}=|\cdot|$.
\epp

\bpp[Huber rings and Tate rings]\label{pseudo-uniformizer} 
Let $R$ be a topological ring. We say that
\begin{itemize}
\item[-]$R$ is \emph{adic}, if it has an ideal $I\subset R$ such that $\{I^n\}_{n=1}^{+\infty}$ form a basis of open neighborhoods of $0\in R$;
\item[-] $R$ is \emph{Huber}, if it has an open subring $R_0$ with a finitely generated ideal $I\subset R_0$ making $R_0$ adic;
\item[-] $R$ is \emph{Tate}, if it is Huber and has a \emph{topologically nilpotent unit} $\varpi\in R\backslash\{0\}$, that is, $\lim_{n\ra +\infty}\varpi^n=0$.
\end{itemize}
\epp

Now, we present a relation (\emph{cf.}~\cite{Hub96}*{I, Def.~1.1.4}) between valuation topologies and the notions above.

\bprop\label{val-top}
Let $(K,\nu)$ be a valued field with valuation ring $V$.
The following are equivalent:
\benumr 
\item\label{vt-i} $V$ has a prime ideal of height one; 
\item\label{vt-ii} the valuation topology on $K$ is induced by a valuation of rank one;
\item\label{vt-iii} $K$ is a Tate ring under its valuation topology;
\item\label{vt-iv} $K$ has a topologically nilpotent unit for the valuation topology.
\eenum
In particular, there exist nonzero topologically nilpotent elements $\varpi\in V$, and every such $\varpi$ satisfies that $\sqrt{(\varpi)}$ is the prime ideal of height one in $V$.
\bpf 
Before proving the equivalences, first note that the set of all ideals of $V$ ordered by inclusion is totally ordered.
For two ideals $I,J\subset V$, if there is an element $j\in J$ such that $j\not\in I$, then $ji^{-1}\not\in V$ for all $i\in I\backslash \{0\}$.
By the definition of valuation rings, $ij^{-1}\in V$ for all $i\in I$.
This implies that $I\subset (j)\subset J$.

\ref{vt-i}$\Rightarrow$\ref{vt-iv}. For the prime $\fp\subset V$ of height one, we claim that any $\varpi\in \fp\backslash \{0\}$ is topologically nilpotent. 
For any $\gamma\in \GG$, it suffices to find an $n\in \b{Z}_+$ such that $\varpi^n\in U_{\gamma}=\{x\in K\,|\,\nu(x)>\gamma \}$.
Since $\nu\colon K\surjects \GG$ is surjective, we show that for any $\f{a}{b}\in K$ where $a,b\in V\backslash \{0\}$, there is $n\in \b{Z}_+$ such that $\nu(\varpi^n)>\nu(a)-\nu(b)$, in particular, such that $\nu(\varpi^n)>\nu(a)$ suffices. 
If $\nu(a)\geq \nu(\varpi^n)$ holds for all $n$, then $a/\varpi^n\in V$ holds for all $n$, that is, $a\in \bigcap_n(\varpi^n)$.
But $\bigcap_n(\varpi^n)=0$ (\cite{FK18}*{Ch.~0, Prop.~6.7.2}), so $a=0$, a contradiction.

\ref{vt-i}$\Rightarrow$\ref{vt-iii}. As above, there is a topologically nilpotent unit $\varpi\in V$ of $K$.
Take $V$ as an open subring of $K$, it suffices to show that $\{(\varpi^n)\}_{n=1}^{+\infty}$ form a basis of open neighborhoods of $0\in V$. 
We have known that every $U_{\gamma}$ contains some $(\varpi^n)$. 
Conversely, for a fixed $n\in \b{Z}_+$, there is $\gamma\in \GG$ such that $U_{\gamma}\subset (\varpi^n)$.
To see this, we need to find $\gamma\in \GG$ such that the condition $\nu(x)>\gamma$ implies that $\nu(x)>\nu(\varpi^n)$.
It suffices to let $\gamma>\nu(\varpi^n)=n\nu(\varpi)$, say, $\gamma=(n+1)\nu(\varpi)$.

\ref{vt-iii}$\Rightarrow$\ref{vt-iv}. By definition of Tate rings, it is obvious.

\ref{vt-i}$\Rightarrow$\ref{vt-ii}.
The argument for \ref{vt-i}$\Rightarrow$\ref{vt-iii} implies that $\{(\varpi^n)\}_n$ form a basis of open neighborhoods of $0\in V$.
As $\varpi$ lies in the height-one prime ideal, the valuation topology on $K$ is induced by its rank-one valuation.

\ref{vt-ii}$\Rightarrow$\ref{vt-i}.
The rank-one valuation corresponds to the height-one prime ideal of $V$, since all nonequivalent valuations of $K$ are in one-to-one correspondence with the prime ideals of $V$ (\cite{FK18}*{Ch.~0, Prop.~6.2.9}).

\ref{vt-iv}$\Rightarrow$\ref{vt-i}. 
For a topologically nilpotent unit $\varpi\in K$, we prove that $\fp\ce \sqrt{(\varpi)}$ is the prime ideal of height one. 
If $a,b\in V$ such that $ab\in \fp$ and $b\not\in \fp$, then there are an integer $n>0$ and $c\in V$ such that $a^nb^n=\varpi c$, and $\varpi/b^m\in V$ holds for every integer $m>0$. 
It follows that $a^{2n}=\varpi(\varpi/b^{2n})c^2\in (\varpi)$, so $a\in \fp$ and we see that $\fp$ is a prime.
To see that $\fp$ is of height one, note that the set of ideals of $V$ is totally ordered under inclusion and $\varpi^n$ tends to zero, every nonzero prime ideal $\mathfrak{q}$ between $(0)$ and $\fp$ satisfies $(\varpi^N)\subset \mathfrak{q}\subset \fp$ for some $N$. 
Taking radicals of these inclusions, we find that $\mathfrak{q}=\fp$, thus $\fp$ is of height one.
\epf
\eprop

\bpp[Nonarchimedean fields]\label{nonarch} 
A \emph{nonarchimedean field} is a topological field $K$ whose topology is induced by a nontrivial valuation of rank one on $K$.\footnote{Some authors additionally require the completeness of $K$, for instance, Scholze \cite{Sch12}*{Def.~2.1}.} 
By the end of \ref{abv}, a topological field $K$ is nonarchimedean if and only if its topology is induced by a nonarchimedean absolute value on $K$.
If an absolute value on $K$ is not nonarchimedean, then it is archimedean.
We note that the existence of absolute values on the topological field $K$ is a prerequisite for our discussion of Archimedean properties.
\epp

\bpp[$a$-adic topologies]\label{a-adic-top} 
For a valuation ring $V$ and an element $a\in \fm_V\backslash \{0\}$, the $a$-adic topologies on $V$ and on $V[\f{1}{a}]$ are determined by the respective fundamental systems of open neighborhoods of $0$:
\[
   \textstyle \x{$\{a^nV\}_{n\geq 0}$\q and \q $\{\mathrm{Im}(a^nV\ra V[\f{1}{a}])\}_{n\geq 0}$.}
\]
Note that the $a$-adic topology on $V[\f{1}{a}]$ is not defined by ideals, since such topology is only $V$-linear (\cite{GR18}*{Def.~8.3.8(iii)}). 
Then, the $a$-\emph{adic completion}s $\hva$ and $\vaih$ are the following inverse limits:
\[
 \textstyle  \x{$\hva\ce \varprojlim_{n>0}V/a^n$\q and \q $\vaih\ce \varprojlim_{n>0}(V[\f{1}{a}]/\mathrm{Im}(a^nV\ra V[\f{1}{a}]))$.}
\]
\epp
\bprop\label{a-completion}
For a valuation ring $V$ and a nonzero element $a\in \fm_V$, 
\benumr
\item\label{a-prime} $\sqrt{(a)}$ is the minimal one among all the prime ideals containing $(a)$, while $\bigcap_{n>0}(a^n)$ is the maximal one among all the prime ideals contained in $(a)$;
\item\label{comp-sep} the $a$-adic completion $V\ra \hva$ factors through the $a$-adically separated quotient $V/\bigcap_{n>0}(a^n)$;
\item\label{loc-com-val} the rings $V[\f{1}{a}]$ and $\hva$ are valuation rings, and we have $\vaih=\hvai$;
\item\label{rank-val} if $V$ has finite rank $n\geq 1$ and $(a)$ is between the primes of heights $r-1$ and $r$ for $1\leq r\leq n$, then $\mathrm{rank} (\hva)=n-r+1$ and $\mathrm{rank} (V[\f{1}{a}])=r-1$; 
\item\label{frac-res} we have $\hvai=\Frac \hva$, which is also the $a$-adic completion of the residue field of $V[\f{1}{a}]$;
\item\label{comp-res} the valuative completion $\wh{V}$, the $a$-adic completion $\hva$, and $V$ share the same residue field;
\item\label{fp-a-comp} we have an isomorphism to a fiber product of rings $V\isoto \via\times_{\nhka}\hva$, where $\hka$ is the $a$-adic completion of $K=\Frac V$.
\eenum
\eprop
\bpf
For \ref{a-prime}, see \cite{FK18}*{Ch.~0, Prop.~6.2.3 and 6.7.1}.
For \ref{comp-sep}, see the end of \cite{FK18}*{Ch.~0, Cor.~9.1.5}.
For \ref{loc-com-val}, by \cite{FK18}*{Ch.~0, Cor.~9.1.5}, $\hva$ is a valuation ring.
Let $\f{\alpha}{\beta}\in K\ce\Frac V$ be an element which is not in $V[\f{1}{a}]$.
Hence, $a^n\f{\alpha}{\beta}\not\in V$ for every $n>0$, which means that $\f{\beta}{\alpha}\in (a^n)$ for every $n>0$.
So $\f{\beta}{\alpha}$ lies in  $\bigcap_{n>0}(a^n)$, the maximal ideal of $V[\f{1}{a}]$ by \ref{a-prime}.
The relation $\hvai=\vaih$ is due to \cite{BC22}*{Ex.~2.1.10~(2)} and the fact that $V$ is $a$-torsion-free.
For \ref{rank-val}, by \ref{a-prime}, the rank of $V[\f{1}{a}]$ is $r-1$; also, $\fq\ce \bigcap_{n>0}(a^n)$ is the prime ideal of height $r-1$. 
Note that $\hva$ is the $a$-adic completion of the $a$-adically separated quotient $V/\fq$, whose rank is $n-r+1$. 
By \cite{FK18}*{Ch.~0, Thm.~9.1.1~(5)}, we conclude that $\hva$ is also of rank $n-r+1$.  
For \ref{frac-res}, by \cite{FK18}*{Ch.~0, Prop.~6.7.2}, $\hvai$ is the fraction field of $\hva$.
By \ref{a-prime}, the residue field $\kappa$ of $V[\f{1}{a}]$ is $V[\f{1}{a}]/\bigcap_{n>0}a^nV$, hence the $a$-adic completion of $\kappa$ is $\vaih$, which is $\hvai$ by \ref{loc-com-val}.
For \ref{comp-res}, see \cite{EP05}*{Prop.~2.4.4}, \ref{comp-sep} and \cite{FK18}*{Ch.~0, Thm.~9.1.1~(2)}.
For \ref{fp-a-comp}, we apply \SP{0BNR} to the $a$-adic completion $V\ra \hva$: note that $V/a^nV\simeq \hva/a^n\hva$ for every positive integer $n$ (\cite{FK18}*{Ch.~0, 7.2.8}), also, $V[a^{\infty}]=\ker (V\ra \via)=0$ and $\hva[a^{\infty}]=\ker(\hva\ra \hvai)=0$; the exactness of $0\ra V\ra \via\oplus \hva\ra \hvai\ra 0$ implies the desired isomorphism $V\isoto \via\times_{\nhka}\hva$.
\epf
\bpp[Comparison of topologies]\label{comparison}
We have compared different valuation topologies to some extent (\Cref{val-top}). 
Now, consider three kinds of topologies on a valuation ring $V$: the $a$-adic topology, the valuation topology, and the $\fm_V$-adic topology, where $\fm_V\subset V$ is the maximal ideal.
First, the $\fm_V$-adic topology is usually non-Hausdorff and does not coincide with any $a$-adic topology: for the rank-one valuative completion $\b{C}_p$ of the algebraic closure $\ov{\b{Q}_p}$ of $\b{Q}_p$, the maximal ideal $\fm$ of the valuation ring $\sO_{\b{C}_p}$ of $\b{C}_p$ satisfies $\fm=\fm^2$.
Thus, for every nonzero $a\in \fm$ and every $n>0$, we have $(a)\not\supset\fm^n=\fm$.
Secondly, for $a,b\in \fm_V\backslash\{0\}$, the comparison of $a$-adic and $b$-adic topologies is \cite{FK18}*{Ch.~0, Prop.~7.2.1}:
\[
   \x{\!\!\!\!\emph{the $a$-adic and $b$-adic topologies coincide} \q $\Leftrightarrow$ \q  $\sqrt{(a)}=\sqrt{(b)}$,}
\]
and in such case, the $a$-adic completion is equal to the $b$-adic completion; also, the Henselizations of pairs $(V,a)$ and $(V,b)$ coincide (\SP{0F0L}). 
Thirdly, to compare $a$-adic topologies and valuation topologies, by \Cref{val-top},  $V$ has a prime ideal of height one $\fp$ if and only if there is a topologically nilpotent $\varpi\in V\backslash \{0\}$ such that the valuation topology on $V$ is $\varpi$-adic and $\sqrt{(\varpi)}=\fp$.
In conclusion, 
\[
   \x{\emph{the valuation topology is nonarchimedean $\Leftrightarrow$ it is $a$-adic for an $a\in \fm_V$ such that $\sqrt{(a)}$ is height-one}.}
\]
Of course, valuation topologies and $a$-adic topologies do not coincide in general since each kind of both has aforementioned internal differences.
Lastly, a valuation ring $V$ equipped with an $a$-adic topology for some $a\in \fm_V\backslash\{0\}$ may not have any prime ideal of height one, so its valuation topology can not be $a$-adic.
\epp
\bcor\label{completion-NA}
For a valuation ring $V$, an element $a\in \fm_V\backslash\{0\}$, and the $a$-adic completion $\hva$ of $V$, the fraction field $\hka\ce \Frac\hva$ is a nonarchimedean field with respect to the $a$-adic topology.
\ecor
\bpf
Let $\GG$ be the value group of $\hka$.
If there is a $\gamma\in \GG$ such that $\nu(a^n)\leq \gamma$ for all $n\in \b{Z}_+$, then there is a $b\in \hva$ such that $\nu(b)=\gamma$ and $b\in \bigcap_n(a^n)$.
Since $\hva$ is $a$-adically separated, we have $\bigcap_n(a^n)=0$ so $b=0$, that is, $\gamma=\infty\not\in \GG$. 
Thus every $U_{\gamma}$ contains some $a^n$, that is, $a$ is topologically nilpotent for the valuation topology, hence $\hka$ is a Tate ring with its open subring $\hva$.
By \Cref{val-top}, $\sqrt{(a)}$ is of height one in $\hva$, the valuation topology on $\hka$ is $a$-adic hence nonarchimedean by \ref{comparison}.
\epf
We end this appendix by a comparison of Henselianity and completeness of valuation rings.
\bprop\label{hens-val-field} 
For a valuation ring $V$ equipped with an $a$-adic topology for an element $a\in \fm_V\backslash \{0\}$.
If $V$ is $a$-adically complete, then the pair $(V,a)$ is Henselian.
If $V$ has finite rank $n$ and $a$ is not in the unique prime $\fp\subset V$ that is of height $n-1$, then the $a$-adic completion $\hva$ is a Henselian local ring.
\eprop
\bpf
If $V$ is $a$-adically complete, then the Henselianity of $(V,a)$ follows from \cite{FK18}*{Ch.~0, Prop.~7.3.5~(1)}. 
Now we show the second part. 
By \Cref{a-completion}~\ref{rank-val}, $\hva$ is of rank one. 
Since $(\hva,a\hva)$ is a Henselian pair and \Cref{a-completion}~\ref{a-prime} implies that $\sqrt{(a)}=\fm_V$, by \SP{0F0L}, the local ring $\hva$ is Henselian.  
\epf


\end{appendix}

\end{document}